\documentclass[12pt,a4paper,notitlepage]{article}
\usepackage[left=3cm,right=3cm,top=1cm,bottom=2cm]{geometry} 
\usepackage{amsmath,amssymb,amsthm} 
\usepackage{hyperref}

\usepackage{pgfplots}
\pgfplotsset{compat=1.15}
\usepackage{mathrsfs}
\usetikzlibrary{arrows}

\theoremstyle{definition}
\newtheorem{Exmp}{Example}[section]
\newtheorem*{Exmp*}{Example}
\newtheorem{Rem}{Remark}[section]
\newtheorem*{Rem*}{Remark}
\newtheorem{Lem}{Lemma}[section]
\newtheorem*{Lem*}{Lemma}
\newtheorem{Thm}{Theorem}
\newtheorem{Def}{Definition}[section]
\newtheorem*{Def*}{Definition}
\newtheorem{Prop}{Proposition}[section]
\newtheorem*{Prop*}{Proposition}

\DeclareFontEncoding{LS1}{}{}
\DeclareFontSubstitution{LS1}{stix}{m}{n}
\DeclareSymbolFont{symbols3}{LS1}{stixbb}{m}{n}
\SetSymbolFont{symbols3}{bold}{LS1}{stixbb}{b}{n}
\DeclareMathSymbol{\Wedge}{\mathbin}{symbols3}{"A3}

\title{The local Schwarz paradox and the gradient of strongly differentiable functions}
\author{Paolo Roselli}
\date{\today}

\begin{document}
\maketitle
    \begin{abstract}
We show that the gradient of a multi-variable strongly differentiable function at a point is the limit of a single coordinate-free Clifford quotient between a multi-difference pseudo-vector and a pseudo-scalar,
or of a sum of Clifford quotients between scalars (as numerators) and vectors (as denominators), both built on a same non degenerate simplex contracting to that point.
Then, we provide some conjectures implied by the foregoing result.
    \end{abstract}

\section{Introduction}
The derivative $f'(x_0)$ of a scalar function $f:\Omega\subseteq\mathbb{R}\to \mathbb{R}$ at a point $x_0$ internal to $\Omega$ is the limit of the quotient $\displaystyle \frac{\Delta f_{(a,x_0)}}{\Delta_{(a,x_0)}}$, between the differences $\Delta f_{(a,x_0)}=f(a)-f(x_0)$, and $\Delta_{(a,x_0)}=a-x_0$
\[
\lim_{a\to x_0}
\frac{\Delta f_{(a,x_0)}}{\Delta_{(a,x_0)}}
=
f'(x_0)\ .
\]
The strong derivative $f^*(x_0)$ is defined\footnote{See for example~\cite{Peano1892}, \cite{Leach}, \cite{Esser}.} by a stronger, but fully symmetric, limit 
\arraycolsep=1.pt\def\arraystretch{0.5}
\[
\lim_{
\begin{array}{c}
	\scriptstyle (a,b)\to (x_0,x_0)\\
	\scriptstyle a, b \textrm{ distinct}
\end{array}
}
\frac{\Delta f_{(a,b)}}{\Delta_{(a,b)}}
=
f^*(x_0) 
\]
of the difference quotient $\displaystyle \frac{f(a)-f(b)}{a-b}$.
So, $f^*(x_0)=f'(x_0)$ when the strong derivative exists.
The existence of those limits corresponds to a well known geometric phenomenon: the secant line of the graph of $f$ at points $\big(a,f(a)\big)$, $\big(b,f(b)\big)$
\[
y= f(a) + \frac{\Delta f_{(a,b)}}{\Delta_{(a,b)}} (x-a)
\]
assumes, as the non degenerate segment joining $a$ and $b$ contracts to $x_0$, the limit position
\[
y= f(x_0) + f'(x_0) (x-x_0)
\]
which is the line tangent the graph of $f$ at point $\big(x_0,f(x_0)\big)$.\\
\subsection{The local Schwarz paradox}
By analogy, one would expect that planes secant the graph of a two-variable real function at three non collinear points always assume as limit position, as those three points converge on the graph 
to a same limit point, the position of the plane tangent the graph at that limit point. Amazingly, this is not the case, even for smooth functions.
\begin{Exmp}\label{exmp:Schwarz}
Let $f(x,y)=\sqrt{1-x^2}$. The plane secant the graph $\Big\{\big(x,y,f(x,y)\big)\ :\ (x,y)\in[-1,1]\times\mathbb{R}\Big\}$ at points $(0,0,1)$, $(-\alpha,\beta,\sqrt{1-\alpha^2})$, and $(\alpha,\beta,\sqrt{1-\alpha^2})$, is defined by the relation between $(x,y,z)\in\mathbb{R}^3$
\begin{equation}
\label{eq:Schwarz secant}
z=1
-\frac{\alpha^2}{\beta(1+\sqrt{1-\alpha^2})}y
\end{equation}

\begin{itemize}
	\item If $\beta=\alpha$ and $\alpha\to 0$, then the limit position of the secant plane~(\ref{eq:Schwarz secant}) is the plane $z=1$, which is tangent the graph of $f$ at $(0,0,1)$;
	\item If $\beta=\alpha^2$ and $\alpha\to 0$, then the limit position of the secant plane~(\ref{eq:Schwarz secant}) is the plane $z=1-\frac{1}{2}y$, which is not tangent to the graph of $f$ at point $(0,0,1)$;
		\item If $\beta=\alpha^3$ and $\alpha\to 0$, then the limit position of the secant plane~(\ref{eq:Schwarz secant}) is the plane $y=0$, which is even orthogonal to the tangent plane!
\end{itemize}
Such local divergent phenomenon implies a global divergent one, concerning the area of smooth surfaces: the so called ``Schwarz paradox''(see for example~\cite{Peano1890}, \cite{Zames}, \cite{Gandon}, and \cite{Roselli}).
\end{Exmp}
\subsection{Some notations for the following}
In this work we perform coordinate-free vector computations. In order to better distinguish dimensionless numbers (i.e., scalars) from vectors, we adopt the following notations:
\begin{itemize}
	\item vectors are denoted by \textbf{bold} Latin lower case letters;
	\item real numbers are denoted by non-bold Latin or Greek lower case letters;
	\item $\mathbb{E}_n$ denotes an $n$-dimensional Euclidean space: a real vector space with a positive definite symmetric bilinear form; such bilinear form is denoted by~$\mathbf{u}\cdot\mathbf{v}$, for each  
$\mathbf{u},\mathbf{v}\in\mathbb{E}_n$.
\end{itemize}
In particular, if $\{\mathbf{e}_1,\mathbf{e}_2\}$ is an orthonormal basis for $\mathbb{E}_2$ (i.e., $\mathbf{e}_1\cdot\mathbf{e}_2=0$, and $\mathbf{e}_1\cdot\mathbf{e}_1=\mathbf{e}_2\cdot\mathbf{e}_2=1$), $\mathbf{x}=x\mathbf{e}_1+y\mathbf{e}_2$, $\mathbf{a}=\mathbf{0}$, $\mathbf{b}=-\alpha\mathbf{e}_1+\beta\mathbf{e}_2$, $\mathbf{c}=\alpha\mathbf{e}_1+\beta\mathbf{e}_2$, then the foregoing Example~\ref{exmp:Schwarz} can be resumed as follows: 
the secant plane $\displaystyle z= f(\mathbf{a})-\frac{\alpha^2}{\beta(1+\sqrt{1-\alpha^2})}y$ has no limit position when the non degenerate triangle (having vertices $\mathbf{a}$, $\mathbf{b}$ and $\mathbf{c})$ contracts to the point $\mathbf{0}$, because $\displaystyle \lim_{(\alpha,\beta)\to(0,0)}\frac{\alpha^2}{\beta(1+\sqrt{1-\alpha^2})}$ does not exist.

In general, for a non linear two variable function $f$, the local Schwarz paradox corresponds to the non existence of the strong limit
\arraycolsep=1.pt\def\arraystretch{0.5}
\[
\lim_{
\begin{array}{c}
	\scriptstyle (\mathbf{a},\mathbf{b},\mathbf{c})\to (\mathbf{x_0},\mathbf{x_0},\mathbf{x_0})\\
	\scriptstyle \mathbf{a},\mathbf{b},\mathbf{c} \textrm{ not collinear}
\end{array}
}
\mathbf{r}_{f_{(\mathbf{a},\mathbf{b},\mathbf{c})}}\ ,
\]
where $\mathbf{r}_{f_{(\mathbf{a},\mathbf{b},\mathbf{c})}}$ is the vector in $\mathbb{E}_2$ such that the secant plane to the graph of $f$ at points 
$\big(\mathbf{a}, f(\mathbf{a})\big)$, $\big(\mathbf{b}, f(\mathbf{b})\big)$, and $\big(\mathbf{c}, f(\mathbf{c})\big)$ has Cartesian equation in~$\mathbb{E}_2\oplus\mathbb{R}$
\[
z= f(\mathbf{a})\ +\ 
\mathbf{r}_{f_{(\mathbf{a},\mathbf{b},\mathbf{c})}} \cdot (\mathbf{x}-\mathbf{a})\ ,
\]
when $\mathbf{a}$, $\mathbf{b}$, and $\mathbf{c}$ are vertices of a non degenerate triangle internal to the domain $\Omega$ of $f$.
\begin{Rem*}
The letter ``r'' in $\mathbf{r}_{f_{(\mathbf{a},\mathbf{b},\mathbf{c})}}$ stands for ``ratio'', because we will show that such vector can be considered as a ratio of a pseudo-vector and a pseudo-scalar in a suitable Clifford algebra.
\end{Rem*}
\begin{Rem*}
Vectors of $\mathbb{E}_n$ will also be called ``points'', because we identify the Euclidean vector space $\mathbb{E}_n$ with an Euclidean affine space $\mathcal{E}_n$ (modeled on $\mathbb{E}_n$) where an arbitrary point $O\in\mathcal{E}_n$ is considered as reference point and identified with the zero vector $\mathbf{0}\in\mathbb{E}_n$. This allows us to interpret geometrically some subsets of $\mathbb{E}_n$. For instance, 
\begin{itemize}
	\item a ``line'' in $\mathbb{E}_n$ is the set 
\begin{center}
$
\mathcal{L}_{(\mathbf{a},\mathbf{b})}
=
\{
\alpha\mathbf{a} +\beta \mathbf{b}\ : \ \alpha,\beta\in\mathbb{R},\ \alpha+\beta=1
\}
$
\end{center}
for some distinct points $\mathbf{a}$, $\mathbf{b}$ in $\mathbb{E}_n$;
\item three distinct points $\mathbf{a}$, $\mathbf{b}$ and $\mathbf{c}$ in $\mathbb{E}_n$ are ``collinear'' if they all belong to a same line;
	\item a ``plane'' in $\mathbb{E}_n$ is the set 
\begin{center}
$
\mathcal{L}_{(\mathbf{a},\mathbf{b},\mathbf{c})}
=
\{
\alpha\mathbf{a} +\beta \mathbf{b}+\gamma \mathbf{c}\ : \ \alpha,\beta,\gamma\in\mathbb{R},\ \alpha+\beta+\gamma=1
\}
$
\end{center}
for some distinct and non collinear points $\mathbf{a}$, $\mathbf{b}$, and $\mathbf{c}$ in $\mathbb{E}_n$;
\item four distinct points $\mathbf{a}_1$, $\mathbf{a}_2$, $\mathbf{a}_3$, and $\mathbf{a}_4$ in $\mathbb{E}_n$ are ``coplanar'' if they all belong to a same plane.
\end{itemize}
\end{Rem*}

\subsection{Summary of this work}
In this work we first provide an explicit coordinate-free expression of vector~$\mathbf{r}_{f_{(\mathbf{a},\mathbf{b},\mathbf{c})}}\in\mathbb{E}_2$, both 
\begin{itemize}
	\item as a linear combination of vectors $\mathbf{a}$, $\mathbf{b}$, $\mathbf{c}$ (in Section~\ref{sec:r as linear combination}), and
	\item as a non-commutative quotient
	\[
	\Big(\mathbf{\Delta}f_{(\mathbf{a},\mathbf{b},\mathbf{c})}\Big)
	\big(\Delta_{(\mathbf{a},\mathbf{b},\mathbf{c})}\big)^{-1}
	\]
	of a multi-difference vector $\mathbf{\Delta} f_{(\mathbf{a},\mathbf{b},\mathbf{c})}\in\mathbb{E}_2$, and the oriented area
	\[
	\Delta_{(\mathbf{a},\mathbf{b},\mathbf{c})}
	= \mathbf{a}\wedge \mathbf{b}+\mathbf{b}\wedge\mathbf{c}+\mathbf{c}\wedge\mathbf{a}
	=
	(\mathbf{a}-\mathbf{b})\wedge (\mathbf{b}-\mathbf{c})
	\]
	with respect to the geometric product of the Clifford algebra
	\[
	\mathbb{G}_2 = \mathcal{C}\ell(\mathbb{E}_2)=\mathcal{C}\ell(2,0,0)
	\]
	generated by the two-dimensional Euclidean space $\mathbb{E}_2$ (in Section~\ref{sec:r as quotient n=2}).
\end{itemize}
Then, in Section~\ref{sec:mean multi diff vector} we introduce a new multi-difference vector $\overline{\mathbf{\Delta}} f_{(\mathbf{a},\mathbf{b},\mathbf{c})}$, such that the corresponding plane
\[
z= f(\mathbf{a})\ +\  
\left[\Big(\overline{\mathbf{\Delta}} f_{(\mathbf{a},\mathbf{b},\mathbf{c})}\Big)
	\big(\Delta_{(\mathbf{a},\mathbf{b},\mathbf{c})}\big)^{-1}\right]\cdot (\mathbf{x}-\mathbf{a})
\]
(called ``mean secant plane'') always assume, as limit position, that of the tangent plane
\[
z= f(\mathbf{x_0})\ +\  
\nabla f(x_0) \cdot (\mathbf{x}-\mathbf{x_0})\ ,
\]
because we prove that the limit
\arraycolsep=1.pt\def\arraystretch{0.5}
\[
\lim_{
\begin{array}{c}
	\scriptstyle (\mathbf{a},\mathbf{b},\mathbf{c})\to (\mathbf{x_0},\mathbf{x_0},\mathbf{x_0})\\
	\scriptstyle \mathbf{a},\mathbf{b},\mathbf{c} \textrm{ not collinear}
\end{array}
}
\Big(\overline{\mathbf{\Delta}} f_{(\mathbf{a},\mathbf{b},\mathbf{c})}\Big)
	\big(\Delta_{(\mathbf{a},\mathbf{b},\mathbf{c})}\big)^{-1}
\ ,
\]
always exists when $f$ is strongly differentiable at~$\mathbf{x}_0$, and it is equal to the gradient~$\nabla f(\mathbf{x}_0)$ (see Theorem~\ref{thm:n=2}).
Moreover, we show in Proposition~\ref{prop:mean multi-difference quotient n=2} that the Clifford ratio 
$
\displaystyle
\Big(\overline{\mathbf{\Delta}} f_{(\mathbf{a},\mathbf{b},\mathbf{c})}\Big)
\big(\Delta_{(\mathbf{a},\mathbf{b},\mathbf{c})}\big)^{-1}
$
can also be written as a sum of quotients between numbers and vectors, strongly resembling the scalar difference quotients.
Then, we extend the foregoing results to three and to arbitrary higher dimensions\footnote{See, for instance, Theorems~\ref{thm:n=3} and~\ref{thm:n}.}.
Finally, in Section~\ref{sec:conclusions} we sketch some possible consequences and conjectures implied by the foregoing results.
\section{Some remainders of Geometric Algebra}
As announced, in this work we use the associative vector algebra
\[
\mathbb{G}_n = \mathcal{C}\ell(\mathbb{E}_n)=\mathcal{C}\ell(n,0,0)\ ,
\]
which is the Clifford geometric algebra\footnote{See for example~\cite{Macdonald}, \cite{HestenesSob}, \cite{Dorst}, \cite{Guerlebeck-Sproessig}, \cite{Choquet-Bruhat-DeWitt-Morette}, \cite{Vaz-Rocha}.} generated by the $n$-dimensional Euclidean space $\mathbb{E}_n$ (the geometric product being denoted by juxtaposition). In particular, in $\mathbb{G}_n$, we have that
\begin{itemize}
	\item for all $\mathbf{v}\in\mathbb{E}_n$ $\mathbf{v}\mathbf{v}=\mathbf{v}^2=\mathbf{v}\cdot \mathbf{v}$, which implies that $\displaystyle \frac{1}{2}(\mathbf{u}\mathbf{v}+\mathbf{v}\mathbf{u})=\mathbf{u}\cdot\mathbf{v}$, for all $\mathbf{u},\mathbf{v}\in\mathbb{E}_n$;
	\item scalars always commute with the geometric product;
	\item if $\{\mathbf{e}_1,\dots,\mathbf{e}_n\}$ is an orthonormal basis for $\mathbb{E}_n$ (i.e., $\mathbf{e}_i\cdot\mathbf{e}_j=0$ when $i\ne j$, and~$\mathbf{e}_i\cdot\mathbf{e}_i=1$), then
	\[
	\{1\}\cup\{\mathbf{e}_{i_1}\cdots\mathbf{e}_{i_k}\}_{1\le i_1<\cdots< i_k\le n}
	\ \ \ \ 
	\textrm{ with } 1\le k\le n
	\]
	is a basis for $\mathbb{G}_n$. So, the associative vector algebra $\mathbb{G}_n$ has dimension $2^n$.  
\end{itemize}
Thus, for instance, 
	\[
	\{1,\ \mathbf{e}_1\ ,\ \mathbf{e}_2\ ,\ \mathbf{e}_1\mathbf{e}_2\}
	\ \ \ \ 
	\textrm{ is a basis for } \mathbb{G}_2\ ;
	\]
	\[
	\{1\ ,\ \mathbf{e}_1\ ,\ \mathbf{e}_2\ ,\ \mathbf{e}_3\ ,\ \mathbf{e}_1\mathbf{e}_2\ ,\ \mathbf{e}_1\mathbf{e}_3\ ,\ \mathbf{e}_2\mathbf{e}_3\ ,\ \mathbf{e}_1\mathbf{e}_2\mathbf{e}_3\}
	\ \ \ \ 
	\textrm{ is a basis for } \mathbb{G}_3\ .
	\]
Defining $\displaystyle \mathbf{u}\wedge\mathbf{v}=\frac{1}{2}(\mathbf{u}\mathbf{v}-\mathbf{v}\mathbf{u})$ for each $\mathbf{u},\mathbf{v}\in\mathbb{E}_n$, you can verify that	$\mathbf{u}\wedge\mathbf{v}=-\mathbf{v}\wedge\mathbf{u}$, $\mathbf{u}\wedge\mathbf{u}=0$, and $\mathbf{u}\wedge\mathbf{v} =\mathbf{u}\mathbf{v}$  when vectors $\mathbf{u}$ and $\mathbf{v}$ are mutually orthogonal (i.e., $\mathbf{u}\cdot\mathbf{v}=0$).\\
\begin{Rem*}
Note that, for each $\mathbf{u},\mathbf{v}\in\mathbb{E}_n$, we have that $ \mathbf{u}\mathbf{v}=(\mathbf{u}\cdot\mathbf{v})+(\mathbf{u}\wedge\mathbf{v})$. In what follows, to limit the use of parentheses, we adopt a rule of precedence between operations: geometric product is used first, secondly the scalar product ``$\cdot$'', then ``$\wedge$'' is performed, and finally the sum ``$+$''.
\end{Rem*}
\subsection{The $2\times 2$ determinant as a Clifford quotient and as a scalar product}
Clifford geometric algebra allows coordinate-free vector computations having interesting geometric interpretations. For example, we can give a coordinate-free interpretation to the determinant of a $2\times 2$ real matrix
\arraycolsep=2.5pt\def\arraystretch{1.}
\[
\left( 
\begin{array}{cc}
	\mu_1 & \mu_2 \\
	\nu_1 & \nu_2
\end{array}
\right)
\]
as a Clifford ratio. As a matter of fact, let us fix any ordered couple $\mathbf{e}_1$, $\mathbf{e}_2$ of orthonormal vectors in $\mathbb{E}_n$ (with $n\ge 2$), and let us consider $\mathbf{u}=\mu_1 \mathbf{e}_1+\mu_2\mathbf{e}_2$ and $\mathbf{v}=\nu_1\mathbf{e}_1+\nu_2\mathbf{e}_2\in\mathbb{E}_2\subseteq \mathbb{E}_n$, then
\begin{align*}
(\mathbf{u}\wedge\mathbf{v})(I_2)^{-1}
 & =
\big[(\mu_1 \mathbf{e}_1+\mu_2\mathbf{e}_2)\wedge(\nu_1\mathbf{e}_1+\nu_2\mathbf{e}_2)\big](\mathbf{e}_1\mathbf{e}_2)^{-1}\\
& =
\det\left( 
\begin{array}{cc}
	\mu_1 & \mu_2 \\
	\nu_1 & \nu_2
\end{array}
\right)
\mathbf{e}_1\mathbf{e}_2(\mathbf{e}_2\mathbf{e}_1)
=
\det\left( 
\begin{array}{cc}
	\mu_1 & \mu_2 \\
	\nu_1 & \nu_2
\end{array}
\right)
\end{align*} 
as $I_2=\mathbf{e}_1\wedge\mathbf{e}_2=\mathbf{e}_1\mathbf{e}_2$, and $(I_2)^{-1}=\mathbf{e}_2\mathbf{e}_1=-\mathbf{e}_1\mathbf{e}_2=-I_2$.
Note that, as $(\mathbf{u}\wedge\mathbf{v})(I_2)^{-1}=(I_2)^{-1}
(\mathbf{u}\wedge\mathbf{v})$, one could also write 
$
\displaystyle
\det\left( 
\begin{array}{cc}
	\mu_1 & \mu_2 \\
	\nu_1 & \nu_2
\end{array}
\right)
=
\frac{\mathbf{u}\wedge\mathbf{v}}{I_2}$.
\begin{Rem*}
The geometric product $I_2$ does not depend on the particular orthonormal basis $\{\mathbf{e}_1,\mathbf{e}_2\}$ of $span \{\mathbf{e}_1,\mathbf{e}_2\}$, but only on its orientation. More precisely, if $\{\mathbf{g}_1,\mathbf{g}_2\}$ is any other orthonormal basis of $span \{\mathbf{e}_1,\mathbf{e}_2\}$, then $\mathbf{g}_1\mathbf{g}_2\ (=\mathbf{g}_1\wedge\mathbf{g}_2)$ is equal to $I_2$ or $-I_2$. That is why $I_2$ is called an ``orientation'' of $span \{\mathbf{e}_1,\mathbf{e}_2\}$. An analogue property and definition is valid for the product $\mathbf{e}_1\cdots\mathbf{e}_k=I_k$ of any ordered list of mutually orthonormal vectors in $\mathbb{E}_n$ (with $n\ge k$).
\end{Rem*}
\begin{Rem*}
We recall  that $span S$ means the smallest linear subspace that contains the set $S\subseteq\mathbb{E}_n$.
\end{Rem*}
Thus, a $2\times 2$ determinant can be considered as the Clifford  ratio between the two blades $\mathbf{u}\wedge\mathbf{v}$ and $I_2$ (a ``blade'' being the geometric product of non zero mutually orthogonal vectors). Those elements are also called ``$\mathbb{G}_2$-pseudo-scalars'', because they are scalar multiples of the orientation $I_2$ of $\mathbb{E}_2$ (which generates $\mathbb{G}_2$), and can be interpreted as oriented areas in $span\{\mathbf{e}_1,\mathbf{e}_2\}=\mathbb{E}_2$. Later we will see that this Clifford geometric interpretation of a $2\times 2$ determinant holds for every $k\times k$ determinant. Let us also observe that, if $\mathbf{v}=\nu_1\mathbf{e}_1+\nu_2\mathbf{e}_2$, then 
\begin{align*}
\mathbf{v} I_2 
& 
= 
(\nu_1\mathbf{e}_1+\nu_2\mathbf{e}_2)\mathbf{e}_1\mathbf{e}_2
=
\nu_1\mathbf{e}_1\mathbf{e}_1\mathbf{e}_2+
\nu_2\mathbf{e}_2\mathbf{e}_1\mathbf{e}_2
=
-\nu_1\mathbf{e}_1\mathbf{e}_2\mathbf{e}_1
-\nu_2\mathbf{e}_1\mathbf{e}_2\mathbf{e}_2\\
& =
-\mathbf{e}_1\mathbf{e}_2(\nu_1\mathbf{e}_1+\mathbf{e}_2) = -I_2\mathbf{v}=-\nu_2 \mathbf{e}_1+\nu_1\mathbf{e}_2\in span \{\mathbf{e}_1,\mathbf{e}_2\}\ .
\end{align*} 
Thus, we can also write a $2\times 2$ determinant as a scalar product
\begin{align*}
\det\left( 
\begin{array}{cc}
	\mu_1 & \mu_2 \\
	\nu_1 & \nu_2
\end{array}
\right)
& =
(\mathbf{u}\wedge\mathbf{v})(I_2)^{-1}
=
-\frac{1}{2}(\mathbf{u}\mathbf{v}-\mathbf{v}\mathbf{u})I_2
=
\frac{1}{2}\big(-\mathbf{u}\mathbf{v}I_2+\mathbf{v}\mathbf{u}I_2\big)\\
& =
\frac{1}{2}\big(\mathbf{u}I_2\mathbf{v}+\mathbf{v}\mathbf{u}I_2\big)
=
\big(\mathbf{u}I_2\big)\cdot\mathbf{v}\ ,
\end{align*}
where $\mathbf{u}I_2$ is the vector obtained by rotating vector $\mathbf{u}$ of a right angle counterclockwise in $span\{\mathbf{e}_1,\mathbf{e}_2\}$ \Big(provided $\mathbf{e}_1$ and $\mathbf{e}_2$ are mutually oriented like this  
\begin{tikzpicture}[line cap=round,line join=round,>=triangle 45,x=0.7cm,y=0.7cm]
\clip(-0.5,-0.1) rectangle (1.1,1.1);
\draw [->,line width=1.pt,-stealth] (0.,0.) -- (0.,1.);
\draw [->,line width=1.pt,-stealth] (0.,0.) -- (1.,0.);
\begin{scriptsize}
\draw[color=black] (0.8817931144633162,0.2) node {$\mathbf{e}_1$};
\draw[color=black] (-0.3,0.8726775794321988) node {$\mathbf{e}_2$};
\end{scriptsize}
\end{tikzpicture}
\Big).
\section{The case of a two-variable function}
\subsection{Coordinate-free expression of a plane secant the graph of a two-variable function}\label{sec:r as quotient n=2}
Let us now use the previous coordinate-free formalism to write the equation of a plane secant the graph of a two-variable function, that is a function $f:\Omega\subseteq\mathbb{E}_2\to\mathbb{R}$ defined on a subset of the two-dimensional Euclidean space $\mathbb{E}_2$.
A plane passing through the three points of $\mathbb{E}_2\oplus\mathbb{R}\simeq\mathbb{R}^3$
\[
\big(\mathbf{a},f(\mathbf{a})\big)
=
\big(\alpha_1\mathbf{e}_1+\alpha_2\mathbf{e}_2,f(\mathbf{a})\big)
\ \ , \ \
\big(\mathbf{b},f(\mathbf{b})\big)
=
\big(\beta_1\mathbf{e}_1+\beta_2\mathbf{e}_2,f(\mathbf{b})\big)
\]
\[
\big(\mathbf{c},f(\mathbf{c})\big)
=
\big(\gamma_1\mathbf{e}_1+\gamma_2\mathbf{e}_2,f(\mathbf{c})\big)
\]
can be represented by the Cartesian relation
\begin{equation}
\label{eq:secant plane}
\det
\left(
\begin{array}{ccc}
	x-\alpha_1 & y-\alpha_2 & z-f(\mathbf{a})\\
	\beta_1-\alpha_1 & \beta_2-\alpha_2 & f(\mathbf{b})-f(\mathbf{a})\\
	\gamma_1-\alpha_1 & \gamma_2-\alpha_2 & f(\mathbf{c})-f(\mathbf{a})
\end{array}
\right)
=0
\end{equation}
between the real variables $x,y,z\in\mathbb{R}$. That determinant can be rewritten by a  Laplace expansion as follows
\begin{align*}
&
\big[z-f(\mathbf{a})\big]
\det
\left(
\begin{array}{cc}
	\beta_1-\alpha_1 & \beta_2-\alpha_2 \\
	\gamma_1-\alpha_1 & \gamma_2-\alpha_2 
\end{array}
\right)
+\\
& -
\big[f(\mathbf{b})-f(\mathbf{a})\big]
\det
\left(
\begin{array}{cc}
	x-\alpha_1 & y-\alpha_2 \\
	\gamma_1-\alpha_1 & \gamma_2-\alpha_2 
\end{array}
\right)
+\\
& +
\big[f(\mathbf{c})-f(\mathbf{a})\big]
\det
\left(
\begin{array}{cc}
	x-\alpha_1 & y-\alpha_2 \\
	\beta_1-\alpha_1 & \beta_2-\alpha_2 
\end{array}
\right)
\end{align*}
Then, in $\mathbb{G}_2$ the equation~(\ref{eq:secant plane}) becomes 
\begin{align*}
&
\big[z-f(\mathbf{a})\big]
\big[(\mathbf{b}-\mathbf{a})\wedge(\mathbf{c}-\mathbf{a})\big]
(I_2)^{-1}
-
\big[f(\mathbf{b})-f(\mathbf{a})\big]
\big[(\mathbf{x}-\mathbf{a})\wedge(\mathbf{c}-\mathbf{a})\big]
(I_2)^{-1}
+\\
& +
\big[f(\mathbf{c})-f(\mathbf{a})\big]
\big[(\mathbf{x}-\mathbf{a})\wedge(\mathbf{b}-\mathbf{a})\big]
(I_2)^{-1}
=
0\ ,
\end{align*}
being $\mathbf{x}=x\mathbf{e}_1+y\mathbf{e}_2\in\mathbb{E}_2$, and $(\mathbf{x},z)\in\mathbb{E}_2\oplus\mathbb{R}$. The foregoing relation is equivalent, in $\mathbb{G}_2$, to
\begin{align*}
\big[z-f(\mathbf{a})\big]
\big[(\mathbf{b}-\mathbf{a})\wedge(\mathbf{c}-\mathbf{a})\big]
=
&
\big[f(\mathbf{c})-f(\mathbf{a})\big]
\big[(\mathbf{b}-\mathbf{a})\wedge(\mathbf{x}-\mathbf{a})\big]
+\\
& 
-
\big[f(\mathbf{b})-f(\mathbf{a})\big]
\big[(\mathbf{c}-\mathbf{a})\wedge(\mathbf{x}-\mathbf{a})\big]
\end{align*}
Let us define 
\begin{center}
$
\Delta_{(\mathbf{a},\mathbf{b},\mathbf{c})}
=
(\mathbf{b}-\mathbf{a})\wedge(\mathbf{c}-\mathbf{a})
\ 
\Big(=
(\mathbf{a}-\mathbf{b})\wedge(\mathbf{b}-\mathbf{c})
\Big)
$. 
\end{center}
We observe that 
\[
\displaystyle 
\tau_2 
=
\frac{1}{2}
\det
\left(
\begin{array}{cc}
	\beta_1-\alpha_1 & \beta_2-\alpha_2 \\
	\gamma_1-\alpha_1 & \gamma_2-\alpha_2 
\end{array}
\right)
=
\frac{1}{2}
\Delta_{(\mathbf{a},\mathbf{b},\mathbf{c})}
(I_2)^{-1}
\]
is the oriented area of the triangle having vertices $\mathbf{a}$, $\mathbf{b}$, and $\mathbf{c}$ (the subscript ``$2$'' in $\tau_2$ anticipate the use of $\tau_n$ as the hypervolume of a $n$-dimensional simplex). So, we can write $\Delta_{(\mathbf{a},\mathbf{b},\mathbf{c})}= 2\tau_2 I_2$, and $\big(\Delta_{(\mathbf{a},\mathbf{b},\mathbf{c})}\big)^{-1}= \frac{1}{2\tau_2}(I_2)^{-1}$. Then, the equation of the secant plane~(\ref{eq:secant plane}) becomes
\begin{align*}
& z =\\ 
= & 
f(\mathbf{a})
+ \frac{1}{2\tau_2}
\Big\{
\big[f(\mathbf{c})
\kern-3pt - \kern-3pt 
f(\mathbf{a})\big]
\big[(\mathbf{b}-\mathbf{a})\wedge(\mathbf{x}-\mathbf{a})\big]
-
\big[f(\mathbf{b})
\kern-3pt - \kern-3pt 
f(\mathbf{a})\big]
\big[(\mathbf{c}-\mathbf{a})\wedge(\mathbf{x}-\mathbf{a})\big]
\Big\}(I_2)^{-1}\\
= & 
f(\mathbf{a})
\kern-2pt + \kern-2pt
 \frac{1}{2\tau_2}
\Big\{
\big[f(\mathbf{b})
\kern-3pt - \kern-3pt 
f(\mathbf{a})\big]
\big[(\mathbf{c}-\mathbf{a})(I_2)^{-1}\big]
\kern-3pt \cdot \kern-2pt 
(\mathbf{x}
\kern-3pt - \kern-3pt 
\mathbf{a})
-
\big[f(\mathbf{c})
\kern-3pt - \kern-3pt 
f(\mathbf{a})\big]
\big[(\mathbf{b}-\mathbf{a})(I_2)^{-1}\big]
\kern-3pt \cdot \kern-2pt 
(\mathbf{x}
\kern-3pt - \kern-3pt 
\mathbf{a})
\kern-2pt \Big\}\\
= & 
f(\mathbf{a})
+ \frac{1}{2\tau_2}
\Big\{
\big[f(\mathbf{b})-f(\mathbf{a})\big]
\big[(\mathbf{c}-\mathbf{a})(I_2)^{-1}\big]
-
\big[f(\mathbf{c})-f(\mathbf{a})\big]
\big[(\mathbf{b}-\mathbf{a})(I_2)^{-1}\big]
\Big\}\cdot(\mathbf{x}-\mathbf{a})\\
= & 
f(\mathbf{a})
+ \frac{1}{2\tau_2}
\Big\{\Big[
\big[f(\mathbf{b})-f(\mathbf{a})\big]
(\mathbf{c}-\mathbf{a})
-
\big[f(\mathbf{c})-f(\mathbf{a})\big]
(\mathbf{b}-\mathbf{a})
\Big](I_2)^{-1}
\Big\}\cdot(\mathbf{x}-\mathbf{a})\\
= & 
f(\mathbf{a})
+ 
\Big\{\Big[
\big[f(\mathbf{b})-f(\mathbf{a})\big]
(\mathbf{c}-\mathbf{a})
-
\big[f(\mathbf{c})-f(\mathbf{a})\big]
(\mathbf{b}-\mathbf{a})\Big]
\big(\Delta_{(\mathbf{a},\mathbf{b},\mathbf{c})}\big)^{-1}
\Big\}\cdot(\mathbf{x}-\mathbf{a})
\end{align*} 
Thus, we have that 
\[
\mathbf{r}_{f_{(\mathbf{a},\mathbf{b},\mathbf{c})}}
=
\Big[
\big[f(\mathbf{b})-f(\mathbf{a})\big]
(\mathbf{c}-\mathbf{a})
-
\big[f(\mathbf{c})-f(\mathbf{a})\big]
(\mathbf{b}-\mathbf{a})
\Big]\big(\Delta_{(\mathbf{a},\mathbf{b},\mathbf{c})}\big)^{-1}\ ,
\]
which is, in fact, the Clifford quotient between the multi-difference vector 
\[
\mathbf{\Delta} f_{(\mathbf{a},\mathbf{b},\mathbf{c})}
=
\big[f(\mathbf{b})-f(\mathbf{a})\big]
(\mathbf{c}-\mathbf{a})
-
\big[f(\mathbf{c})-f(\mathbf{a})\big]
(\mathbf{b}-\mathbf{a})\in\mathbb{E}_2\ ,
\]
and the bivector
\[
\Delta_{(\mathbf{a},\mathbf{b},\mathbf{c})}
=
(\mathbf{b}-\mathbf{a})\wedge(\mathbf{c}-\mathbf{a})
=
(\mathbf{a}-\mathbf{b})\wedge(\mathbf{b}-\mathbf{c})
=
2\tau_2 I_2\ .
\]
\begin{Rem}\label{rem:symmetry of the ratio n=2}
As 
$
\mathbf{r}_{f_{(\mathbf{v}_1,\mathbf{v}_2,\mathbf{v}_3)}}
=
\Big(\mathbf{\Delta} f_{(\mathbf{v}_1,\mathbf{v}_2,\mathbf{v}_3)}\Big)
\big(\Delta_{(\mathbf{v}_1,\mathbf{v}_2,\mathbf{v}_3)}\big)^{-1}
$, you can verify that 
\[
\mathbf{r}_{f_{(\mathbf{v}_{\sigma_1},\mathbf{v}_{\sigma_2},\mathbf{v}_{\sigma_3})}}
=
\mathbf{r}_{f_{(\mathbf{v}_1,\mathbf{v}_2,\mathbf{v}_3)}}
\]
for every triple of non collinear vectors $\mathbf{v}_1,\mathbf{v}_2,\mathbf{v}_3\in\Omega$, and every permutation $\sigma\in\mathcal{S}_3$ of the set $\{1,2,3\}$. In other words, vector $\mathbf{r}_{f_{(\mathbf{a},\mathbf{b},\mathbf{c})}}$ is a totally symmetric function of its vector arguments $\mathbf{a},\mathbf{b},\mathbf{c}\in\mathbb{E}_2$, as the scalar difference quotient $\displaystyle \frac{\Delta f_{(a,b)}}{\Delta_{(a,b)}}$ is a totally symmetric function of its scalar arguments $a,b\in\mathbb{R}$.
\end{Rem}
\begin{Rem}\label{rem:multi diff quotient n=2}
In the case of the Example~\ref{exmp:Schwarz}, we have that
\begin{itemize}
	\item 
	$
	\mathbf{b}-\mathbf{a}
	=-\alpha \mathbf{e}_1+\beta\mathbf{e}_2
	$
	\ , \ 
	$
	\mathbf{c}-\mathbf{a}
	=\alpha \mathbf{e}_1+\beta\mathbf{e}_2
	$
	\item
	$
	f(\mathbf{b})-f(\mathbf{a})
	=
	\sqrt{1-\alpha^2}-1
	=
	f(\mathbf{c})-f(\mathbf{a})
	$
	\item
	$
	\mathbf{\Delta} f_{(\mathbf{a},\mathbf{b},\mathbf{c})}
	=
	2\alpha \big(\sqrt{1-\alpha^2}-1\big) \mathbf{e}_1
	$
	\item
	$
	\Delta_{(\mathbf{a},\mathbf{b},\mathbf{c})}
	=
	-2\alpha\beta \mathbf{e}_1\mathbf{e}_2
	$
\end{itemize}
So, 
\[
\mathbf{r}_{f_{(\mathbf{a},\mathbf{b},\mathbf{c})}}
=
\frac{\sqrt{1-\alpha^2}-1}{\beta}\mathbf{e}_2
\]
\end{Rem}

\subsection{The vector $\mathbf{r}_{f_{(\mathbf{a},\mathbf{b},\mathbf{c})}}$ as linear combination of vectors $\mathbf{a}$, $\mathbf{b}$, and $\mathbf{c}$}\label{sec:r as linear combination}
Given an oriented triangle in $\mathbb{E}_2$ whose vertices are the ordered vectors $\mathbf{a}$, $\mathbf{b}$, and $\mathbf{c}$ in $\mathbb{E}_2$, we define
\[
\mathbf{\Delta}\mathbf{b}
=
\mathbf{b}-\mathbf{a}
\ , \ 
\mathbf{\Delta}\mathbf{c}
=
\mathbf{c}-\mathbf{a}
\ , \ 
\Delta f_{(\mathbf{b})}
=
f(\mathbf{b})-f(\mathbf{a})
\ , \ 
\textrm{ and }\ 
\Delta f_{(\mathbf{c})}
=
f(\mathbf{c})-f(\mathbf{a})\ .
\]
Let us recall that 
\begin{itemize}
	\item 
	$\Delta_{(\mathbf{a},\mathbf{b},\mathbf{c})}
=
(\mathbf{\Delta}\mathbf{b})\wedge(\mathbf{\Delta}\mathbf{c})=2\tau_2 I_2$
	\item 
	$\left(\Delta_{(\mathbf{a},\mathbf{b},\mathbf{c})}\right)^2
	= -4\tau_2^2=|\mathbf{\Delta}\mathbf{b}|^2 |\mathbf{\Delta}\mathbf{c}|^2-\big(\mathbf{\Delta}\mathbf{b}\cdot \mathbf{\Delta}\mathbf{c}\big)^2$
	\item 
	$\displaystyle \left(\Delta_{(\mathbf{a},\mathbf{b},\mathbf{c})}\right)^{-1}= \frac{1}{\left(\Delta_{(\mathbf{a},\mathbf{b},\mathbf{c})}\right)^2}\Delta_{(\mathbf{a},\mathbf{b},\mathbf{c})}
=-\frac{1}{4\tau_2^2}\Delta_{(\mathbf{a},\mathbf{b},\mathbf{c})}= 
-\frac{1}{2\tau_2}I_2 
$
\item $\mathbf{\Delta} f_{(\mathbf{a},\mathbf{b},\mathbf{c})}
=
\Delta f_{(\mathbf{b})}\mathbf{\Delta}\mathbf{c}
-
\Delta f_{(\mathbf{c})}\mathbf{\Delta}\mathbf{b}
$
\item 
$
\mathbf{u}(\mathbf{v}\wedge\mathbf{w})
=
(\mathbf{u}\cdot\mathbf{v})\mathbf{w}
-
(\mathbf{u}\cdot\mathbf{w})\mathbf{v}
$ for all $\mathbf{u}$, $\mathbf{v}$, and $\mathbf{w}$ in $\mathbb{E}_2$.
\end{itemize}
So, we can write
\begin{align*}
\mathbf{r}_{f_{(\mathbf{a},\mathbf{b},\mathbf{c})}}
= & 
\Big(\mathbf{\Delta} f_{(\mathbf{a},\mathbf{b},\mathbf{c})}\Big)
\big(\Delta_{(\mathbf{a},\mathbf{b},\mathbf{c})}\big)^{-1}
=
\Big(
\Delta f_{(\mathbf{b})}\mathbf{\Delta}\mathbf{c}
-
\Delta f_{(\mathbf{c})}\mathbf{\Delta}\mathbf{b}
\Big)\big[(\mathbf{\Delta}\mathbf{b})\wedge(\mathbf{\Delta}\mathbf{c})\big]^{-1}\\
= & 
-\frac{1}{4\tau_2^2}
\Big(
\Delta f_{(\mathbf{b})}\mathbf{\Delta}\mathbf{c}
-
\Delta f_{(\mathbf{c})}\mathbf{\Delta}\mathbf{b}
\Big)\big[(\mathbf{\Delta}\mathbf{b})\wedge(\mathbf{\Delta}\mathbf{c})\big]\\
= &
-\frac{1}{4\tau_2^2}
\Big\{
\Delta f_{(\mathbf{b})}\mathbf{\Delta}\mathbf{c}
\big[(\mathbf{\Delta}\mathbf{b})\wedge(\mathbf{\Delta}\mathbf{c})\big]
-
\Delta f_{(\mathbf{c})}\mathbf{\Delta}\mathbf{b}
\big[(\mathbf{\Delta}\mathbf{b})\wedge(\mathbf{\Delta}\mathbf{c})\big]
\Big\}\\
= & 
-\frac{1}{4\tau_2^2}
\Big\{
\Delta f_{(\mathbf{b})}
(\mathbf{\Delta}\mathbf{c}\cdot \mathbf{\Delta}\mathbf{b})\mathbf{\Delta}\mathbf{c})
-
\Delta f_{(\mathbf{b})}
|\mathbf{\Delta}\mathbf{c}|^2
\mathbf{\Delta}\mathbf{b}\Big\}+\\
& +
\frac{1}{4\tau_2^2}
\Big\{
\Delta f_{(\mathbf{c})}|\mathbf{\Delta}\mathbf{b}|^2
\mathbf{\Delta}\mathbf{c}
-
\Delta f_{(\mathbf{c})}
(\mathbf{\Delta}\mathbf{b}\cdot \mathbf{\Delta}\mathbf{c})\mathbf{\Delta}\mathbf{b})\big]\\
= &
\frac{\Delta f_{(\mathbf{b})}
|\mathbf{\Delta}\mathbf{c}|^2-\Delta f_{(\mathbf{c})}
(\mathbf{\Delta}\mathbf{b}\cdot \mathbf{\Delta}\mathbf{c})}{|\mathbf{\Delta}\mathbf{b}|^2 |\mathbf{\Delta}\mathbf{c}|^2-\big(\mathbf{\Delta}\mathbf{b}\cdot \mathbf{\Delta}\mathbf{c}\big)^2}
\mathbf{\Delta}\mathbf{b}
+
\frac{\Delta f_{(\mathbf{c})}
|\mathbf{\Delta}\mathbf{b}|^2-\Delta f_{(\mathbf{b})}
(\mathbf{\Delta}\mathbf{b}\cdot \mathbf{\Delta}\mathbf{c})}{|\mathbf{\Delta}\mathbf{b}|^2 |\mathbf{\Delta}\mathbf{c}|^2-\big(\mathbf{\Delta}\mathbf{b}\cdot \mathbf{\Delta}\mathbf{c}\big)^2}
\mathbf{\Delta}\mathbf{c}
\end{align*}

\subsection{Mirroring vectors and points I}
\subsubsection{Vector mirrored by a $1$-dimensional linear subspace}
We recall that to each non zero vector $\mathbf{u}\in\mathbb{E}_n$ we can associate the $1$-dimensional linear subspace of $\mathbb{E}_n$ $span\{\mathbf{u}\}=\{\lambda\mathbf{u}\ : \ \lambda\in\mathbb{R}\}=\mathbb{R}\mathbf{u}$.
Besides, every non zero vector $\mathbf{u}\in\mathbb{E}_n$ is invertible in $\mathbb{G}_n$, and $\displaystyle \mathbf{u}^{-1}=\frac{1}{|\mathbf{u}|^2}\mathbf{u}$. 
Then, given a vector $\mathbf{v} \in\mathbb{E}_n$, we can write
\begin{equation}
\label{eq:decomposition 1}
\mathbf{v}
=
\mathbf{v}\mathbf{u}\mathbf{u}^{-1}
=
(\mathbf{v}\mathbf{u})\mathbf{u}^{-1}
=
(\mathbf{v}\cdot\mathbf{u}+\mathbf{v}\wedge\mathbf{u})\mathbf{u}^{-1}
=
\frac{\mathbf{v}\cdot\mathbf{u}}{|\mathbf{u}|^2}\mathbf{u}+(\mathbf{v}\wedge\mathbf{u})\mathbf{u}^{-1}
\end{equation}
\begin{Rem}\label{rem:quotient between scalars and vectors}
If $\mathbf{v}\in\mathbb{E}_n$ is invertible, and $\alpha\in\mathbb{R}$ the expression
\[
\frac{\alpha}{\mathbf{v}}
\]
is unambiguous, because in $\mathbb{G}_n$ scalars commute with vectors (and with any other element, indeed). As a matter of fact,
\[
\frac{\alpha}{\mathbf{v}}
=
\alpha \mathbf{v}^{-1}
=
\mathbf{v}^{-1}\alpha
=
\frac{\alpha}{\mathbf{v}\cdot\mathbf{v}}\mathbf{v}
=
\frac{\alpha}{|\mathbf{v}|^2}\mathbf{v}
=
\cdots
\] 
\end{Rem}
\begin{Rem}\label{rem:rejection is orthogonal 1}
Notice that $(\mathbf{v}\wedge\mathbf{u})\mathbf{u}^{-1}$ is orthogonal to $\mathbf{u}$. As a matter of fact,
\begin{align*}
4\big[(\mathbf{v}\wedge\mathbf{u})\mathbf{u}^{-1}\big]\cdot \mathbf{u}
& =
2(\mathbf{v}\wedge\mathbf{u})\mathbf{u}^{-1} \mathbf{u}
+
2\mathbf{u}(\mathbf{v}\wedge\mathbf{u})\mathbf{u}^{-1}\\
& =
2(\mathbf{v}\wedge\mathbf{u})
+
\mathbf{u}(\mathbf{v}\mathbf{u}-\mathbf{u}\mathbf{v})\mathbf{u}^{-1}\\
& =
2(\mathbf{v}\wedge\mathbf{u})
+
\mathbf{u}\mathbf{v}-\mathbf{u}\mathbf{u}\mathbf{v}\mathbf{u}^{-1}\\
& =
2(\mathbf{v}\wedge\mathbf{u})
+
\mathbf{u}\mathbf{v}-\mathbf{v}\mathbf{u}
=
2(\mathbf{v}\wedge\mathbf{u})
+
2(\mathbf{u}\wedge\mathbf{v})
=
0
\end{align*}
\end{Rem}
As $\displaystyle \frac{\mathbf{v}\cdot\mathbf{u}}{|\mathbf{u}|^2}\mathbf{u}$ is parallel to $\mathbf{u}$, and  $(\mathbf{v}\wedge\mathbf{u})\mathbf{u}^{-1}$ is orthogonal to $\mathbf{u}$, we can consider relation~(\ref{eq:decomposition 1}) as the decomposition of $\mathbf{v}=\mathbf{v}_{\parallel}+\mathbf{v}_{\bot}$ into its orthogonal projection $\displaystyle \mathbf{v}_{\parallel}= \frac{\mathbf{v}\cdot\mathbf{u}}{|\mathbf{u}|^2}\mathbf{u}=
(\mathbf{u}\cdot\mathbf{v})\mathbf{u}^{-1}$ parallel to the line $\mathbb{R}\mathbf{u}=\mathcal{L}_{(\mathbf{0},\mathbf{u})}$, and its rejection $\mathbf{v}_{\bot}=(\mathbf{v}\wedge\mathbf{u})\mathbf{u}^{-1}$ orthogonal to $\mathbb{R}\mathbf{u}$.
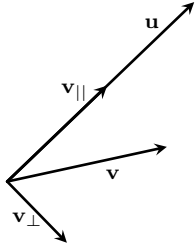
\begin{center}
\begin{tikzpicture}[line cap=round,line join=round,>=triangle 45,x=1.0cm,y=1.0cm]
\clip(-0.1,-0.9) rectangle (3.,2.5);
\draw [->,line width=1.pt,-stealth] (0.,0.) -- (2.46,2.38);
\draw [->,line width=1.pt,-stealth] (0.,0.) -- (2.1,0.46);
\draw [->,line width=1.1pt,-stealth] (0.,0.) -- (1.3145756230795496,1.271825196312735);
\draw [->,line width=1.pt,-stealth] (0.,0.) -- (0.7854243769204505,-0.811825196312735);
\begin{scriptsize}
\draw[color=black] (1.9,2.1) node {$\mathbf{u}$};
\draw[color=black] (1.3898934781246584,0.1) node {$\mathbf{v}$};
\draw[color=black] (0.887511680376797,1.183604933727941) node {$\mathbf{v}_{||}$};
\draw[color=black] (0.26659035731764236,-0.5) node {$\mathbf{v}_{\bot}$};
\end{scriptsize}
\end{tikzpicture}
\end{center}
Thus, the vector $\mathbf{\hat{v}}$, obtained by mirroring $\mathbf{v}$ through $\mathcal{L}_{(\mathbf{0},\mathbf{u})}$, 
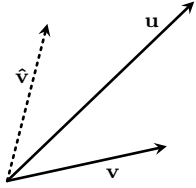
\begin{center}
\begin{tikzpicture}[line cap=round,line join=round,>=triangle 45,x=1.0cm,y=1.0cm]
\clip(-0.1,-0.1) rectangle (3,2.6);
\draw [->,line width=1.pt,-stealth] (0.,0.) -- (2.46,2.38);
\draw [->,line width=1.pt,-stealth] (0.,0.) -- (2.1,0.46);
\draw [->,line width=1.pt,-stealth,dotted] (0.,0.) -- (0.529151246159099,2.08365039262547);
\begin{scriptsize}
\draw[color=black] (1.9,2.1) node {$\mathbf{u}$};
\draw[color=black] (1.3898934781246586,0.1) node {$\mathbf{v}$};
\draw[color=black] (0.2,1.392460287847839) node {$\mathbf{\hat{v}}$};
\end{scriptsize}
\end{tikzpicture}
\end{center}
can be written as
\begin{align*}
\mathbf{\hat{v}}
= &
\mathbf{v}_{\parallel}-\mathbf{v}_{\bot}
=
(\mathbf{u}\cdot\mathbf{v})\mathbf{u}^{-1}
-
(\mathbf{v}\wedge\mathbf{u})\mathbf{u}^{-1}
=
(\mathbf{u}\cdot\mathbf{v})\mathbf{u}^{-1}
+
(\mathbf{u}\wedge\mathbf{v})\mathbf{u}^{-1}\\
= &
(\mathbf{u}\cdot\mathbf{v}+\mathbf{u}\wedge\mathbf{v})\mathbf{u}^{-1}
= 
\mathbf{u}\mathbf{v}\mathbf{u}^{-1}\\
= &
\big[2(\mathbf{u}\cdot\mathbf{v})-\mathbf{v}\mathbf{u}\big]\mathbf{u}^{-1}
=
2(\mathbf{u}\cdot\mathbf{v})\mathbf{u}^{-1}-\mathbf{v}
=
2\frac{\mathbf{u}\cdot\mathbf{v}}{|\mathbf{u}|^2}\mathbf{u}-\mathbf{v}\in span\{\mathbf{u},\mathbf{v}\}\ ,
\end{align*}
as for each $\mathbf{u},\mathbf{v}\in\mathbb{E}_n$\ \ $\mathbf{u}\mathbf{v}=2(\mathbf{u}\cdot\mathbf{v})-\mathbf{v}\mathbf{u}$. Moreover, $|\mathbf{\hat{v}}|= |\mathbf{v}|$; as a matter of fact,
\begin{center}
$
|\mathbf{\hat{v}}|^2
=
\mathbf{\hat{v}}\mathbf{\hat{v}}
=
\mathbf{u}\mathbf{v}\mathbf{u}^{-1}
\mathbf{u}\mathbf{v}\mathbf{u}^{-1}
=
\mathbf{u}|\mathbf{v}|^2\mathbf{u}^{-1}
=
|\mathbf{v}|^2
$.
\end{center}
\subsubsection{Point mirrored by a line in $\mathbb{E}_n$ (with $n\ge 2$)}\label{subsubsec: point mirrored by a line}
Given three non collinear points $\mathbf{a}$, $\mathbf{b}$, and $\mathbf{c}$ in $\mathbb{E}_n$, we want to mirror point~$\mathbf{a}$ by the line $\mathcal{L}_{(\mathbf{b},\mathbf{c})}$ passing through the points $\mathbf{b}$, and $\mathbf{c}$. We denote by $\bar{\mathbf{a}}$ that mirrored point.
We can express  $\bar{\mathbf{a}}$ by computing the vector $\mathbf{\hat{v}}$, obtained by mirroring vector $\mathbf{v}=\mathbf{a}-\mathbf{b}$ by the non zero vector $\mathbf{c}-\mathbf{b}$.
Thus, the  reflected point $\bar{\mathbf{a}}$
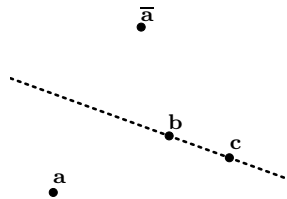
\begin{center}
\begin{tikzpicture}[line cap=round,line join=round,>=triangle 45,x=0.6cm,y=0.5cm]
\clip(1.,2.) rectangle (7.,7.1);
\draw [line width=1.pt,dotted,domain=1.:7.] plot(\x,{(--7.4296-0.58*\x)/1.32});
\begin{scriptsize}
\draw [fill=black](1.94,2.16) circle (1.5pt);
\draw[color=black] (2.08,2.53) node {$\mathbf{a}$};
\draw [fill=black] (4.48,3.66) circle (1.5pt);
\draw[color=black] (4.62,4.03) node {$\mathbf{b}$};
\draw [fill=black] (5.8,3.08) circle (1.5pt);
\draw[color=black] (5.94,3.45) node {$\mathbf{c}$};
\draw[color=black] (-7.7,8.85) node {$f$};
\draw [fill=black] (3.8669347700596517,6.545437752549548) circle (1.5pt);
\draw[color=black] (4.,6.9) node {$\overline{\mathbf{a}}$};
\end{scriptsize}
\end{tikzpicture}
\end{center}
can be expressed using the geometric Clifford product in $\mathbb{G}_n$
\begin{align*}
\bar{\mathbf{a}}
= &
\mathbf{b}
+
(\mathbf{c}-\mathbf{b})(\mathbf{a}-\mathbf{b})(\mathbf{c}-\mathbf{b})^{-1}
=
\mathbf{b}
-
2\big[(\mathbf{c}-\mathbf{b})\cdot(\mathbf{\Delta}\mathbf{b})\big](\mathbf{c}-\mathbf{b})^{-1}
+(\mathbf{b}-\mathbf{a})\\
= &
2\mathbf{b}
-
2\big[(\mathbf{c}-\mathbf{b})\cdot(\mathbf{\Delta}\mathbf{b})\big](\mathbf{c}-\mathbf{b})^{-1}
-\mathbf{a}
=
2\mathbf{b}
-2
\frac{(\mathbf{c}-\mathbf{b})\cdot(\mathbf{\Delta}\mathbf{b})}
{\mathbf{c}-\mathbf{b}}
-\mathbf{a}\\
= &
2\mathbf{b}
-2
\frac{(\mathbf{c}-\mathbf{b})\cdot(\mathbf{\Delta}\mathbf{b})}
{|\mathbf{c}-\mathbf{b}|^2}
(\mathbf{c}-\mathbf{b})
-\mathbf{a}
\in \mathcal{L}_{(\mathbf{a},\mathbf{b},\mathbf{c})}\ .
\end{align*}
Thus, 
\begin{align*}
\bar{\mathbf{a}}-\mathbf{a}
& =
2(\mathbf{\Delta}\mathbf{b})
-2\big[(\mathbf{c}-\mathbf{b})\cdot(\mathbf{\Delta}\mathbf{b})\big](\mathbf{c}-\mathbf{b})^{-1}\\
& =
2(\mathbf{\Delta}\mathbf{b})
(\mathbf{c}-\mathbf{b})(\mathbf{c}-\mathbf{b})^{-1}
-2\big[(\mathbf{\Delta}\mathbf{b})\cdot(\mathbf{c}-\mathbf{b})\big](\mathbf{c}-\mathbf{b})^{-1}\\
& =
2\Big\{(\mathbf{\Delta}\mathbf{b})(\mathbf{c}-\mathbf{b})
-\big[(\mathbf{\Delta}\mathbf{b})\cdot(\mathbf{c}-\mathbf{b})\big]\Big\}(\mathbf{c}-\mathbf{b})^{-1}\\
& =
2\big[(\mathbf{\Delta}\mathbf{b})\wedge(\mathbf{c}-\mathbf{b})\big]
(\mathbf{c}-\mathbf{b})^{-1}
\end{align*}
and $|\bar{\mathbf{a}}-\mathbf{b}|=|\mathbf{a}-\mathbf{b}|$.
\begin{Rem}\label{rem:ortho-diagonals I}
By using Remark~\ref{rem:rejection is orthogonal 1}, you can verify that vectors $\bar{\mathbf{a}}-\mathbf{a}$ and $\mathbf{b}-\mathbf{c}$ are mutually orthogonal. Moreover,
\begin{align*}
(\bar{\mathbf{a}}-\mathbf{a})\wedge (\mathbf{c}-\mathbf{b})
& =
(\bar{\mathbf{a}}-\mathbf{a})(\mathbf{c}-\mathbf{b})
=
2\big[(\mathbf{\Delta}\mathbf{b})\wedge(\mathbf{c}-\mathbf{b})\big]
(\mathbf{c}-\mathbf{b})^{-1}(\mathbf{c}-\mathbf{b})\\
& =
2(\mathbf{b}-\mathbf{a})\wedge(\mathbf{c}-\mathbf{b})
=
2
\Delta_{(\mathbf{a},\mathbf{b},\mathbf{c})}
\end{align*}
\end{Rem}
\subsection{The mean multi-difference vector}\label{sec:mean multi diff vector}  
Let us define the ``mean multi-difference vector'',
\[
\overline{\Delta}f_{(\mathbf{a},\mathbf{b},\mathbf{c})}
=
\frac{1}{2}
\Big(
\Delta f_{(\mathbf{a},\mathbf{b},\mathbf{c})}
+
\Delta f_{(\bar{\mathbf{a}},\mathbf{c},\mathbf{b})}
\Big)=
\frac{1}{2}
\Big(
\Delta f_{(\mathbf{a},\mathbf{b},\mathbf{c})}
-
\Delta f_{(\bar{\mathbf{a}},\mathbf{b},\mathbf{c})}
\Big)
\ .
\]
\begin{Lem}\label{lem:mean multi-difference vector}
\[
\overline{\Delta}f_{(\mathbf{a},\mathbf{b},\mathbf{c})}
=
\frac{1}{2}
\Big\{
\big[f(\bar{\mathbf{a}})-f(\mathbf{a})\big]
(\mathbf{c}-\mathbf{b})
-
\big[f(\mathbf{c})-f(\mathbf{b})\big]
(\bar{\mathbf{a}}-\mathbf{a})
\Big\}\ .
\]
\end{Lem}
{\bf Proof of Lemma~\ref{lem:mean multi-difference vector}.}
\begin{align*}
& 
\overline{\Delta}f_{(\mathbf{a},\mathbf{b},\mathbf{c})}=\\
= & 
\frac{1}{2}
\Big\{
\big[f(\mathbf{b})-f(\mathbf{a})\big]
(\mathbf{c}-\mathbf{a})
-
\big[f(\mathbf{c})-f(\mathbf{a})\big]
(\mathbf{b}-\mathbf{a})
\Big\}+\\
& -
\frac{1}{2}
\Big\{
\big[f(\mathbf{b})-f(\bar{\mathbf{a}})\big]
(\mathbf{c}-\bar{\mathbf{a}})
-
\big[f(\mathbf{c})-f(\bar{\mathbf{a}})\big]
(\mathbf{b}-\bar{\mathbf{a}})
\Big\}\\
= &
\frac{1}{2}
\Big\{
\big[f(\bar{\mathbf{a}})-f(\mathbf{a})\big]\mathbf{c}
+
\big[f(\mathbf{c})-f(\mathbf{b})\big]\mathbf{a}
+
\big[f(\mathbf{a})-f(\bar{\mathbf{a}})\big]\mathbf{b}
+
\big[f(\mathbf{b})-f(\mathbf{c})\big]\bar{\mathbf{a}}
\Big\}\ \square
\end{align*} 
The foregoing Lemma and Remark~\ref{rem:quotient between scalars and vectors} allows us to obtain a simple expression for the vector corresponding to the mean multi-difference quotient\\
$\overline{\mathbf{r}}_{f_{(\mathbf{a},\mathbf{b},\mathbf{c})}}
= \Big(\overline{\mathbf{\Delta}}f_{(\mathbf{a},\mathbf{b},\mathbf{c})}\Big)\big(\Delta_{(\mathbf{a},\mathbf{b},\mathbf{c})}\big)^{-1}$, which strongly recall the usual difference quotients (as  anticipated in the abstract and in the introduction of this work).
\begin{Prop}\label{prop:mean multi-difference quotient n=2}
\[
\overline{\mathbf{r}}_{f_{(\mathbf{a},\mathbf{b},\mathbf{c})}}
=
\frac{f(\bar{\mathbf{a}})-f(\mathbf{a})}{\bar{\mathbf{a}}-\mathbf{a}} 
+
\frac{f(\mathbf{c})-f(\mathbf{b})}{\mathbf{c}-\mathbf{b}}\ .
\]
\end{Prop}
{\bf Proof of Proposition~\ref{prop:mean multi-difference quotient n=2}.}
From Remark~\ref{rem:ortho-diagonals I} we have that $\bar{\mathbf{a}}-\mathbf{a}$ is orthogonal to  $\mathbf{c}-\mathbf{b}$, and 
$
(\bar{\mathbf{a}}-\mathbf{a})\wedge(\mathbf{c}-\mathbf{b})
=  
2 \Delta_{(\mathbf{a},\mathbf{b},\mathbf{c})}
$. So, we can write
\begin{center}
$
\Delta_{(\mathbf{a},\mathbf{b},\mathbf{c})}
=
\frac{1}{2}(\bar{\mathbf{a}}-\mathbf{a})\wedge(\mathbf{c}-\mathbf{b})
=
\frac{1}{2}(\bar{\mathbf{a}}-\mathbf{a})(\mathbf{c}-\mathbf{b})
$ 
and 
$
\big(\Delta_{(\mathbf{a},\mathbf{b},\mathbf{c})}\big)^{-1}
=
2(\mathbf{c}-\mathbf{b})^{-1}(\bar{\mathbf{a}}-\mathbf{a})^{-1}
=
-2(\bar{\mathbf{a}}-\mathbf{a})^{-1}(\mathbf{c}-\mathbf{b})^{-1}
$.
\end{center}
Then, by Lemma~\ref{lem:mean multi-difference vector}, we can write 
\begin{align*}
\overline{\mathbf{r}}_{f_{(\mathbf{a},\mathbf{b},\mathbf{c})}}
& = 
\Big(\overline{\mathbf{\Delta}}f_{(\mathbf{a},\mathbf{b},\mathbf{c})}\Big)\big(\Delta_{(\mathbf{a},\mathbf{b},\mathbf{c})}\big)^{-1}=\\
& =
\Big\{
\big[f(\bar{\mathbf{a}})-f(\mathbf{a})\big]
(\mathbf{c}-\mathbf{b})
-
\big[f(\mathbf{c})-f(\mathbf{b})\big]
(\bar{\mathbf{a}}-\mathbf{a})
\Big\}(\mathbf{c}-\mathbf{b})^{-1}(\bar{\mathbf{a}}-\mathbf{a})^{-1}\\
& =
\big[f(\bar{\mathbf{a}})-f(\mathbf{a})\big]
(\bar{\mathbf{a}}-\mathbf{a})^{-1}
+
\big[f(\mathbf{c})-f(\mathbf{b})\big]
(\mathbf{c}-\mathbf{b})^{-1}\\
& =
\frac{f(\bar{\mathbf{a}})-f(\mathbf{a})}{\bar{\mathbf{a}}-\mathbf{a}}
+
\frac{f(\mathbf{c})-f(\mathbf{b})}{\mathbf{c}-\mathbf{b}}\ ,
\end{align*}
coherently with Remark~\ref{rem:quotient between scalars and vectors}.~$\square$

\subsection{Convergence of the mean secant plane to the tangent plane}
As we have seen, in the Example, the local Schwarz paradox is due to the non existence of the limit
\arraycolsep=1.pt\def\arraystretch{0.5}
\[
\lim_{\begin{array}{c}
	\scriptstyle (\mathbf{a},\mathbf{b},\mathbf{c})\to (\mathbf{x_0},\mathbf{x_0},\mathbf{x_0})\\
	\scriptstyle \mathbf{a},\mathbf{b},\mathbf{c} \textrm{ not collinear}
\end{array}
}
\Big(\mathbf{\Delta} f_{(\mathbf{a},\mathbf{b},\mathbf{c})}\Big)
\big(\Delta_{(\mathbf{a},\mathbf{b},\mathbf{c})}\big)^{-1}\ .
\]
Here we recall\footnote{See for example~\cite{Peano1892}, \cite{Leach}, \cite{Esser}.} the definition of strong (or strict) differentiability of a multi-variable function at an internal point of its domain.
\begin{Def}\label{def:strongly diff func}
A function $f:\Omega\subseteq\mathbb{E}_n\to\mathbb{R}$ is strongly differentiable at $\mathbf{x}_0$ (a point internal to $\Omega$) if there exists a vector $\mathbf{f^*}(\mathbf{x}_0)\in\mathbb{E}_n$ such that for each $\epsilon>0$ there exists a $\delta>0$ such that if $|\mathbf{u}-\mathbf{x}_0|<\delta$ and $|\mathbf{v}-\mathbf{x}_0|<\delta$, then
\[
\Big|f(\mathbf{u})-f(\mathbf{v})-\mathbf{f^*}(\mathbf{x}_0) \cdot(\mathbf{u}-\mathbf{v})\Big|< \epsilon\ \big|\mathbf{u}-\mathbf{v}\big|\ ,
\]
being $\mathbf{u},\mathbf{v}\in\Omega$.
\end{Def}
\begin{Rem*}
We recall that, if a function is strongly differentiable at $\mathbf{x}_0$, then it is also differentiable, and the vector $\mathbf{f^*}(\mathbf{x}_0)$ coincides with the gradient $\nabla f(\mathbf{x}_0)$.
\end{Rem*}

\begin{Thm}\label{thm:n=2}
If the function $f:\Omega\subseteq\mathbb{E}_2\to\mathbb{R}$ is strongly differentiable at $\mathbf{x}_0$ (a point internal to $\Omega$), then
\arraycolsep=1.pt\def\arraystretch{0.5}
\[
\lim_{\begin{array}{c}
	\scriptstyle (\mathbf{a},\mathbf{b},\mathbf{c})\to (\mathbf{x_0},\mathbf{x_0},\mathbf{x_0})\\
	\scriptstyle \mathbf{a},\mathbf{b},\mathbf{c} \textrm{ not collinear}
\end{array}
}
\Big(\overline{\mathbf{\Delta}}f_{(\mathbf{a},\mathbf{b},\mathbf{c})}\Big)
\big(\Delta_{(\mathbf{a},\mathbf{b},\mathbf{c})}\big)^{-1}
=\nabla f(\mathbf{x}_0)\ .
\]
\end{Thm}
\begin{Rem*}
The foregoing result state that the ``mean secant plane''
\[
z= f(\mathbf{a})\ +\ \overline{\mathbf{r}}_{f_{(\mathbf{a},\mathbf{b},\mathbf{c})}} \cdot (\mathbf{x}-\mathbf{a})\ ,
\]
where $\overline{\mathbf{r}}_{f_{(\mathbf{a},\mathbf{b},\mathbf{c})}}
= \Big(\overline{\mathbf{\Delta}}f_{(\mathbf{a},\mathbf{b},\mathbf{c})}\Big)
\big(\Delta_{(\mathbf{a},\mathbf{b},\mathbf{c})}\big)^{-1}$, always converges to the tangent plane  
\[
z= f(\mathbf{x}_0)\ +\ \nabla f(\mathbf{x}_0) \cdot (\mathbf{x}-\mathbf{x}_0)\ ,
\]
as the non degenerate triangle with vertices $\mathbf{a}$, $\mathbf{b}$, $\mathbf{c}$ contracts to point $\mathbf{x}_0$.
\end{Rem*}
{\bf Proof of Theorem~\ref{thm:n=2}.}\\
Let us recall that, given three vectors $\mathbf{u}$, $\mathbf{v}$, and $\mathbf{w}$ in $\mathbb{E}_n$, if $\mathbf{u}\in span\{\mathbf{v},\mathbf{w}\}$, then $\mathbf{u}(\mathbf{v}\wedge\mathbf{w})=(\mathbf{u}\cdot\mathbf{v})\mathbf{w}-(\mathbf{u}\cdot\mathbf{w})\mathbf{v}$. So, 
we can write
\begin{align*}
\nabla f(\mathbf{x}_0)
\big(\Delta_{(\mathbf{a},\mathbf{b},\mathbf{c})}\big)
& =
\frac{1}{2}
\nabla f(\mathbf{x}_0)
\big[
(\bar{\mathbf{a}}-\mathbf{a})\wedge(\mathbf{c}-\mathbf{b})
\big]\\
& =
\frac{1}{2}
\big[\nabla f(\mathbf{x}_0)\cdot(\bar{\mathbf{a}}-\mathbf{a})\big]
(\mathbf{c}-\mathbf{b})
-
\frac{1}{2}
\big[\nabla f(\mathbf{x}_0)\cdot (\mathbf{c}-\mathbf{b})\big]
(\bar{\mathbf{a}}-\mathbf{a})
\end{align*}
As 
$
\big(\Delta_{(\mathbf{a},\mathbf{b},\mathbf{c})}\big)^{-1}
=
2(\mathbf{c}-\mathbf{b})^{-1}(\bar{\mathbf{a}}-\mathbf{a})^{-1}
=
-2(\bar{\mathbf{a}}-\mathbf{a})^{-1}(\mathbf{c}-\mathbf{b})^{-1}
$, we can write 
\begin{align*}
\nabla f(\mathbf{x}_0)
&=
\nabla f(\mathbf{x}_0)
\big(\Delta_{(\mathbf{a},\mathbf{b},\mathbf{c})}\big)
\big(\Delta_{(\mathbf{a},\mathbf{b},\mathbf{c})}\big)^{-1}
\\
& =
\big[\nabla f(\mathbf{x}_0)\cdot(\bar{\mathbf{a}}-\mathbf{a})\big]
(\bar{\mathbf{a}}-\mathbf{a})^{-1}
+
\big[\nabla f(\mathbf{x}_0)\cdot (\mathbf{c}-\mathbf{b})\big]
(\mathbf{c}-\mathbf{b})^{-1}\\
& =
\frac{\nabla f(\mathbf{x}_0)\cdot(\bar{\mathbf{a}}-\mathbf{a})}
{\bar{\mathbf{a}}-\mathbf{a}}
+
\frac{\nabla f(\mathbf{x}_0)\cdot (\mathbf{c}-\mathbf{b})}
{\mathbf{c}-\mathbf{b}}
\end{align*}
So, by Proposition~\ref{prop:mean multi-difference quotient n=2}, we can write
\begin{align*}
\overline{\mathbf{r}}_{f_{(\mathbf{a},\mathbf{b},\mathbf{c})}}
-
\nabla f(\mathbf{x}_0)
= &
\big[f(\bar{\mathbf{a}})-f(\mathbf{a})-\nabla f(\mathbf{x}_0)\cdot(\bar{\mathbf{a}}-\mathbf{a})\big]
(\bar{\mathbf{a}}-\mathbf{a})^{-1}
+\\
& +
\big[f(\mathbf{c})-f(\mathbf{b})-\nabla f(\mathbf{x}_0)\cdot(\mathbf{c}-\mathbf{b})\big]
(\mathbf{c}-\mathbf{b})^{-1}\\
= &
\frac{f(\bar{\mathbf{a}})-f(\mathbf{a})-\nabla f(\mathbf{x}_0)\cdot(\bar{\mathbf{a}}-\mathbf{a})}
{\bar{\mathbf{a}}-\mathbf{a}}+\\
& +
\frac{f(\mathbf{c})-f(\mathbf{b})-\nabla f(\mathbf{x}_0)\cdot(\mathbf{c}-\mathbf{b})}{\mathbf{c}-\mathbf{b}}\ .
\end{align*} 
As 
$
\displaystyle
\left|
\frac{\alpha}{\mathbf{v}}
\right|
=
\frac{|\alpha|}{|\mathbf{v}|}
$
for every $\alpha\in\mathbb{R}$ and invertible vector $\mathbf{v}\in\mathbb{E}_n$, we have that 
\begin{align*}
& 
\big|
\overline{\mathbf{r}}_{f_{(\mathbf{a},\mathbf{b},\mathbf{c})}}
-
\nabla f(\mathbf{x}_0)
\big|
\le\\
\le &
\frac{\big|f(\bar{\mathbf{a}})-f(\mathbf{a})-\nabla f(\mathbf{x}_0)\cdot(\bar{\mathbf{a}}-\mathbf{a})\big|}{|\bar{\mathbf{a}}-\mathbf{a}|}
+
\frac{\big|f(\mathbf{c})-f(\mathbf{b})-\nabla f(\mathbf{x}_0)\cdot(\mathbf{c}-\mathbf{b})\big|}{|\mathbf{c}-\mathbf{b}|}\ .
\end{align*}
As $f$ is strongly differentiable at $\mathbf{x}_0$, we know that, given $\epsilon>0$ there exists $\delta_\epsilon>0$ such that if  $|\mathbf{u}-\mathbf{x}_0|<\delta_\epsilon$ and $|\mathbf{v}-\mathbf{x}_0|<\delta_\epsilon$, then
\begin{center}
$
|f(\mathbf{u})-f(\mathbf{v})-\nabla f(\mathbf{x}_0)\cdot(\mathbf{u}-\mathbf{v})|<\epsilon |\mathbf{u}-\mathbf{v}|\ .
$
\end{center}
So, if we choose $\mathbf{a}$, $\mathbf{b}$, and $\mathbf{c}$ non collinear, and such that
\begin{itemize}
	\item $|\mathbf{a}-\mathbf{x}_0|< \frac{1}{3}\left(\delta_{\frac{\epsilon}{2}}\right)$,
	\item $|\mathbf{b}-\mathbf{x}_0|< \frac{1}{3}\left(\delta_{\frac{\epsilon}{2}}\right)$,
	\item $|\mathbf{c}-\mathbf{x}_0|< \delta_{\frac{\epsilon}{2}}$,
\end{itemize}
then 
\[
\big|
\overline{\mathbf{r}}_{f_{(\mathbf{a},\mathbf{b},\mathbf{c})}}
-
\nabla f(\mathbf{x}_0)
\big|<\epsilon\ ,
\]
as 
\begin{itemize}
	\item $|\mathbf{a}-\mathbf{b}|\le |\mathbf{a}-\mathbf{x}_0|+ |\mathbf{b}-\mathbf{x}_0|< \frac{2}{3}\left(\delta_{\frac{\epsilon}{2}}\right)$, 
	\item $|\bar{\mathbf{a}}-\mathbf{x}_0|\le|\bar{\mathbf{a}}-\mathbf{b}|+|\mathbf{b}-\mathbf{x}_0|=|\mathbf{a}-\mathbf{b}|+|\mathbf{b}-\mathbf{x}_0| < \delta_{\frac{\epsilon}{2}}$.~$\square$
\end{itemize}
\begin{Rem*}
The divergence phenomenon of the local Schwarz paradox is due to the fact that vectors $\mathbf{b}-\mathbf{a}$  and~$\mathbf{c}-\mathbf{a}$ are not mutually orthogonal, in general. On the contrary, vectors $\mathbf{a}-\bar{\mathbf{a}}$ and $\mathbf{b}-\mathbf{c}$ are always orthogonal. Moreover, the following crucial identities hold
\begin{center}
$
(\bar{\mathbf{a}}-\mathbf{a})(\mathbf{c}-\mathbf{b})
=
(\bar{\mathbf{a}}-\mathbf{a})\wedge(\mathbf{c}-\mathbf{b})
=
2
(\mathbf{b}-\mathbf{a})\wedge(\mathbf{c}-\mathbf{a})
$.
\end{center}
Such key properties of the mirrored points will be used in extending the foregoing convergence result to higher dimensions.
\end{Rem*}
\begin{Rem*}
In the case of the Example~\ref{exmp:Schwarz}, we have that
\begin{itemize}
	\item 
	$
	\bar{\mathbf{a}}
	=
	2\beta\mathbf{e}_2
	$
	\ , \ 
	$
	f(\bar{\mathbf{a}})
	=1
	$,
	\item
	$
	\overline{\mathbf{\Delta}}f_{(\mathbf{a},\mathbf{b},\mathbf{c})}
	=
	\mathbf{0}
	$.
\end{itemize}
So, $\mathbf{r}_{f_{(\mathbf{a},\mathbf{b},\mathbf{c})}}
=
\mathbf{0}
$, and $\nabla f(\mathbf{a})=\mathbf{0}$, indeed.\\
In the same example we could also choose the same points but in a different order. For example, one could choose 
\begin{center}
$\mathbf{a}=-\alpha\mathbf{e}_1+\beta\mathbf{e}_2$, $\mathbf{b}=\alpha\mathbf{e}_1+\beta\mathbf{e}_2$, $\mathbf{c}=\mathbf{0}$.
\end{center}
In this case, the mirrored point would be different, and we would have
\begin{itemize}
	\item 
	$
	\mathbf{\Delta}\mathbf{b}
	=
	\mathbf{b}-\mathbf{a}
	=
	2\alpha \mathbf{e}_1
	$
	\ , \ 
	$
	\mathbf{\Delta}\mathbf{c}=\mathbf{c}-\mathbf{a}
	=\alpha \mathbf{e}_1-\beta\mathbf{e}_2
	$
	\ , \ 
	$
	\mathbf{c}-\mathbf{b}
	=
	-\mathbf{b}
	$
	\item
	$
	\Delta_{(\mathbf{a},\mathbf{b},\mathbf{c})}
	=
	(\mathbf{b}-\mathbf{a})\wedge(\mathbf{c}-\mathbf{a})
	=
	-2\alpha\beta\mathbf{e}_1\mathbf{e}_2
	$
	\item
	$
	f(\mathbf{b})-f(\mathbf{a})
	=
	0
	$
	\ , \
	$
	f(\mathbf{c})-f(\mathbf{a})
	=
	1-\sqrt{1-\alpha^2}
	$
	\item
	$
	\mathbf{\Delta} f_{(\mathbf{a},\mathbf{b},\mathbf{c})}
	=
	\big[f(\mathbf{b})-f(\mathbf{a})\big]
	(\mathbf{c}-\mathbf{a})
	-
	\big[f(\mathbf{c})-f(\mathbf{a})\big]
	(\mathbf{b}-\mathbf{a})
	=
	2\alpha (\sqrt{1-\alpha^2}-1)\mathbf{e}_1
	$
	\item
	$
	\displaystyle
	\mathbf{r}_{f_{(\mathbf{a},\mathbf{b},\mathbf{c})}}
	=
	\frac{1-\sqrt{1-\alpha^2}}{2\beta}
	\mathbf{e}_2
	$,
\end{itemize}
as one would expect by Remark~\ref{rem:symmetry of the ratio n=2} and Remark~\ref{rem:multi diff quotient n=2}; besides,
\begin{itemize}
	\item 
	$
	\displaystyle
	\bar{\mathbf{a}}
	=
	2\mathbf{b}
	-
	2\big[(\mathbf{c}-\mathbf{b})\cdot(\mathbf{\Delta}\mathbf{b})\big](		\mathbf{c}-\mathbf{b})^{-1}-\mathbf{a}
	=
	\frac{3\beta^2-\alpha^2}{\alpha^2+\beta^2}
	\alpha \mathbf{e}_1
	+
	\frac{\beta^2-3\alpha^2}{\alpha^2+\beta^2}
	\beta \mathbf{e}_2
	$,
	\item
	$
	\displaystyle
	\bar{\mathbf{a}}-\mathbf{a}
	=
	\frac{4\alpha\beta^2}{\alpha^2+\beta^2}
	\mathbf{e}_1
	-
	\frac{4\alpha^2\beta}{\alpha^2+\beta^2}
	\mathbf{e}_2
	$
	\ \ , \ \
	$
	\displaystyle
	|\bar{\mathbf{a}}-\mathbf{a}|^2
	=
	\frac{16\alpha^2\beta^2}{\alpha^2+\beta^2}
	$
	\item 
	$
	\displaystyle
	f(\bar{\mathbf{a}})
	=
	\sqrt{
	1
	-
	\alpha^2
	\left(
	\frac{3\beta^2-\alpha^2}{\alpha^2+\beta^2}
	\right)^2}
	$.
\end{itemize}
So, we obtain
	\begin{align*}
	&
	\bar{\mathbf{r}}_{f_{(\mathbf{a},\mathbf{b},\mathbf{c})}}
	=
	\frac{f(\bar{\mathbf{a}})-f(\mathbf{a})}{|\bar{\mathbf{a}}-\mathbf{a}|^2}
	(\bar{\mathbf{a}}-\mathbf{a})
	+
	\frac{f(\mathbf{c})-f(\mathbf{b})}{|\mathbf{c}-\mathbf{b}|^2}
	(\mathbf{c}-\mathbf{b})=\\
= & 
	\left\{
	\alpha\beta^2
	\frac{\alpha^2-\beta^2}{(\alpha^2+\beta^2)^2}
	\frac{2}
	{\sqrt{1-\alpha^2
	\left(
	\frac{3\beta^2-\alpha^2}{\alpha^2+\beta^2}
	\right)^2}
	+\sqrt{1-\alpha^2}}
	-
	\frac{\alpha^3}{(\alpha^2+\beta^2)\big[1+\sqrt{1-\alpha^2}\big]}
	\right\}
	\mathbf{e}_1+\\
-	& 
	\left\{
	\alpha^2\beta
	\frac{\alpha^2-\beta^2}{(\alpha^2+\beta^2)^2}
	\frac{2}
	{\sqrt{1-\alpha^2
	\left(
	\frac{3\beta^2-\alpha^2}{\alpha^2+\beta^2}
	\right)^2}
	+\sqrt{1-\alpha^2}}
	+
	\frac{\alpha^2\beta}{(\alpha^2+\beta^2)\big[1+\sqrt{1-\alpha^2}\big]}
	\right\}
	\mathbf{e}_2\\
= & 
	\bar{\mathbf{r}}_{(\alpha,\beta)}\ ,
\end{align*}
and you can verify that
\[
\lim_{(\alpha,\beta)\to(0,0)}
\bar{\mathbf{r}}_{(\alpha,\beta)}
=
\mathbf{0}=\nabla f(\mathbf{0})\ .
\]
\end{Rem*}

\section{The case of a three-variable function}

\subsection{The $3\times 3$ determinant as a Clifford quotient and as a scalar product}
Let us now recall that the determinant of a $3\times 3$ real matrix
\arraycolsep=2.5pt\def\arraystretch{1.}
\[
\left( 
\begin{array}{ccc}
	\mu_{1,1} & \mu_{1,2} & \mu_{1,3} \\
	\mu_{2,1} & \mu_{2,2} & \mu_{2,3} \\
	\mu_{3,1} & \mu_{3,2} & \mu_{3,3} 
\end{array}
\right)
\]
can be written as a coordinate-free Clifford quotient. As usual, let us fix any ordered triple $\mathbf{e}_1,\mathbf{e}_2,\mathbf{e}_3$ of mutually orthonormal vectors in $\mathbb{E}_n$ (with $n\ge 3$), and let us consider $\mathbf{u}_1=\mu_{1,1} \mathbf{e}_1+\mu_{1,2}\mathbf{e}_2+\mu_{1,3}\mathbf{e}_3$, $\mathbf{u}_2=\mu_{2,1} \mathbf{e}_1+\mu_{2,2}\mathbf{e}_2+\mu_{2,3}\mathbf{e}_3$, and $\mathbf{u}_3=\mu_{3,1} \mathbf{e}_1+\mu_{3,2}\mathbf{e}_2+\mu_{3,3}\mathbf{e}_3$, then you can verify that
\[
(\mathbf{u}_1\wedge\mathbf{u}_2\wedge\mathbf{u}_3)(I_3)^{-1}
 =
\det
\left( 
\begin{array}{ccc}
	\mu_{1,1} & \mu_{1,2} & \mu_{1,3} \\
	\mu_{2,1} & \mu_{2,2} & \mu_{2,3} \\
	\mu_{3,1} & \mu_{3,2} & \mu_{3,3} 
\end{array}
\right)\ ,
\]
because 
\[
\mathbf{u}_{\sigma_1}\wedge\mathbf{u}_{\sigma_2}\wedge\mathbf{u}_{\sigma_3}
=
\epsilon_\sigma\
\mathbf{u}_1\wedge\mathbf{u}_2\wedge\mathbf{u}_3\ ,
\]
for each permutation $\sigma\in\mathcal{S}_3$ of the set $\{1,2,3\}$, having parity $\epsilon_\sigma\in\{-1,1\}$; where
\begin{align*}
\mathbf{u}_1\wedge\mathbf{u}_2\wedge\mathbf{u}_3 & =
\frac{1}{6}
\big(
	\mathbf{u}_1\mathbf{u}_2\mathbf{u}_3
	-\mathbf{u}_1\mathbf{u}_3\mathbf{u}_2
	+\mathbf{u}_3\mathbf{u}_1\mathbf{u}_2
	-\mathbf{u}_3\mathbf{u}_2\mathbf{u}_1
	+\mathbf{u}_2\mathbf{u}_3\mathbf{u}_1
	-\mathbf{u}_2\mathbf{u}_1\mathbf{u}_3
\big)\\
& =
\frac{1}{2}
\Big[
(\mathbf{u}_1\wedge\mathbf{u}_2)\mathbf{u}_3
+
\mathbf{u}_3(\mathbf{u}_1\wedge\mathbf{u}_2)
\Big]\\
& =
\frac{1}{2}
(\mathbf{u}_1\mathbf{u}_2\mathbf{u}_3-\mathbf{u}_3\mathbf{u}_2\mathbf{u}_1)\ ,\\ \\
\textrm{and }\ I_3 & = \mathbf{e}_1\mathbf{e}_2\mathbf{e}_3=\mathbf{e}_1\wedge \mathbf{e}_2\wedge \mathbf{e}_3\ , \ \textrm{ so that }\  (I_3)^{-1}= \mathbf{e}_3\mathbf{e}_2\mathbf{e}_1=-I_3\ .
\end{align*}
\begin{Rem*}
As $I_2$ before, also $I_3$ does not depend on the particular orthonormal basis chosen to define it in $span \{\mathbf{e}_1,\mathbf{e}_2,\mathbf{e}_3\}\subseteq \mathbb{E}_n$, but only on the orientation of that basis. More precisely, if $\{\mathbf{g}_1,\mathbf{g}_2, \mathbf{g}_3\}$ is any other orthonormal basis of $span \{\mathbf{e}_1,\mathbf{e}_2,\mathbf{e}_3\}$, then $\mathbf{g}_1\wedge\mathbf{g}_2\wedge \mathbf{g}_3=\mathbf{g}_1\mathbf{g}_2\mathbf{g}_3$ is equal to $I_3$ or $-I_3$. That is why $I_3$ is called an ``orientation'' of~$span \{\mathbf{e}_1,\mathbf{e}_2,\mathbf{e}_3\}=\mathbb{E}_3$.
\end{Rem*}
Thus, a $3\times 3$ determinant can be considered as the Clifford  ratio between the two ``$3$-blades'' $\mathbf{u}_1\wedge\mathbf{u}_2\wedge\mathbf{u}_3$ and $I_3$ (a ``$k$-blade'' being the geometric product of $k$ non zero and mutually orthogonal vectors). Those elements can also be called ``$\mathbb{G}_3$-pseudo-scalars'', as $\mathbb{G}_3$ is generated by~$\mathbb{E}_3=span \{\mathbf{e}_1,\mathbf{e}_2,\mathbf{e}_3\}$, and can be interpreted as oriented volumes in $\mathbb{E}_3$.
Let us observe that, if 
\begin{center}
$V=\nu_{1,2}\mathbf{e}_1\mathbf{e}_2+\nu_{1,3}\mathbf{e}_1\mathbf{e}_3+\nu_{2,3}\mathbf{e}_2\mathbf{e}_3\in span\{\mathbf{e}_1\mathbf{e}_2\ ,\ \mathbf{e}_1\mathbf{e}_3\ , \ \mathbf{e}_2\mathbf{e}_3\}=\mathbb{G}_{3 \choose 2}$,
\end{center}
(such elements in $\mathbb{G}_3$ are also called the ``$2$-vectors''), then 
\begin{align*}
V I_3 
= & 
(\nu_{1,2}\mathbf{e}_1\mathbf{e}_2+\nu_{1,3}\mathbf{e}_1\mathbf{e}_3+\nu_{2,3}\mathbf{e}_2\mathbf{e}_3)\mathbf{e}_1\mathbf{e}_2\mathbf{e}_3\\
= & 
\nu_{1,2}\mathbf{e}_1\mathbf{e}_2\mathbf{e}_1\mathbf{e}_2\mathbf{e}_3
+\nu_{1,3}\mathbf{e}_1\mathbf{e}_3\mathbf{e}_1\mathbf{e}_2\mathbf{e}_3
+\nu_{2,3}\mathbf{e}_2\mathbf{e}_3\mathbf{e}_1\mathbf{e}_2\mathbf{e}_3\\
= & 
\nu_{1,2}\mathbf{e}_1\mathbf{e}_2\mathbf{e}_3\mathbf{e}_1\mathbf{e}_2
+\nu_{1,3}\mathbf{e}_1\mathbf{e}_2\mathbf{e}_3\mathbf{e}_1\mathbf{e}_3
+\nu_{2,3}\mathbf{e}_1\mathbf{e}_2\mathbf{e}_3\mathbf{e}_2\mathbf{e}_3\\
= &
\mathbf{e}_1\mathbf{e}_2\mathbf{e}_3
(\nu_{1,2}\mathbf{e}_1\mathbf{e}_2+\nu_{1,3}\mathbf{e}_1\mathbf{e}_3+\nu_{2,3}\mathbf{e}_2\mathbf{e}_3)
=
I_3 V\\
= &
-\nu_{1,2}\mathbf{e}_3+\nu_{1,3}\mathbf{e}_2-\nu_{2,3}\mathbf{e}_1\in\mathbb{E}_3\ ,
\end{align*} 
that is why the elements in $\mathbb{G}_{3\choose 2}$ are also called ``$\mathbb{G}_3$-pseudo-vectors'': geometric multiplication by $I_3$ establishes a duality between vectors of $\mathbb{E}_3$ and elements of $\mathbb{G}_{3\choose 2}$. In a similar way, you can verify that $\mathbf{u}I_3=I_3\mathbf{u}$ for all $\mathbf{u}\in\mathbb{E}_3$. The foregoing properties allows us to write a $3\times 3$ determinant as a scalar product
\begin{align*}
&
\det
\left( 
\begin{array}{ccc}
	\mu_{1,1} & \mu_{1,2} & \mu_{1,3} \\
	\mu_{2,1} & \mu_{2,2} & \mu_{2,3} \\
	\mu_{3,1} & \mu_{3,2} & \mu_{3,3} 
\end{array}
\right)
= 
(\mathbf{u}_1\wedge\mathbf{u}_2\wedge\mathbf{u}_3)(I_3)^{-1}=\\
= &
-\frac{1}{2}\Big[
\big(\mathbf{u}_1\wedge\mathbf{u}_2\big)\mathbf{u}_3+\mathbf{u}_3\big(\mathbf{u}_1\wedge\mathbf{u}_2\big)\Big]I_3
= 
-\frac{1}{2}\Big[
\big(\mathbf{u}_1\wedge\mathbf{u}_2\big)I_3\mathbf{u}_3+\mathbf{u}_3\big(\mathbf{u}_1\wedge\mathbf{u}_2\big)I_3\Big]=\\
= &
-\big[(\mathbf{u}_1\wedge\mathbf{u}_2)I_3\big]\cdot\mathbf{u}_3
=
\big[(\mathbf{u}_1\wedge\mathbf{u}_2)(I_3)^{-1}\big]\cdot\mathbf{u}_3\ ,
\end{align*}
where $(\mathbf{u}_1\wedge\mathbf{u}_2)(I_3)^{-1}$ is a vector orthogonal to $span\{\mathbf{u}_1,\mathbf{u}_2\}$, because it is orthogonal to both $\mathbf{u}_1$ and $\mathbf{u}_2$, as
\begin{center}
$
\big[(\mathbf{u}_1\wedge\mathbf{u}_2)(I_3)^{-1}\big]\cdot\mathbf{u}_1
=
(\mathbf{u}_1\wedge\mathbf{u}_2\wedge\mathbf{u}_1)(I_3)^{-1}
=
0
$
and 
$
\big[(\mathbf{u}_1\wedge\mathbf{u}_2)(I_3)^{-1}\big]\cdot\mathbf{u}_2
=
(\mathbf{u}_1\wedge\mathbf{u}_2\wedge\mathbf{u}_2)(I_3)^{-1}
=
0
$.
\end{center}
As a matter of fact, $(\mathbf{u}_1\wedge\mathbf{u}_2)(I_3)^{-1}$ corresponds to the classical Gibbs and Heaviside cross product of $\mathbf{u}_1$ and $\mathbf{u}_2$, when ~$\mathbb{E}_3=span \{\mathbf{e}_1,\mathbf{e}_2,\mathbf{e}_3\}$ is identified with~$\mathbb{R}^3$.

\subsection{Coordinate-free expression of a hyperplane secant the graph of a three-variable function}
Let us write the equation of a hyperplane secant the graph of a three-variable function $f:\Omega\subseteq\mathbb{E}_3\to\mathbb{R}$ at four non coplanar points of that graph.
A hyperplane passing through the four non coplanar points 
\[
\big(\mathbf{a}_i,f(\mathbf{a}_i)\big)
=
\big(\alpha_{i,1}\mathbf{e}_1+\alpha_{i,2}\mathbf{e}_2+\alpha_{i,3}\mathbf{e}_3,f(\mathbf{a}_i)\big)\in \mathbb{E}_3\oplus\mathbb{R}\ ,
\]
with $i=1,2,3,4$, can be represented by the Cartesian relation
\begin{equation}
\label{eq:secant hyperplane}
\det
\left(
\begin{array}{cccc}
	\chi_1-\alpha_{1,1} & \chi_2-\alpha_{1,2} & \chi_3-\alpha_{1,3} & z-f(\mathbf{a}_1)\\
	\alpha_{2,1}-\alpha_{1,1} & \alpha_{2,2}-\alpha_{1,2} & \alpha_{2,3}-\alpha_{1,3} & f(\mathbf{a}_2)-f(\mathbf{a}_1)\\
	\alpha_{3,1}-\alpha_{1,1} & \alpha_{3,2}-\alpha_{1,2} & \alpha_{3,3}-\alpha_{1,3} & f(\mathbf{a}_3)-f(\mathbf{a}_1)\\
	\alpha_{4,1}-\alpha_{1,1} & \alpha_{4,2}-\alpha_{1,2} & \alpha_{4,3}-\alpha_{1,3} & f(\mathbf{a}_4)-f(\mathbf{a}_1)
\end{array}
\right)
=0
\end{equation}
between the real variables $\chi_1,\chi_2,\chi_3,z\in\mathbb{R}$. That determinant can be rewritten by a  Laplace expansion as follows
\begin{align*}
&
\big[z-f(\mathbf{a}_1)\big]
\det
\left(
\begin{array}{ccc}
\alpha_{2,1}-\alpha_{1,1} & \alpha_{2,2}-\alpha_{1,2} & \alpha_{2,3}-\alpha_{1,3} \\
	\alpha_{3,1}-\alpha_{1,1} & \alpha_{3,2}-\alpha_{1,2} & \alpha_{3,3}-\alpha_{1,3} \\
	\alpha_{4,1}-\alpha_{1,1} & \alpha_{4,2}-\alpha_{1,2} & \alpha_{4,3}-\alpha_{1,3} 
\end{array}
\right)
+\\
& -
\big[f(\mathbf{a}_2)-f(\mathbf{a}_1)\big]
\det
\left(
\begin{array}{ccc}
	\chi_1-\alpha_{1,1} & \chi_2-\alpha_{1,2} & \chi_3-\alpha_{1,3}\\
	\alpha_{3,1}-\alpha_{1,1} & \alpha_{3,2}-\alpha_{1,2} & \alpha_{3,3}-\alpha_{1,3} \\
	\alpha_{4,1}-\alpha_{1,1} & \alpha_{4,2}-\alpha_{1,2} & \alpha_{4,3}-\alpha_{1,3} 
\end{array}
\right)
+\\
& +
\big[f(\mathbf{a}_3)-f(\mathbf{a}_1)\big]
\det
\left(
\begin{array}{ccc}
	\chi_1-\alpha_{1,1} & \chi_2-\alpha_{1,2} & \chi_3-\alpha_{1,3} \\
	\alpha_{2,1}-\alpha_{1,1} & \alpha_{2,2}-\alpha_{1,2} & \alpha_{2,3}-\alpha_{1,3} \\
\alpha_{4,1}-\alpha_{1,1} & \alpha_{4,2}-\alpha_{1,2} & \alpha_{4,3}-\alpha_{1,3} 
\end{array}
\right)
+\\
& -
\big[f(\mathbf{a}_4)-f(\mathbf{a}_1)\big]
\det
\left(
\begin{array}{ccc}
	\chi_1-\alpha_{1,1} & \chi_2-\alpha_{1,2} & \chi_3-\alpha_{1,3} \\
	\alpha_{2,1}-\alpha_{1,1} & \alpha_{2,2}-\alpha_{1,2} & \alpha_{2,3}-\alpha_{1,3} \\
	\alpha_{3,1}-\alpha_{1,1} & \alpha_{3,2}-\alpha_{1,2} & \alpha_{3,3}-\alpha_{1,3} 
\end{array}
\right)
\end{align*}
Then, in $\mathbb{G}_3$ the equation~(\ref{eq:secant hyperplane}) becomes 
\begin{align*}
&
\big[z-f(\mathbf{a}_1)\big]
\big[(\mathbf{a}_2-\mathbf{a}_1)\wedge(\mathbf{a}_3-\mathbf{a}_1)\wedge(\mathbf{a}_4-\mathbf{a}_1)\big]
(I_3)^{-1}+\\
- &
\big[f(\mathbf{a}_2)-f(\mathbf{a}_1)\big]
\big[(\mathbf{x}-\mathbf{a}_1)\wedge(\mathbf{a}_3-\mathbf{a}_1)\wedge(\mathbf{a}_4-\mathbf{a}_1)\big]
(I_3)^{-1}+\\
+ &
\big[f(\mathbf{a}_3)-f(\mathbf{a}_1)\big]
\big[(\mathbf{x}-\mathbf{a}_1)\wedge(\mathbf{a}_2-\mathbf{a}_1)\wedge(\mathbf{a}_4-\mathbf{a}_1)\big]
(I_3)^{-1}+\\
- &
\big[f(\mathbf{a}_4)-f(\mathbf{a}_1)\big]
\big[(\mathbf{x}-\mathbf{a}_1)\wedge(\mathbf{a}_2-\mathbf{a}_1)\wedge(\mathbf{a}_3-\mathbf{a}_1)\big]
(I_3)^{-1}
=
0\ ,
\end{align*}
being $\mathbf{x}=\chi_1\mathbf{e}_1+\chi_2\mathbf{e}_2+\chi_3\mathbf{e}_3\in\mathbb{E}_3$, $(\mathbf{x},z)\in\mathbb{E}_3\oplus\mathbb{R}$. The foregoing relation is equivalent, in $\mathbb{G}_3$, to
\begin{align*}
&
\big[z-f(\mathbf{a}_1)\big]
\big[(\mathbf{a}_2-\mathbf{a}_1)\wedge(\mathbf{a}_3-\mathbf{a}_1)\wedge(\mathbf{a}_4-\mathbf{a}_1)\big]
=\\
= &
\big[f(\mathbf{a}_2)-f(\mathbf{a}_1)\big]
\big[(\mathbf{x}-\mathbf{a}_1)\wedge(\mathbf{a}_3-\mathbf{a}_1)\wedge(\mathbf{a}_4-\mathbf{a}_1)\big]
+\\
& -
\big[f(\mathbf{a}_3)-f(\mathbf{a}_1)\big]
\big[(\mathbf{x}-\mathbf{a}_1)\wedge(\mathbf{a}_2-\mathbf{a}_1)\wedge(\mathbf{a}_4-\mathbf{a}_1)\big]
+\\
& +
\big[f(\mathbf{a}_4)-f(\mathbf{a}_1)\big]
\big[(\mathbf{x}-\mathbf{a}_1)\wedge(\mathbf{a}_2-\mathbf{a}_1)\wedge(\mathbf{a}_3-\mathbf{a}_1)\big]\ ,
\end{align*}
Let us define 
\begin{center}
$
\displaystyle
\Delta_{(\mathbf{a}_1,\mathbf{a}_2,\mathbf{a}_3,\mathbf{a}_4)}
=
(\mathbf{a}_2-\mathbf{a}_1)\wedge(\mathbf{a}_3-\mathbf{a}_1)\wedge(\mathbf{a}_4-\mathbf{a}_1)
\Big(
=
-
(\mathbf{a}_1-\mathbf{a}_2)
\wedge
(\mathbf{a}_2-\mathbf{a}_3)
\wedge
(\mathbf{a}_3-\mathbf{a}_4)
\Big)
$.
\end{center}
We observe that 
\begin{center}
$
\displaystyle \tau_3 =\frac{1}{6}
\det
\left(
\begin{array}{ccc}
\alpha_{2,1}-\alpha_{1,1} & \alpha_{2,2}-\alpha_{1,2} & \alpha_{2,3}-\alpha_{1,3} \\
	\alpha_{3,1}-\alpha_{1,1} & \alpha_{3,2}-\alpha_{1,2} & \alpha_{3,3}-\alpha_{1,3} \\
	\alpha_{4,1}-\alpha_{1,1} & \alpha_{4,2}-\alpha_{1,2} & \alpha_{4,3}-\alpha_{1,3}\end{array}
\right)
=
\frac{1}{6}
\Delta_{(\mathbf{a}_1,\mathbf{a}_2,\mathbf{a}_3,\mathbf{a}_4)}
(I_3)^{-1}
$
\end{center}
is the oriented volume of the tetrahedron having vertices $\mathbf{a}_1$, $\mathbf{a}_2$, $\mathbf{a}_3$, and $\mathbf{a}_4$. 
So, we can write 
$
\Delta_{(\mathbf{a}_1,\mathbf{a}_2,\mathbf{a}_3,\mathbf{a}_4)}
=
6\tau_3 I_3
$, and $\big(\Delta_{(\mathbf{a}_1,\mathbf{a}_2,\mathbf{a}_3,\mathbf{a}_4)}\big)^{-1}= \frac{1}{6\tau}(I_3)^{-1}$. By denoting $\Delta f_{(\mathbf{a}_i)} = f(\mathbf{a}_i)-f(\mathbf{a}_1)$ and $\mathbf{\Delta}\mathbf{a}_i=\mathbf{a}_i-\mathbf{a}_1$, when $i=2,3,4$, the equation of the secant hyperplane~(\ref{eq:secant hyperplane}) becomes
\begin{align*}
& z =\\ 
= & 
f(\mathbf{a}_1)
+ \frac{1}{6\tau_3}
\Big\{
\Delta f_{(\mathbf{a}_2)}
(\mathbf{x}-\mathbf{a}_1)\wedge 
\mathbf{\Delta}\mathbf{a}_3\wedge \mathbf{\Delta}\mathbf{a}_4
-
\Delta f_{(\mathbf{a}_3)}
(\mathbf{x}-\mathbf{a}_1)\wedge 
\mathbf{\Delta}\mathbf{a}_2\wedge \mathbf{\Delta}\mathbf{a}_4
+\\
& +
\Delta f_{(\mathbf{a}_4)}
(\mathbf{x}-\mathbf{a}_1)\wedge 
\mathbf{\Delta}\mathbf{a}_2\wedge \mathbf{\Delta}\mathbf{a}_3\Big\}(I_3)^{-1}\\
= & 
f(\mathbf{a}_1)
+ \frac{1}{6\tau_3}
\Big\{
\Delta f_{(\mathbf{a}_2)}
\mathbf{\Delta}\mathbf{a}_3\wedge \mathbf{\Delta}\mathbf{a}_4 \wedge (\mathbf{x}-\mathbf{a}_1)
-
\Delta f_{(\mathbf{a}_3)}
\mathbf{\Delta}\mathbf{a}_2\wedge \mathbf{\Delta}\mathbf{a}_4\wedge(\mathbf{x}-\mathbf{a}_1)
+\\
& +
\Delta f_{(\mathbf{a}_4)}
\mathbf{\Delta}\mathbf{a}_2\wedge \mathbf{\Delta}\mathbf{a}_3\wedge (\mathbf{x}-\mathbf{a}_1)
\Big\}(I_3)^{-1}\\
= & 
f(\mathbf{a}_1)
+ \frac{1}{6\tau_3}
\Big\{
\Delta f_{(\mathbf{a}_2)}
(\mathbf{\Delta}\mathbf{a}_3\wedge \mathbf{\Delta}\mathbf{a}_4)(I_3)^{-1}
-
\Delta f_{(\mathbf{a}_3)}
(\mathbf{\Delta}\mathbf{a}_2\wedge \mathbf{\Delta}\mathbf{a}_4)(I_3)^{-1}+\\
& +
\Delta f_{(\mathbf{a}_4)}
(\mathbf{\Delta}\mathbf{a}_2\wedge \mathbf{\Delta}\mathbf{a}_3)(I_3)^{-1}\Big\}\cdot (\mathbf{x}-\mathbf{a}_1) \\
= & 
f(\mathbf{a}_1)+\\ 
+ &
\frac{1}{6\tau_3}
\Big\{\Big[
\Delta f_{(\mathbf{a}_2)}
\mathbf{\Delta}\mathbf{a}_3\wedge \mathbf{\Delta}\mathbf{a}_4
-
\Delta f_{(\mathbf{a}_3)}
\mathbf{\Delta}\mathbf{a}_2\wedge \mathbf{\Delta}\mathbf{a}_4
+
\Delta f_{(\mathbf{a}_4)}
\mathbf{\Delta}\mathbf{a}_2\wedge \mathbf{\Delta}\mathbf{a}_3\Big](I_3)^{-1}\Big\}
\kern-3pt \cdot\kern-3pt
(\mathbf{x}-\mathbf{a}_1) 
\end{align*} 
Thus, we have that the vector $\mathbf{r}_{f_{(\mathbf{a}_1,\mathbf{a}_2,\mathbf{a}_3,\mathbf{a}_4)}}\in\mathbb{E}_3$, 
characterizing the equation of the hyperplane secant the graph of $f$ as 
$
z
=
f(\mathbf{a}_1)+\mathbf{r}_{f_{(\mathbf{a}_1,\mathbf{a}_2,\mathbf{a}_3,\mathbf{a}_4)}}\cdot (\mathbf{x}-\mathbf{a}_1)
$,
is
\[
\mathbf{r}_{f_{(\mathbf{a}_1,\mathbf{a}_2,\mathbf{a}_3,\mathbf{a}_4)}}
=
\left[
\sum_{i=2}^4
(-1)^i
\Delta f_{(\mathbf{a}_i)}
\Delta^i_{(\mathbf{a}_1,\mathbf{a}_2,\mathbf{a}_3,\mathbf{a}_4)}
\right]
\big(\Delta_{(\mathbf{a}_1,\mathbf{a}_2,\mathbf{a}_3,\mathbf{a}_4)}\big)^{-1}\ ,
\]
(where $\Delta^2_{(\mathbf{a}_1,\mathbf{a}_2,\mathbf{a}_3,\mathbf{a}_4)}=\mathbf{\Delta}\mathbf{a}_3\wedge \mathbf{\Delta}\mathbf{a}_4
$, $\Delta^3_{(\mathbf{a}_1,\mathbf{a}_2,\mathbf{a}_3,\mathbf{a}_4)}=\mathbf{\Delta}\mathbf{a}_2\wedge \mathbf{\Delta}\mathbf{a}_4
$, $\Delta^4_{(\mathbf{a}_1,\mathbf{a}_2,\mathbf{a}_3,\mathbf{a}_4)}=\mathbf{\Delta}\mathbf{a}_2\wedge \mathbf{\Delta}\mathbf{a}_3
$), which is, in fact, the Clifford quotient between the multi-difference $\mathbb{G}_3$-pseudo-vector 
\[
\Delta f_{(\mathbf{a}_1,\mathbf{a}_2,\mathbf{a}_3,\mathbf{a}_4)}
=
\Delta f_{(\mathbf{a}_2)}
\mathbf{\Delta}\mathbf{a}_3\wedge \mathbf{\Delta}\mathbf{a}_4
-
\Delta f_{(\mathbf{a}_3)}
\mathbf{\Delta}\mathbf{a}_2\wedge \mathbf{\Delta}\mathbf{a}_4
+
\Delta f_{(\mathbf{a}_4)}
\mathbf{\Delta}\mathbf{a}_2\wedge \mathbf{\Delta}\mathbf{a}_3\in\mathbb{G}_{3\choose 2}\ ,
\]
and the $\mathbb{G}_3$-pseudo-scalar
\[
\Delta_{(\mathbf{a}_1,\mathbf{a}_2,\mathbf{a}_3,\mathbf{a}_4)}
=
\mathbf{\Delta}\mathbf{a}_2
\wedge
\mathbf{\Delta}\mathbf{a}_3
\wedge
\mathbf{\Delta}\mathbf{a}_4
=
6\tau_3 I_3\in\mathbb{G}_{3\choose 3}\simeq \mathbb{R}I_3\ .
\]

\subsection{Mirroring vectors and points II}
\subsubsection{Vector mirrored by a $2$-dimensional linear subspace}
We recall that to each pair of linearly independent vectors $\mathbf{u}_1$, $\mathbf{u}_2 \in\mathbb{E}_n$ we can associate the $2$-dimensional linear subspace $span\{\mathbf{u}_1,\mathbf{u}_2\}\subseteq \mathbb{E}_n$.
Moreover, $\mathbf{u}_1\wedge\mathbf{u}_2$ is always a $2$-blade. As a matter of fact, there always exists an orthogonal basis $\{\mathbf{g}_1,\mathbf{g}_2\}$ of  $span\{\mathbf{u}_1,\mathbf{u}_2\}$, and you can verify that $\mathbf{u}_1\wedge\mathbf{u}_2$ is a non zero multiple of the geometric product $\mathbf{g}_1\mathbf{g}_2$.
We recall that 
\begin{itemize}
	\item the square of every $2$-blade is a non zero scalar, 
	\item every $2$-blade~$B$ is invertible in $\mathbb{G}_n$, and $\displaystyle B^{-1}=\frac{1}{B^2}B$.
\end{itemize}
Let us recall that, if $\mathbf{v}\in\mathbb{E}_n$, and $B=\mathbf{u}_1\wedge\mathbf{u}_2$ is a $2$-blade, then
\begin{equation}
\label{eq:2-blade vector decomposition}
\mathbf{v}
=
\mathbf{v}B B^{-1}
=
(\mathbf{v}B)B^{-1}
=
(\mathbf{v}\circ B + \mathbf{v}\Wedge B)B^{-1}
=
(\mathbf{v}\circ B)B^{-1}+(\mathbf{v}\Wedge B)B^{-1}
\end{equation}
where
\begin{align*}
\mathbf{v} \circ B
& =
\frac{1}{2}\big(\mathbf{v}B - B\mathbf{v}\big)
=
- B \circ\mathbf{v}
\in\mathbb{G}_{n\choose 1}\\
\mathbf{v} \Wedge B
& =
\frac{1}{2}\big(\mathbf{v}B+ B\mathbf{v}\big)
=
B \Wedge\mathbf{v}
\in\mathbb{G}_{n\choose 3}\\
\mathbf{v} \Wedge (\mathbf{u}_1\wedge\mathbf{u}_2)
& =
\mathbf{v} \wedge \mathbf{u}_1\wedge\mathbf{u}_2\\
\mathbb{G}_{n \choose k}
& =
\left\{
\begin{array}{ll}
	\mathbb{R} & \textrm{ if } k=0\\
	span\{\mathbf{e}_{i_1}\cdots\mathbf{e}_{i_k}\}_{1\le i_1<\cdots<i_k\le n} & \textrm{ if } 1\le k\le n\\
\end{array}
\right.
\end{align*}
\begin{Rem*}
Elements of $\mathbb{G}_{n\choose k}$ are called ``$k$-vectors''. Notice that the dimension of~$\mathbb{G}_{n\choose k}$ is the binomial coefficient$\displaystyle {n\choose k} = \frac{n!}{(n-k)!\ k!}$. Moreover $\displaystyle \mathbb{G}_n=\bigoplus_{k=0}^n \mathbb{G}_{n\choose k}$.
\end{Rem*}
\begin{Rem}\label{rem:comm and antcomm vector-blade}
We recall that the foregoing operations ``$\circ$'' and ``$\Wedge$'' can be extended\footnote{See for example~\cite{Wilmot} or \cite{Delaby}.} to $k$-blades 
$
\displaystyle B
=
\mathbf{u}_1\wedge \cdots\wedge\mathbf{u}_k=\bigwedge_{j=1}^k
\mathbf{u}_j
\in\mathbb{G}_{n\choose k}$ where
\[
\mathbf{u}_1\wedge \cdots\wedge\mathbf{u}_k
=
\frac{1}{k!}
\sum_{\sigma\in\mathcal{S}_k}
\epsilon_\sigma 
\mathbf{u}_{\sigma_1}\cdots \mathbf{u}_{\sigma_k}
\]
\big($\mathcal{S}_k$ being the group of all permutations of $\{1,\dots,k\}$, $\epsilon_\sigma\in\{-1,1\}$ the parity of permutation $\sigma$\big),
as follows
\begin{align*}
\mathbf{v} \circ B
& =
\frac{1}{2}\big(\mathbf{v}B-(-1)^k B\mathbf{v}\big)
=
(-1)^{k+1} B \circ\mathbf{v}
\in\mathbb{G}_{n\choose k-1}\ ,\\
\mathbf{v} \Wedge B
& =
\frac{1}{2}\big(\mathbf{v}B+(-1)^k B\mathbf{v}\big)
=
(-1)^{k} B \Wedge\mathbf{v}
\in\mathbb{G}_{n\choose k+1}\ .
\end{align*}
We also recall that 
\begin{align*}
\mathbf{v}\circ (\mathbf{u}_1\wedge \cdots\wedge\mathbf{u}_k)
& =
\sum_{i=1}^k (-1)^{i+1}(\mathbf{v}\cdot \mathbf{u}_i) 
\bigwedge_{
\begin{array}{c}
\scriptstyle j=1 \\
\scriptstyle j\ne i
\end{array}}^k
\mathbf{u}_j\\
\mathbf{v}\Wedge (\mathbf{u}_1\wedge \cdots\wedge\mathbf{u}_k)
& = 
\mathbf{v}\wedge \mathbf{u}_1\wedge \cdots\wedge\mathbf{u}_k\ .
\end{align*}
So, we have that
\[
B\mathbf{v}
=
2(B\circ\mathbf{v})+(-1)^k\mathbf{v}B
=
(B\circ \mathbf{v})
+
(B\Wedge \mathbf{v})
\ .
\]
\end{Rem}
\begin{Rem*}
Notice that, if $\mathbf{u}$ and $\mathbf{v}$ are vectors in $\mathbb{E}_n$, then $\mathbf{u}\circ \mathbf{v}=\mathbf{u}\cdot \mathbf{v}$, and $\mathbf{u}\Wedge \mathbf{v}=\mathbf{u}\wedge \mathbf{v}$. 
\end{Rem*}
\begin{Rem*}
The operations ``$\circ$'' and ``$\Wedge$'' are particular cases of more general operations between blades. More precisely, if $H$ is a $h$-blade, and $K$ is a $k$-blade, then
\begin{align*}
H \circ K
& =
\frac{1}{2}\big(HK-(-1)^{hk} KH\big)
=
(-1)^{hk+1} K \circ H 
\ \textrm{ is called ``graded commutator''}\\
H \Wedge K
& =
\frac{1}{2}\big(HK+(-1)^k KH\big)
=
(-1)^{hk} K \Wedge H
\textrm{ is called ``graded anti-commutator''.}
\end{align*}
Such operations can then be extended, by linearity, to linear combinations of blades that is, to every element of the geometric algebra~$\mathbb{G}_n$. 
\end{Rem*}
\begin{Rem}\label{rem:rejection is orthogonal 2}
If $B=\mathbf{u}_1\wedge\mathbf{u}_2$ is a $2$-blade, then $(\mathbf{v}\Wedge B)B^{-1}$ is a vector which is orthogonal both to $\mathbf{u}_1$ and $\mathbf{u}_2$ that is, to all two-dimensional $span\{ \mathbf{u}_1, \mathbf{u}_2\}$. As a matter of fact, for $i=1,2$ we have that 
\[
2\big[(\mathbf{v}\Wedge B)B^{-1}\big]\cdot \mathbf{u}_i
=
(\mathbf{v}\Wedge B)B^{-1} \mathbf{u}_i
+
\mathbf{u}_i(\mathbf{v}\Wedge B)B^{-1}
\]
You can verify that $B\mathbf{u}_i = -B \mathbf{u}_i$ (for $i=1,2$ ). So we have that 
\[
2\big[(\mathbf{v}\Wedge B)B^{-1}\big]\cdot \mathbf{u}_i
=
\big[
-(\mathbf{v}\Wedge B) \mathbf{u}_i 
+
\mathbf{u}_i(\mathbf{v}\Wedge B)
\big]B^{-1}
\]
As $\mathbf{v}\Wedge B= \mathbf{v}\wedge \mathbf{u}_1 \wedge \mathbf{u}_2$ is either zero or a $3$-blade, then we can write
\[
4\big[(\mathbf{v}\Wedge B)B^{-1}\big]\cdot \mathbf{u}_i
=
\big[\mathbf{u}_i\Wedge (\mathbf{v}\Wedge B)
\big]B^{-1}
=
\big(\mathbf{u}_i\wedge \mathbf{v}\wedge \mathbf{u}_1\wedge \mathbf{u}_2
\big)B^{-1}
=
0\ .
\]
\end{Rem}
In a similar way you can verify that for every vector $\mathbf{v}\in\mathbb{E}_n$ we have that 
\begin{center}
$\mathbf{v}_\parallel= (\mathbf{v}\circ B)B^{-1}\in span\{\mathbf{u}_1,\mathbf{u}_2\}$.
\end{center}
So, relation~(\ref{eq:2-blade vector decomposition}) corresponds to the orthogonal  decomposition of vector $\mathbf{v}$ with respect to the two-dimensional linear subspace $span\{\mathbf{u}_1,\mathbf{u}_2\}=\mathcal{L}_{(\mathbf{0},\mathbf{u}_1,\mathbf{u}_2)}$.\\
Then, the vector $\mathbf{\hat{v}}$, mirrored of vector $\mathbf{v}$ by the $2$-dimensional linear subspace $span\{\mathbf{u}_1,\mathbf{u}_2\}$, can be written as
\begin{align*}
\mathbf{\hat{v}}
= &
\mathbf{v}_{\parallel}-\mathbf{v}_{\bot}
=
(\mathbf{v}\circ B)B^{-1}
-
(\mathbf{v}\Wedge B)B^{-1}
=
\big[(\mathbf{v}\circ B)-(\mathbf{v}\Wedge B)\big]B^{-1}\\
= &
\big[-(B\circ\mathbf{v})-(B\Wedge \mathbf{v})\big]B^{-1}
=
-\big[(B\circ\mathbf{v})+(B\Wedge \mathbf{v})\big]B^{-1}\\
= &
-B\mathbf{v}B^{-1}
=
-\Big[2(B\circ\mathbf{v})+\mathbf{v}B\big]B^{-1}
= 
-2(B\circ\mathbf{v})B^{-1}-\mathbf{v}\ .
\end{align*}
Moreover, $|\mathbf{\hat{v}}|= |\mathbf{v}|$; as a matter of fact,
\begin{center}
$
|\mathbf{\hat{v}}|^2
=
\mathbf{\hat{v}}\mathbf{\hat{v}}
=
B\mathbf{v}B^{-1}
B\mathbf{v}B^{-1}
=
B\mathbf{v}\mathbf{v}B^{-1}
=
|\mathbf{v}|^2BB^{-1}
=
|\mathbf{v}|^2
$.
\end{center}
\begin{Rem*}
Notice that, $\mathbf{\hat{v}}$ is in $span\{\mathbf{u}_1,\mathbf{u}_2,\mathbf{v}\}$.
\end{Rem*}
\subsubsection{Point mirrored by a plane in $\mathbb{E}_n$ (with $n\ge 3$)}
Given four non coplanar points $\mathbf{a}_1$, $\mathbf{a}_2$, , $\mathbf{a}_3$ and $\mathbf{a}_4$ in $\mathbb{E}_n$ (with $n\ge 3$), we mirror the point~$\mathbf{a}_1$ through the plane $\mathcal{L}_{(\mathbf{a}_2,\mathbf{a}_3,\mathbf{a}_4)}$ generated by the points $\mathbf{a}_2$, $\mathbf{a}_3$, and $\mathbf{a}_4$. We denote $\overline{\mathbf{a}_1}$ that mirrored point.\\
We can express  $\overline{\mathbf{a}_1}$ by computing in $\mathbb{G}_n$ the vector $\mathbf{\hat{v}}$, obtained by mirroring vector $\mathbf{v}=\mathbf{a}_1-\mathbf{a}_2$ by the $2$-blade 
$
L_1
=
L_{(\mathbf{a}_2,\mathbf{a}_3,\mathbf{a}_4)}
=
(\mathbf{a}_3-\mathbf{a}_2)\wedge(\mathbf{a}_4-\mathbf{a}_2)
$.\\
More precisely, the  mirrored point $\overline{\mathbf{a}_1}$
can be computed by using the geometric Clifford product in $\mathbb{G}_n$ as follows
\begin{align*}
\overline{\mathbf{a}_1}
= &
\mathbf{a}_2
-
L_1(\mathbf{a}_1-\mathbf{a}_2)(L_1)^{-1}
=
\mathbf{a}_2
-
2\big[L_1\circ (\mathbf{a}_1-\mathbf{a}_2)\big](L_1)^{-1}
-(\mathbf{a}_1-\mathbf{a}_2)\\
= &
2\mathbf{a}_2
-
2\big[L_1\circ (\mathbf{a}_1-\mathbf{a}_2)\big](L_1)^{-1}
-\mathbf{a}_1\in span \{ \mathbf{a}_1,\mathbf{a}_2,\mathbf{a}_3,\mathbf{a}_4\}\ .
\end{align*}
Thus, 
\begin{align*}
\overline{\mathbf{a}_1}-\mathbf{a}_1
& =
2(\mathbf{a}_2-\mathbf{a}_1)
-
2\big[L_1\circ (\mathbf{a}_1-\mathbf{a}_2)\big](L_1)^{-1}\\
& =
2(\mathbf{a}_2-\mathbf{a}_1)
+
2\big[(\mathbf{a}_1-\mathbf{a}_2)\circ L_1\big](L_1)^{-1}\\
& =
2(\mathbf{a}_2-\mathbf{a}_1)L_1(L_1)^{-1}
-
2\big[(\mathbf{a}_2-\mathbf{a}_1)\circ L_1\big](L_1)^{-1}\\
& =
2\Big\{
(\mathbf{a}_2-\mathbf{a}_1)L_1
-
\big[(\mathbf{a}_2-\mathbf{a}_1)\circ L_1\big]
\Big\}(L_1)^{-1}\\
& =
2
\big[(\mathbf{a}_2-\mathbf{a}_1)\Wedge L_1\big](L_1)^{-1}\\
& =
2
\Big[
(\mathbf{a}_2-\mathbf{a}_1)\wedge (\mathbf{a}_3-\mathbf{a}_2)\wedge(\mathbf{a}_4-\mathbf{a}_2)
\Big]
\big[(\mathbf{a}_3-\mathbf{a}_2)\wedge(\mathbf{a}_4-\mathbf{a}_2)\big]^{-1}\\
& =
2
\Big[
(\mathbf{a}_2-\mathbf{a}_1)\wedge (\mathbf{a}_3-\mathbf{a}_2)\wedge(\mathbf{a}_4-\mathbf{a}_2)
\Big]
\big[(\mathbf{a}_3-\mathbf{a}_2)\wedge(\mathbf{a}_4-\mathbf{a}_3)\big]^{-1}\ ,
\end{align*}
and you can verify that $|\overline{\mathbf{a}_1}-\mathbf{a}_2|=|\mathbf{a}_1-\mathbf{a}_2|$.
\begin{Rem}\label{rem:ortho-diagonals 2}
By using Remark~\ref{rem:rejection is orthogonal 2}, you can verify that the vector $\overline{\mathbf{a}_1}-\mathbf{a}_1$ is orthogonal to both $\mathbf{a}_2-\mathbf{a}_3$ and $\mathbf{a}_3-\mathbf{a}_4$. 
\end{Rem}
Then, we mirror the point $\mathbf{a}_2$ through the line $\mathcal{L}_{(\mathbf{a}_3,\mathbf{a}_4)}$, obtaining $\overline{\mathbf{a}_2}$ such that 
\begin{align*}
\overline{\mathbf{a}_2}-\mathbf{a}_2
& =
2\big[(\mathbf{a}_3-\mathbf{a}_2)\wedge(\mathbf{a}_4-\mathbf{a}_3)\big](\mathbf{a}_4-\mathbf{a}_3)^{-1}
\in span\{\mathbf{a}_3-\mathbf{a}_2\ ,\ \mathbf{a}_4-\mathbf{a}_3 \}\ ,
\end{align*}
as shown in Section~\ref{subsubsec: point mirrored by a line}. 
Thus we have now that vectors
$\overline{\mathbf{a}_1}-\mathbf{a}_1$, $\overline{\mathbf{a}_2}-\mathbf{a}_2$, and $\mathbf{a}_3-\mathbf{a}_4$ are mutually orthogonal. 
Moreover
\begin{align*}
(\overline{\mathbf{a}_1}-\mathbf{a}_1)
\wedge
(\overline{\mathbf{a}_2}-\mathbf{a}_2)
\wedge
(\mathbf{a}_4-\mathbf{a}_3)
& =
(\overline{\mathbf{a}_1}-\mathbf{a}_1)
(\overline{\mathbf{a}_2}-\mathbf{a}_2)
(\mathbf{a}_4-\mathbf{a}_3)\\
& =
2
(\overline{\mathbf{a}_1}-\mathbf{a}_1)
\big[(\mathbf{a}_3-\mathbf{a}_2)\wedge(\mathbf{a}_4-\mathbf{a}_3)\big]
(\mathbf{a}_4-\mathbf{a}_3)^{-1}(\mathbf{a}_4-\mathbf{a}_3)\\
& =
2
(\overline{\mathbf{a}_1}-\mathbf{a}_1)
\big[(\mathbf{a}_3-\mathbf{a}_2)\wedge(\mathbf{a}_4-\mathbf{a}_3)\big]\\
& =
4
(\mathbf{a}_2-\mathbf{a}_1)
\wedge(\mathbf{a}_3-\mathbf{a}_2)
\wedge(\mathbf{a}_4-\mathbf{a}_2)\\
& =
4
(\mathbf{a}_2-\mathbf{a}_1)
\wedge
(\mathbf{a}_3-\mathbf{a}_1)
\wedge
(\mathbf{a}_4-\mathbf{a}_1)\\
& =
4
\Delta_{(\mathbf{a}_1,\mathbf{a}_2,\mathbf{a}_3,\mathbf{a}_4)}
\end{align*}

\subsection{Convergence of the mean secant hyperplane plane to the tangent hyperplane I}
\begin{Def}\label{def:multi diff pseudo-vector n=3}
Let us define the ``mean multi-difference $\mathbb{G}_{3\choose 2}$-pseudo-vector'',
\begin{align*}
& 
\overline{\Delta}f_{(\mathbf{a}_1,\mathbf{a}_2,\mathbf{a}_3,\mathbf{a}_4)}
=\\
= &
\frac{1}{4}
\Big\{
\big[
f(\overline{\mathbf{a}_1})
-
f(\mathbf{a}_1)
\big]
(\overline{\mathbf{a}_2}-\mathbf{a}_2)
\wedge
(\mathbf{a}_4-\mathbf{a}_3)
-
\big[
f(\overline{\mathbf{a}_2})
-
f(\mathbf{a}_2)
\big]
(\overline{\mathbf{a}_1}-\mathbf{a}_1)
\wedge
(\mathbf{a}_4-\mathbf{a}_3)+\\
& 
+\big[
f(\mathbf{a}_4)
-
f(\mathbf{a}_3)
\big]
(\overline{\mathbf{a}_1}-\mathbf{a}_1)
\wedge 
(\overline{\mathbf{a}_2}-\mathbf{a}_2)
\Big\}\in \mathbb{G}_{3\choose 2}\ .
\end{align*}
\end{Def}
\begin{Rem*}
The term ``mean'' in the foregoing definition is due to the fact that 
\[
\overline{\Delta}f_{(\mathbf{a}_1,\mathbf{a}_2,\mathbf{a}_3,\mathbf{a}_4)}
=
\frac{1}{4}
\Big\{
\Delta f_{(\mathbf{a}_1,\mathbf{a}_2,\mathbf{a}_3,\mathbf{a}_4)}
-
\Delta f_{(\overline{\mathbf{a}_1},\mathbf{a}_2,\mathbf{a}_3,\mathbf{a}_4)}
-
\Delta f_{(\mathbf{a}_1,\overline{\mathbf{a}_2},\mathbf{a}_3,\mathbf{a}_4)}
+
\Delta f_{(\overline{\mathbf{a}_1},\overline{\mathbf{a}_2},\mathbf{a}_3,\mathbf{a}_4)}
\Big\}\ ,
\]
as you can verify. 
\end{Rem*}
\begin{Rem}\label{rem:mean 3,2 pseudovector geom prod}
The mean multi-difference $\mathbb{G}_{3\choose 2}$-pseudo-vector can also be writter as
\begin{align*}
& 
\overline{\Delta}f_{(\mathbf{a}_1,\mathbf{a}_2,\mathbf{a}_3,\mathbf{a}_4)}
=\\
= &
\frac{1}{4}
\Big\{
\big[
f(\overline{\mathbf{a}_1})
-
f(\mathbf{a}_1)
\big]
(\overline{\mathbf{a}_2}-\mathbf{a}_2)
(\mathbf{a}_4-\mathbf{a}_3)
-
\big[
f(\overline{\mathbf{a}_2})
-
f(\mathbf{a}_2)
\big]
(\overline{\mathbf{a}_1}-\mathbf{a}_1)
(\mathbf{a}_4-\mathbf{a}_3)+\\
& 
+\big[
f(\mathbf{a}_4)
-
f(\mathbf{a}_3)
\big]
(\overline{\mathbf{a}_1}-\mathbf{a}_1)
(\overline{\mathbf{a}_2}-\mathbf{a}_2)
\Big\}\ ,
\end{align*}
because vectors $\overline{\mathbf{a}_1}-\mathbf{a}_1$, $\overline{\mathbf{a}_2}-\mathbf{a}_2$ and $\mathbf{a}_4-\mathbf{a}_3$ are mutually orthogonal.
\end{Rem}
\begin{Prop}\label{lem:mean multi-difference quotient n=3}
\[
\overline{\mathbf{r}}_{f_{(\mathbf{a}_1,\mathbf{a}_2,\mathbf{a}_3,\mathbf{a}_4)}}
=
\frac{f(\overline{\mathbf{a}_1})-f(\mathbf{a}_1)}{\overline{\mathbf{a}_1}-\mathbf{a}_1}
+
\frac{f(\overline{\mathbf{a}_2})-f(\mathbf{a}_2)}{\overline{\mathbf{a}_2}-\mathbf{a}_2}
+
\frac{f(\mathbf{a}_4)-f(\mathbf{a}_3)}{\mathbf{a}_4-\mathbf{a}_3}\ .
\]
\end{Prop}
{\bf Proof of Proposition~\ref{lem:mean multi-difference quotient n=3}.}
From Remark~\ref{rem:ortho-diagonals 2} we have that $\overline{\mathbf{a}_1}-\mathbf{a}_1$, $\overline{\mathbf{a}_2}-\mathbf{a}_2$, $\mathbf{a}_4-\mathbf{a}_3$ are mutually orthogonal, and 
$
(\overline{\mathbf{a}_1}-\mathbf{a}_1)
(\overline{\mathbf{a}_2}-\mathbf{a}_2)
(\mathbf{a}_4-\mathbf{a}_3)
=  
4 \Delta_{(\mathbf{a}_1,\mathbf{a}_2,\mathbf{a}_3,\mathbf{a}_4)}
$. So, we can write
\[
\Delta_{(\mathbf{a}_1,\mathbf{a}_2,\mathbf{a}_3,\mathbf{a}_4)}
=
\frac{1}{4}
(\overline{\mathbf{a}_1}-\mathbf{a}_1)
\wedge
(\overline{\mathbf{a}_2}-\mathbf{a}_2)
\wedge
(\mathbf{a}_4-\mathbf{a}_3)
=
\frac{1}{4}
(\overline{\mathbf{a}_1}-\mathbf{a}_1)
(\overline{\mathbf{a}_2}-\mathbf{a}_2)
(\mathbf{a}_4-\mathbf{a}_3)
\]
\begin{align}
\big(\Delta_{(\mathbf{a}_1,\mathbf{a}_2,\mathbf{a}_3,\mathbf{a}_4)}\big)^{-1}
& =
4
(\mathbf{a}_4-\mathbf{a}_3)^{-1}
(\overline{\mathbf{a}_2}-\mathbf{a}_2)^{-1}
(\overline{\mathbf{a}_1}-\mathbf{a}_1)^{-1}
\nonumber \\
&=
-4
(\mathbf{a}_4-\mathbf{a}_3)^{-1}
(\overline{\mathbf{a}_1}-\mathbf{a}_1)^{-1}
(\overline{\mathbf{a}_2}-\mathbf{a}_2)^{-1}
\label{eq:volume as geom prod}\\
&=
4
(\overline{\mathbf{a}_2}-\mathbf{a}_2)^{-1}
(\mathbf{a}_4-\mathbf{a}_3)^{-1} \nonumber
\end{align}
Then, by Remark~\ref{rem:mean 3,2 pseudovector geom prod}, we can write
\begin{align*}
\overline{\mathbf{r}}_{f_{(\mathbf{a}_1,\mathbf{a}_2,\mathbf{a}_3,\mathbf{a}_4)}}
= &  
\Big(\overline{\mathbf{\Delta}}
f_{(\mathbf{a}_1,\mathbf{a}_2,\mathbf{a}_3,\mathbf{a}_4)}\Big)
\big(\Delta_{(\mathbf{a}_1,\mathbf{a}_2,\mathbf{a}_3,\mathbf{a}_4)}\big)^{-1}=\\
= & 
\big[
f(\overline{\mathbf{a}_1})
-
f(\mathbf{a}_1)
\big]
(\overline{\mathbf{a}_1}-\mathbf{a}_1)^{-1}
+
\big[
f(\overline{\mathbf{a}_2})
-
f(\mathbf{a}_2)
\big]
(\overline{\mathbf{a}_2}-\mathbf{a}_2)^{-1}+\\
& 
+
\big[
f(\mathbf{a}_4)
-
f(\mathbf{a}_3)
\big]
(\mathbf{a}_4-\mathbf{a}_3)^{-1}\ \square
\end{align*}
\begin{Thm}\label{thm:n=3}
If the function $f:\Omega\subseteq\mathbb{E}_3\to\mathbb{R}$ is strongly differentiable at $\mathbf{x}_0$ (a point internal to $\Omega$), then
\arraycolsep=1.pt\def\arraystretch{0.5}
\[
\lim_{\begin{array}{c}
	\scriptstyle (\mathbf{a}_1,\mathbf{a}_2,\mathbf{a}_3,\mathbf{a}_4)\to (\mathbf{x_0},\mathbf{x_0},\mathbf{x_0},\mathbf{x_0})\\
	\scriptstyle \mathbf{a}_1,\mathbf{a}_2,\mathbf{a}_3,\mathbf{a}_4 \textrm{ not coplanar}
\end{array}
}
\Big(\overline{\Delta}f_{(\mathbf{a}_1,\mathbf{a}_2,\mathbf{a}_3,\mathbf{a}_4)}\Big)
\big(\Delta_{(\mathbf{a}_1,\mathbf{a}_2,\mathbf{a}_3,\mathbf{a}_4)}\big)^{-1}
=\nabla f(\mathbf{x}_0)\ .
\]
\end{Thm}
\begin{Rem*}
The foregoing result state that the ``mean secant hyperplane''
\[
z= f(\mathbf{a}_1)\ +\ \overline{\mathbf{r}}_{f_{(\mathbf{a}_1,\mathbf{a}_2,\mathbf{a}_3,\mathbf{a}_4)}} \cdot (\mathbf{x}-\mathbf{a}_1)\ ,
\]
where 
$\overline{\mathbf{r}}_{f_{(\mathbf{a}_1,\mathbf{a}_2,\mathbf{a}_3,\mathbf{a}_4)}}
=
\Big(\overline{\Delta}f_{(\mathbf{a}_1,\mathbf{a}_2,\mathbf{a}_3,\mathbf{a}_4)}\Big)
\big(\Delta_{(\mathbf{a}_1,\mathbf{a}_2,\mathbf{a}_3,\mathbf{a}_4)}\big)^{-1}
$, always converges to the tangent hyperplane  
\[
z= f(\mathbf{x}_0)\ +\ \nabla f(\mathbf{x}_0) \cdot (\mathbf{x}-\mathbf{x}_0)\ ,
\]
as the non degenerate tetrahedron having vertices $\mathbf{a}_1$, $\mathbf{a}_2$, $\mathbf{a}_3$, and $\mathbf{a}_4$ contracts to the point $\mathbf{x}_0$.
\end{Rem*}
{\bf Proof of Theorem~\ref{thm:n=3}.}\\
Let us observe that 
\begin{align*}
\nabla f(\mathbf{x}_0)
\big(\Delta_{(\mathbf{a}_1,\mathbf{a}_2,\mathbf{a}_3,\mathbf{a}_4)}\big)
= &
\nabla f(\mathbf{x}_0)\circ
\big(\Delta_{(\mathbf{a}_1,\mathbf{a}_2,\mathbf{a}_3,\mathbf{a}_4)}\big)
\\
= &
\frac{1}{4}
\big[\nabla f(\mathbf{x}_0)\cdot(\overline{\mathbf{a}_1}-\mathbf{a}_1)\big]
(\overline{\mathbf{a}_2}-\mathbf{a}_2)
\wedge
(\mathbf{a}_4-\mathbf{a}_3)+\\
& 
-\frac{1}{4}
\big[\nabla f(\mathbf{x}_0)\cdot(\overline{\mathbf{a}_2}-\mathbf{a}_2)\big]
(\overline{\mathbf{a}_1}-\mathbf{a}_1)
\wedge
(\mathbf{a}_4-\mathbf{a}_3)+\\
& 
+ \frac{1}{4}
\big[\nabla f(\mathbf{x}_0)\cdot(\mathbf{a}_4-\mathbf{a}_3)\big]
(\overline{\mathbf{a}_1}-\mathbf{a}_1)
\wedge
(\overline{\mathbf{a}_2}-\mathbf{a}_2)\\
= &
\frac{1}{4}
\big[\nabla f(\mathbf{x}_0)\cdot(\overline{\mathbf{a}_1}-\mathbf{a}_1)\big]
(\overline{\mathbf{a}_2}-\mathbf{a}_2)
(\mathbf{a}_4-\mathbf{a}_3)+\\
& 
-\frac{1}{4}
\big[\nabla f(\mathbf{x}_0)\cdot(\overline{\mathbf{a}_2}-\mathbf{a}_2)\big]
(\overline{\mathbf{a}_1}-\mathbf{a}_1)
(\mathbf{a}_4-\mathbf{a}_3)+\\
& 
+ \frac{1}{4}
\big[\nabla f(\mathbf{x}_0)\cdot(\mathbf{a}_4-\mathbf{a}_3)\big]
(\overline{\mathbf{a}_1}-\mathbf{a}_1)
(\overline{\mathbf{a}_2}-\mathbf{a}_2)
\end{align*}
By relations~(\ref{eq:volume as geom prod}), we can write 
\begin{align*}
\nabla f(\mathbf{x}_0)
= &
\nabla f(\mathbf{x}_0)
\big(\Delta_{(\mathbf{a},\mathbf{b},\mathbf{c})}\big)
\big(\Delta_{(\mathbf{a},\mathbf{b},\mathbf{c})}\big)^{-1}
\\
= & 
\big[\nabla f(\mathbf{x}_0)\cdot(\overline{\mathbf{a}_1}-\mathbf{a}_1)\big]
(\overline{\mathbf{a}_1}-\mathbf{a}_1)^{-1}+\\
& 
+
\big[\nabla f(\mathbf{x}_0)\cdot(\overline{\mathbf{a}_2}-\mathbf{a}_2)\big]
(\overline{\mathbf{a}_2}-\mathbf{a}_2)^{-1}+\\
& 
+
\big[\nabla f(\mathbf{x}_0)\cdot(\mathbf{a}_4-\mathbf{a}_3)\big]
(\mathbf{a}_4-\mathbf{a}_3)^{-1}\\
= & 
\frac{\nabla f(\mathbf{x}_0)\cdot(\overline{\mathbf{a}_1}-\mathbf{a}_1)
}{\overline{\mathbf{a}_1}-\mathbf{a}_1}
+
\frac{\nabla f(\mathbf{x}_0)\cdot(\overline{\mathbf{a}_2}-\mathbf{a}_2)
}{\overline{\mathbf{a}_2}-\mathbf{a}_2}
+ 
\frac{\nabla f(\mathbf{x}_0)\cdot(\mathbf{a}_4-\mathbf{a}_3)}{\mathbf{a}_4-\mathbf{a}_3}\ .
\end{align*}
So, by Proposition~\ref{lem:mean multi-difference quotient n=3}, we can write
\begin{align*}
& 
\overline{\mathbf{r}}_{f_{(\mathbf{a}_1,\mathbf{a}_2,\mathbf{a}_3,\mathbf{a}_4)}}
-
\nabla f(\mathbf{x}_0)
=
\frac{f(\overline{\mathbf{a}_1})
-
f(\mathbf{a}_1)
-
\big[
\nabla f(\mathbf{x}_0)\cdot (\overline{\mathbf{a}_1}-\mathbf{a}_1)
\big]
}{\overline{\mathbf{a}_1}-\mathbf{a}_1}
+\\
& 
+
\frac{
f(\overline{\mathbf{a}_2})
-
f(\mathbf{a}_2)
-
\big[
\nabla f(\mathbf{x}_0)\cdot (\overline{\mathbf{a}_2}-\mathbf{a}_2)
\big]
}{\overline{\mathbf{a}_2}-\mathbf{a}_2}
+
\frac{f(\mathbf{a}_4)
-
f(\mathbf{a}_3)
-
\big[
\nabla f(\mathbf{x}_0)\cdot (\mathbf{a}_4-\mathbf{a}_3)
\big]
}{\mathbf{a}_4-\mathbf{a}_3}
\end{align*} 
Thus
\begin{align*}
& 
\big|
\overline{\mathbf{r}}_{f_{(\mathbf{a}_1,\mathbf{a}_2,\mathbf{a}_3,\mathbf{a}_4)}}
-
\nabla f(\mathbf{x}_0)
\big|
\le 
\frac{\big|f(\overline{\mathbf{a}_1})-f(\mathbf{a}_1)-\nabla f(\mathbf{x}_0)\cdot(\overline{\mathbf{a}_1}-\mathbf{a}_1)\big|}{|\overline{\mathbf{a}_1}-\mathbf{a}_1|}+\\
& +
\frac{\big|f(\overline{\mathbf{a}_2})-f(\mathbf{a}_2)-\nabla f(\mathbf{x}_0)\cdot(\overline{\mathbf{a}_2}-\mathbf{a}_2)\big|}{|\overline{\mathbf{a}_2}-\mathbf{a}_2|}
+
\frac{\big|f(\mathbf{a}_4)-f(\mathbf{a}_3)-\nabla f(\mathbf{x}_0)\cdot(\mathbf{a}_4-\mathbf{a}_3)\big|}{|\mathbf{a}_4-\mathbf{a}_3|}
\end{align*}
As $f$ is strongly differentiable at $\mathbf{x}_0$, we know that, given $\epsilon>0$ there exists $\delta_\epsilon>0$ such that if  $|\mathbf{u}-\mathbf{x}_0|<\delta_\epsilon$ and $|\mathbf{v}-\mathbf{x}_0|<\delta_\epsilon$, then
\begin{center}
$
|f(\mathbf{u})-f(\mathbf{v})-\nabla f(\mathbf{x}_0)\cdot(\mathbf{u}-\mathbf{v})|<\epsilon |\mathbf{u}-\mathbf{v}|\ .
$
\end{center}
So, if we choose $\mathbf{a}_1$, $\mathbf{a}_2$, $\mathbf{a}_3$, and $\mathbf{a}_4$ non collinear, and such that
\begin{center}
$|\mathbf{a}_i-\mathbf{x}_0|< \frac{1}{3}\left(\delta_{\frac{\epsilon}{3}}\right)$, for $i=1,2,3,4$,
\end{center}
then 
\[
\big|
\overline{\mathbf{r}}_{f_{(\mathbf{a}_1,\mathbf{a}_2,\mathbf{a}_3,\mathbf{a}_4)}}
-
\nabla f(\mathbf{x}_0)
\big|<\epsilon\ ,
\]
as, for $i=1,2,3$, we have
\begin{itemize}
	\item $|\mathbf{a}_i-\mathbf{a}_{i+1}|\le |\mathbf{a}_i-\mathbf{x}_0|+ |\mathbf{a}_{i+1}-\mathbf{x}_0|< \frac{2}{3}\left(\delta_{\frac{\epsilon}{3}}\right)$,
	\item $|\overline{\mathbf{a}_i}-\mathbf{x}_0|\le|\overline{\mathbf{a}_i}-\mathbf{a}_{i+1}|+|\mathbf{a}_{i+1}-\mathbf{x}_0|=|\mathbf{a}_i-\mathbf{a}_{i+1}|+|\mathbf{a}_{i+1}-\mathbf{x}_0| < \delta_{\frac{\epsilon}{3}}$.~$\square$
\end{itemize}
\begin{Rem*}
The key property that provide the convergence of the mean secant hyperplane is the following chain of identities 
\begin{center}
$
(\overline{\mathbf{a}_1}-\mathbf{a}_1)
(\overline{\mathbf{a}_2}-\mathbf{a}_2)
(\mathbf{a}_4-\mathbf{a}_3)
=
(\overline{\mathbf{a}_1}-\mathbf{a}_1)
\wedge
(\overline{\mathbf{a}_2}-\mathbf{a}_2)
\wedge
(\mathbf{a}_4-\mathbf{a}_3)
=
4
\Delta_{(\mathbf{a}_1,\mathbf{a}_2,\mathbf{a}_3,\mathbf{a}_4)}
$,
\end{center}
implied by the orthogonality properties of the constructed mirrored points based on the four non coplanar points $\mathbf{a}_1$, $\mathbf{a}_2$, $\mathbf{a}_3$, and $\mathbf{a}_4$.
\end{Rem*}
\section{The case of a multi-variable function}

\subsection{The $k\times k$ determinant as a Clifford quotient and as a scalar product}
The determinant of a $k\times k$ real matrix
\[
\left( 
\begin{array}{ccc}
	\mu_{1,1} & \cdots  & \mu_{1,k} \\
	\vdots & \ddots & \vdots \\
	\mu_{k,1} & \cdots & \mu_{k,k} 
\end{array}
\right)
\]
can be written as a coordinate-free Clifford quotient in $\mathbb{G}_n$ (with 
$n\ge k$). Let us fix any ordered set $\mathbf{e}_1,\dots ,\mathbf{e}_k$ of 
mutually orthonormal vectors in $\mathbb{E}_n$, and let us consider 
\[
\mathbf{u}_i=\sum_{j=1}^k\mu_{i,j}\ \mathbf{e}_j\ ,
\]
with $i=1,\dots,k$, then you can verify that
\[
(\mathbf{u}_1\wedge\cdots \wedge\mathbf{u}_k)(I_k)^{-1}
 =
\det
\left( 
\begin{array}{ccc}
	\mu_{1,1} & \cdots  & \mu_{1,k} \\
	\vdots & \ddots & \vdots \\
	\mu_{k,1} & \cdots & \mu_{k,k} 
\end{array}
\right)
\]
because 
\[
\mathbf{u}_{\sigma_1}\wedge\cdots \wedge\mathbf{u}_{\sigma_k}
=
\epsilon_\sigma\
\mathbf{u}_1\wedge\cdots \wedge\mathbf{u}_k\ ,
\]
for each permutation $\sigma\in\mathcal{S}_k$ of the set $\{1,\dots ,k\}$, 
having parity $\epsilon_\sigma\in\{-1,1\}$; where
\begin{align*}
I_k 
& = 
\mathbf{e}_1\cdots \mathbf{e}_k
=
\mathbf{e}_1\wedge \cdots \wedge \mathbf{e}_k\ , \ \textrm{ so that }\  (I_k)^
{-1}
=
\mathbf{e}_k\cdots \mathbf{e}_1=(-1)^{\frac{k(k-1)}{2}}I_k\ .
\end{align*}
\begin{Rem*}
As $I_2$ and $I_3$ before, also $I_k$ does not depend on the particular 
orthonormal basis of $span \{\mathbf{e}_1,\dots,\mathbf{e}_k\}=\mathbb{E}_k\subseteq \mathbb{E}_n$ to define it, but only on the 
orientation of that basis. More precisely, if $\{\mathbf{g}_1,\dots , \mathbf{g}_k\}$ is any other orthonormal basis of $\mathbb{E}_k$, then 
$
\mathbf{g}_1 \wedge\cdots \wedge \mathbf{g}_k
=
\mathbf{g}_1\cdots \mathbf{g}_k
$ is equal to 
$I_k$ or $-I_k$. That is why $I_k$ is called an ``orientation'' of~$\mathbb{E}_k$.
\end{Rem*}
Thus, a $k\times k$ determinant can be considered as the Clifford  ratio 
between the two ``$k$-blades'' 
$
\mathbf{u}_1\wedge\cdots \wedge\mathbf{u}_k$ and $I_k$ (a ``$k$-blade'' being the geometric product of $k$ non zero 
and mutually orthogonal vectors). Those elements can also be called 
``$\mathbb{G}_k$-pseudo-scalars'', as $\mathbb{G}_k$ is generated by~$\mathbb{E}_k
=span 
\{\mathbf{e}_1,\dots ,\mathbf{e}_k\}$, and can be interpreted as 
oriented hyper-volumes in $\mathbb{E}_k$.\\
In general, if $\{\mathbf{e}_1,\dots ,\mathbf{e}_n\}$  is an orthonormal basis for $\mathbb{E}_n$ and 
\begin{center}
$
\displaystyle 
V
=
\sum_{1\le i_1<\cdots < i_k\le n}
\nu_{i_1,\dots ,i_k}\mathbf{e}_{i_1}\cdots \mathbf{e}_{i_k}
\in 
\mathbb{G}_{n \choose k}$,
\end{center}
(such elements in $\mathbb{G}_n$ are also called the ``$k$-vectors''), then 
\begin{align*}
V I_n 
= &
\sum_{1\le i_1<\cdots < i_k\le n}
\nu_{i_1,\dots ,i_k}\mathbf{e}_{i_1}\cdots \mathbf{e}_{i_k}
\mathbf{e_1}\cdots \mathbf{e_n}\\
= & 
(-1)^{n-1}
\sum_{1\le i_1<\cdots < i_k\le n}
\nu_{i_1,\dots ,i_k}
\mathbf{e}_{i_1}\cdots \mathbf{e}_{i_{k-1}}
\mathbf{e_1}\cdots \mathbf{e_n}\mathbf{e}_{i_{k}}\\
= & 
(-1)^{(n-1)k}
\sum_{1\le i_1<\cdots < i_k\le n}
\nu_{i_1,\dots ,i_k}
\mathbf{e_1}\cdots \mathbf{e_n}
\mathbf{e}_{i_1}\cdots \mathbf{e}_{i_{k}} \in \mathbb{G}_{n\choose n-k}\\
= &
(-1)^{(n-1)k}
I_n V\ ,
\end{align*} 
that is why, when $k=n-1$, the elements in $\mathbb{G}_{n\choose n-1}$ are also called 
``$\mathbb{G}_n$-pseudo-vectors'': geometric multiplication by $I_n$ establishes a duality between elements of~$\mathbb{G}_{n \choose n-1}$ and vectors in $\mathbb{E}_n=\mathbb{G}_{n \choose 1}$. So, we have that $V I_n= (-1)^{(n-1)^2}I_n V = (-1)^{n-1}I_n V$, for each $\mathbb{G}_n$-pseudo-vector $V\in\mathbb{G}_{n\choose n-1}$.\\
The foregoing properties allows us to 
write a $k\times k$ determinant also as a scalar product
\begin{align*}
&
\det
\left( 
\begin{array}{ccc}
	\mu_{1,1} & \cdots  & \mu_{1,k} \\
	\vdots & \ddots & \vdots \\
	\mu_{k,1} & \cdots & \mu_{k,k} 
\end{array}
\right)
= 
(\mathbf{u}_1\wedge\cdots \wedge\mathbf{u}_k)(I_k)^{-1}
= 
(-1)^{\frac{k(k-1)}{2}}
(\mathbf{u}_1\wedge\cdots \wedge\mathbf{u}_k)I_k=\\
= &
(-1)^{\frac{k(k-1)}{2}}
\big[(\mathbf{u}_1\wedge\cdots\wedge\mathbf{u}_{k-1})\Wedge\mathbf{u}_k\big]I_k\\
= &
\frac{(-1)^{\frac{k(k-1)}{2}}}{2}
\Big[
(\mathbf{u}_1\wedge\cdots\wedge\mathbf{u}_{k-1})\mathbf{u}_k
+
(-1)^{k-1}
\mathbf{u}_k
(\mathbf{u}_1\wedge\cdots\wedge\mathbf{u}_{k-1})
\Big]
I_k\\
= &
\frac{(-1)^{\frac{k(k-1)}{2}}}{2}
\Big[
(\mathbf{u}_1\wedge\cdots\wedge\mathbf{u}_{k-1})\mathbf{u}_k I_k
+
(-1)^{k-1}
\mathbf{u}_k
(\mathbf{u}_1\wedge\cdots\wedge\mathbf{u}_{k-1})I_k
\Big]\\
= &
(-1)^{k-1}\frac{(-1)^{\frac{k(k-1)}{2}}}{2}
\Big[
(\mathbf{u}_1\wedge\cdots\wedge\mathbf{u}_{k-1})I_k\mathbf{u}_k
+
\mathbf{u}_k
(\mathbf{u}_1\wedge\cdots\wedge\mathbf{u}_{k-1})I_k
\Big]\\
= &
(-1)^{k-1}
\big[
(\mathbf{u}_1\wedge\cdots\wedge\mathbf{u}_{k-1})(I_k)^{-1}
\big]\cdot \mathbf{u}_k\ .
\end{align*}
We can also verify that the vector 
$
(\mathbf{u}_1\wedge\cdots\wedge\mathbf{u}_{k-1})I_k
$
is a vector orthogonal to $span\{\mathbf{u}_1,\dots ,\mathbf{u}_{k-1}\}$, because it is orthogonal to each $\mathbf{u}_i$ (when $i=1,\dots, k-1$), as you can verify that 
\begin{center}
$
\big[(\mathbf{u}_1\wedge \cdots \wedge\mathbf{u}_{k-1})(I_k)^{-1}\big]\cdot\mathbf{u}_i
=
(\mathbf{u}_1\wedge \cdots \wedge\mathbf{u}_{k-1}\wedge\mathbf{u}_i)(I_k)^{-1}
=
0
$.
\end{center}

\subsection{Coordinate-free expression of a hyperplane secant the graph of a multi-variable function}
Let us write the equation of a hyperplane secant the graph of a multi-variable function $f:\Omega\subseteq\mathbb{E}_n\to\mathbb{R}$ at $n+1$  points of that graph non being on a same $n$-dimensional hyperplane. If ${\mathbf{e}_1,\dots,\mathbf{e}_n}$ is an orthonormal basis for $\mathbb{E}_n$, a hyperplane passing through such $n+1$ points 
\[
\big(\mathbf{a}_i,f(\mathbf{a}_i)\big)
=
\left(
\sum_{j=1}^n
\alpha_{i,j}\mathbf{e}_j,f(\mathbf{a}_i)
\right)
\in \mathbb{E}_n\oplus\mathbb{R}\ ,
\]
with $i=1,\dots,n+1 $, can be represented by the Cartesian relation
\begin{equation}
\label{eq:general secant hyperplane}
\det
\left(
\begin{array}{cccc}
	\chi_1-\alpha_{1,1} & \cdots & \chi_n-\alpha_{1,n} & z-f(\mathbf{a}_1)\\
	\alpha_{2,1}-\alpha_{1,1} & \cdots & \alpha_{2,n}-\alpha_{1,n} & f(\mathbf{a}_2)-f(\mathbf{a}_1)\\
	\vdots  & \ddots  & \vdots & \vdots\\
	\alpha_{n+1,1}-\alpha_{1,1} & \cdots & \alpha_{n+1,n}-\alpha_{1,n} & f(\mathbf{a}_{n+1})-f(\mathbf{a}_1)
\end{array}
\right)
=0
\end{equation}
between the $n+1$ real variables $\chi_1,\dots ,\chi_n,z\in\mathbb{R}$. That determinant can be rewritten by a  Laplace expansion as follows
\begin{align*}
&
\big[z-f(\mathbf{a}_1)\big]
\det
\left(
\begin{array}{ccc}
\alpha_{2,1}-\alpha_{1,1} & \cdots  & \alpha_{2,n}-\alpha_{1,n} \\
	\vdots  & \ddots  & \vdots  \\
	\alpha_{n+1,1}-\alpha_{1,1} & \cdots & \alpha_{n+1,n}-\alpha_{1,n} 
\end{array}
\right)
+\\
& -
\big[f(\mathbf{a}_2)-f(\mathbf{a}_1)\big]
\det
\left(
\begin{array}{ccc}
	\chi_1-\alpha_{1,1} & \cdots  & \chi_n-\alpha_{1,n}\\
	\alpha_{3,1}-\alpha_{1,1} & \cdots & \alpha_{3,n}-\alpha_{1,n} \\
	\vdots  & \ddots & \vdots  \\
	\alpha_{n+1,1}-\alpha_{1,1} & \cdots & \alpha_{n+1,n}-\alpha_{1,n} 
\end{array}
\right)
+\\
& +
\big[f(\mathbf{a}_3)-f(\mathbf{a}_1)\big]
\det
\left(
\begin{array}{ccc}
	\chi_1-\alpha_{1,1} & \cdots  & \chi_n-\alpha_{1,n} \\
	\alpha_{2,1}-\alpha_{1,1} & \cdots & \alpha_{2,n}-\alpha_{1,n} \\
	\alpha_{4,1}-\alpha_{1,1} & \cdots & \alpha_{4,n}-\alpha_{1,n} \\
	\vdots  & \ddots & \vdots  \\
\alpha_{n+1,1}-\alpha_{1,1} & \cdots & \alpha_{n+1,n}-\alpha_{1,n} 
\end{array}
\right)
+\\
& +\cdots +\\
& +(-1)^{n}
\big[f(\mathbf{a}_{n+1})-f(\mathbf{a}_1)\big]
\det
\left(
\begin{array}{ccc}
	\chi_1-\alpha_{1,1} & \cdots  & \chi_n-\alpha_{1,n} \\
	\vdots  & \ddots & \vdots \\
	\alpha_{n,1}-\alpha_{1,1} & \cdots & \alpha_{n,n}-\alpha_{1,n} 
\end{array}
\right)
\end{align*}
Then, in $\mathbb{G}_3$ the equation~(\ref{eq:general secant hyperplane}) becomes 
\begin{align*}
&
\big[z-f(\mathbf{a}_1)\big]
\big[(\mathbf{a}_2-\mathbf{a}_1)\wedge 
\cdots 
\wedge(\mathbf{a}_{n+1}-\mathbf{a}_1)\big]
(I_n)^{-1}+\\
- &
\big[f(\mathbf{a}_2)-f(\mathbf{a}_1)\big]
\big[
(\mathbf{x}-\mathbf{a}_1)
\wedge
(\mathbf{a}_3-\mathbf{a}_1)
\wedge \cdots \wedge
(\mathbf{a}_{n+1}-\mathbf{a}_1)
\big]
(I_n)^{-1}+\\
+ &
\big[f(\mathbf{a}_3)-f(\mathbf{a}_1)\big]
\big[
(\mathbf{x}-\mathbf{a}_1)
\wedge
(\mathbf{a}_2-\mathbf{a}_1)
\wedge
(\mathbf{a}_4-\mathbf{a}_1)
\wedge \cdots\wedge
(\mathbf{a}_{n+1}-\mathbf{a}_1)
\big]
(I_n)^{-1}+\\
+ &\cdots + \\
+ & (-1)^{n} 
\big[f(\mathbf{a}_{n+1})-f(\mathbf{a}_1)\big]
\big[
(\mathbf{x}-\mathbf{a}_1)
\wedge
(\mathbf{a}_2-\mathbf{a}_1)
\wedge\cdots \wedge
(\mathbf{a}_n-\mathbf{a}_1)\big]
(I_3)^{-1}
=
0\ ,
\end{align*}
being $\displaystyle \mathbf{x}=\sum_{i=1}^n\chi_i\mathbf{e}_i\in\mathbb{E}_n$, $(\mathbf{x},z)\in\mathbb{E}_n\oplus\mathbb{R}$. The foregoing relation is equivalent, in~$\mathbb{G}_3$, to
\begin{align*}
&
\big[z-f(\mathbf{a}_1)\big]
\big[(\mathbf{a}_2-\mathbf{a}_1)\wedge 
\cdots 
\wedge(\mathbf{a}_{n+1}-\mathbf{a}_1)\big]
=\\
= &
\big[f(\mathbf{a}_2)-f(\mathbf{a}_1)\big]
\big[
(\mathbf{x}-\mathbf{a}_1)
\wedge
(\mathbf{a}_3-\mathbf{a}_1)
\wedge \cdots \wedge
(\mathbf{a}_{n+1}-\mathbf{a}_1)
\big]+\\
& -
\big[f(\mathbf{a}_3)-f(\mathbf{a}_1)\big]
\big[
(\mathbf{x}-\mathbf{a}_1)
\wedge
(\mathbf{a}_2-\mathbf{a}_1)
\wedge
(\mathbf{a}_4-\mathbf{a}_1)
\wedge \cdots\wedge
(\mathbf{a}_{n+1}-\mathbf{a}_1)
\big]+\\
& + \cdots + \\
& + (-1)^{n+1} 
\big[f(\mathbf{a}_{n+1})-f(\mathbf{a}_1)\big]
\big[
(\mathbf{x}-\mathbf{a}_1)
\wedge
(\mathbf{a}_2-\mathbf{a}_1)
\wedge\cdots \wedge
(\mathbf{a}_n-\mathbf{a}_1)\big]\ ,
\end{align*}
Let us define 
\begin{align*}
\Delta_{(\mathbf{a}_1,\dots ,\mathbf{a}_{n+1})}
= &
(\mathbf{a}_2-\mathbf{a}_1)\wedge(\mathbf{a}_3-\mathbf{a}_1)
\wedge\cdots\wedge
(\mathbf{a}_{n+1}-\mathbf{a}_1)
=
\bigwedge_{i=2}^{n+1}
(\mathbf{a}_i-\mathbf{a}_1)
\\
& \Big(
=
(-1)^{n}
(\mathbf{a}_1-\mathbf{a}_2)
\wedge
(\mathbf{a}_2-\mathbf{a}_3)
\wedge\cdots\wedge
(\mathbf{a}_n-\mathbf{a}_{n+1})
\Big)
\end{align*}
We observe that 
\begin{center}
$
\displaystyle \tau_n =\frac{1}{n!}
\det
\left(
\begin{array}{ccc}
\alpha_{2,1}-\alpha_{1,1} & \cdots  & \alpha_{2,n}-\alpha_{1,n} \\
	\vdots  & \ddots & \vdots  \\
	\alpha_{n+1,1}-\alpha_{1,1} & \cdots & \alpha_{n+1,n}-\alpha_{1,n}
\end{array}
\right)
=
\frac{1}{n!}
\Delta_{(\mathbf{a}_1,\dots,\mathbf{a}_{n+1})}
(I_n)^{-1}
$
\end{center}
is the oriented hyper-volume of the simplex having vertices $\mathbf{a}_1$, \dots , $\mathbf{a}_n$, and $\mathbf{a}_{n+1}$. 
So, we can write 
$
\Delta_{(\mathbf{a}_1,\dots ,\mathbf{a}_{n+1})}
=
n!\tau_n I_n
$, and $\big(\Delta_{(\mathbf{a}_1,\dots,\mathbf{a}_{n+1})}\big)^{-1}= \frac{1}{n!\tau_n}(I_n)^{-1}$. By denoting $\Delta f_{(\mathbf{a}_i)} = f(\mathbf{a}_i)-f(\mathbf{a}_1)$ and $\mathbf{\Delta}\mathbf{a}_i=\mathbf{a}_i-\mathbf{a}_1$, when $i=2,\dots ,n+1$, the equation of the secant hyperplane~(\ref{eq:general secant hyperplane}) becomes
\begin{align*}
z =
 & 
f(\mathbf{a}_1)
+ \frac{1}{n!\tau_n}
\Big\{
\Delta f_{(\mathbf{a}_2)}
(\mathbf{x}-\mathbf{a}_1)
\wedge 
\mathbf{\Delta}\mathbf{a}_3
\wedge \cdots \wedge 
\mathbf{\Delta}\mathbf{a}_{n+1}+\\
& -
\Delta f_{(\mathbf{a}_3)}
(\mathbf{x}-\mathbf{a}_1)\wedge 
\mathbf{\Delta}\mathbf{a}_2\wedge 
\mathbf{\Delta}\mathbf{a}_4
\wedge\cdots\wedge
\mathbf{\Delta}\mathbf{a}_{n+1}
+\cdots +\\
& +(-1)^{n+1}
\Delta f_{(\mathbf{a}_{n+1})}
(\mathbf{x}-\mathbf{a}_1)
\wedge 
\mathbf{\Delta}\mathbf{a}_2
\wedge \cdots\wedge 
\mathbf{\Delta}\mathbf{a}_n\Big\}(I_n)^{-1}
\end{align*}
If we define 
\arraycolsep=1.pt\def\arraystretch{0.5}
\[
\Delta^i_{(\mathbf{a}_1,\dots,\mathbf{a}_{n+1})}
=
\bigwedge_{
\begin{array}{c}
	\scriptstyle j=2\\
	\scriptstyle j\ne i
\end{array}
}^{n+1}
\mathbf{\Delta}\mathbf{a}_j
=
\Delta^i
\in\mathbb{G}_{n\choose n-1}\ ,
\]
the equation of the secant hyperplane~(\ref{eq:general secant hyperplane}) becomes
\begin{align*}
& z = 
f(\mathbf{a}_1)+\\
+ &
\frac{1}{n!\tau_n}
\Big\{
\Delta f_{(\mathbf{a}_2)}
(\mathbf{x}-\mathbf{a}_1)
\Wedge
\Delta^2
-
\Delta f_{(\mathbf{a}_3)}
(\mathbf{x}-\mathbf{a}_1)
\Wedge
\Delta^3
+ \cdots +
\Delta f_{(\mathbf{a}_{n+1})}
(\mathbf{x}-\mathbf{a}_1)
\Wedge
\Delta^{n+1}
\Big\}(I_n)^{-1}\\
& = 
f(\mathbf{a}_1)
+ \frac{1}{n!\tau_n}
\Big\{
(\mathbf{x}-\mathbf{a}_1)
\Wedge
\big[
\Delta f_{(\mathbf{a}_2)}
\Delta^2
-
\Delta f_{(\mathbf{a}_3)}
\Delta^3
+ \cdots +
\Delta f_{(\mathbf{a}_{n+1})}
\Delta^{n+1}
\big]
\Big\}(I_n)^{-1}\\
& = 
f(\mathbf{a}_1)
+ \frac{(-1)^{n-1}}{n!\tau_n}
\Big\{
\big[
\Delta f_{(\mathbf{a}_2)}
\Delta^2
-
\Delta f_{(\mathbf{a}_3)}
\Delta^3
+ \cdots +
\Delta f_{(\mathbf{a}_{n+1})}
\Delta^{n+1}
\big]
\Wedge
(\mathbf{x}-\mathbf{a}_1)
\Big\}
(I_n)^{-1}\\
& = 
f(\mathbf{a}_1)
+ \frac{1}{n!\tau_n}
\Big\{
\big[
\Delta f_{(\mathbf{a}_2)}
\Delta^2
-
\Delta f_{(\mathbf{a}_3)}
\Delta^3
+ \cdots +
\Delta f_{(\mathbf{a}_{n+1})}
\Delta^{n+1}
\big]
(I_n)^{-1}
\Big\}
\cdot
(\mathbf{x}-\mathbf{a}_1)\\
& = 
f(\mathbf{a}_1)
+ 
\Big\{
\big[
\Delta f_{(\mathbf{a}_2)}
\Delta^2
-
\Delta f_{(\mathbf{a}_3)}
\Delta^3
+ \cdots +
\Delta f_{(\mathbf{a}_{n+1})}
\Delta^{n+1}
\big]
(\Delta_{(\mathbf{a}_1,\dots ,\mathbf{a}_{n+1})})^{-1}
\Big\}
\cdot
(\mathbf{x}-\mathbf{a}_1)
\end{align*} 
Thus, we have that the vector $\mathbf{r}_{f_{(\mathbf{a}_1,\dots,\mathbf{a}_{n+1})}}\in\mathbb{E}_n$, 
characterizing the equation of the hyperplane secant the graph of $f$ as 
$
z
=
f(\mathbf{a}_1)+\mathbf{r}_{f_{(\mathbf{a}_1,\dots,\mathbf{a}_{n+1})}}\cdot (\mathbf{x}-\mathbf{a}_1)
$,
is
\[
\mathbf{r}_{f_{(\mathbf{a}_1,\dots ,\mathbf{a}_{n+1})}}
=
\left[
\sum_{i=2}^{n+1}
(-1)^i
\Delta f_{(\mathbf{a}_i)}
\Delta^i_{(\mathbf{a}_1,\dots,\mathbf{a}_{n+1})}
\right]
\big(\Delta_{(\mathbf{a}_1,\dots,\mathbf{a}_{n+1})}\big)^{-1}\ ,
\]
which is, in fact, the Clifford quotient between the multi-difference $\mathbb{G}_n$-pseudo-vector 
\[
\Delta f_{(\mathbf{a}_1,\dots,\mathbf{a}_{n+1})}
=
\sum_{i=2}^{n+1}
(-1)^i
\Delta f_{(\mathbf{a}_i)}
\Delta^i_{(\mathbf{a}_1,\dots,\mathbf{a}_{n+1})}
\in\mathbb{G}_{n\choose n-1}\ ,
\]
and the $\mathbb{G}_n$-pseudo-scalar
$\Delta_{(\mathbf{a}_1,\dots,\mathbf{a}_{n+1})}\in \mathbb{G}_{n\choose n}\simeq \mathbb{R}I_n$.

\subsection{Mirroring vectors and points III}

Here, we simply iterate the process already followed in the two foregoing cases. More precisely, if $n>3$, we consider $n+1$ points $\mathbf{a}_1$, \dots, $\mathbf{a}_{n+1}$ in $\mathbb{E}_n$ as vertices of a non degenerate simplex (that is, $\Delta_{(\mathbf{a}_1,\dots,\mathbf{a}_{n+1})}\ne 0$), we construct $n-1$ mirrored points $\overline{\mathbf{a}_1}$, \dots, $\overline{\mathbf{a}_{n-1}}$ in $\mathbb{E}_n$ in the following way
\begin{itemize}
	\item $\overline{\mathbf{a}_1}$ is obtained by mirroring point $\mathbf{a}_1$ through the hyperplane $\mathcal{L}_{(\mathbf{a}_2,\dots,\mathbf{a}_{n+1})}$;
	\item $\overline{\mathbf{a}_2}$ is obtained by mirroring point $\mathbf{a}_2$ through the hyperplane $\mathcal{L}_{(\mathbf{a}_3,\dots,\mathbf{a}_{n+1})}$;
	\item $\vdots$
	\item $\overline{\mathbf{a}_{n-2}}$ is obtained by mirroring point $\mathbf{a}_{n-2}$ through the plane $\mathcal{L}_{(\mathbf{a}_{n-1},\mathbf{a}_{n},\mathbf{a}_{n+1})}$;
	\item $\overline{\mathbf{a}_{n-1}}$ is obtained by mirroring point $\mathbf{a}_{n-1}$ through the line $\mathcal{L}_{(\mathbf{a}_{n},\mathbf{a}_{n+1})}$,
\end{itemize}
where
$
\mathcal{L}_{(\mathbf{v}_1,\dots,\mathbf{v}_h)}
=
\{
\nu_1\mathbf{v}_1 +\cdots +\nu_h \mathbf{v}_h\ : \ \nu_1,\dots,\nu_h\in\mathbb{R},\ \nu_1+\cdots+\nu_h=1
\}
$, is the hyperplane passing through points 
$
\mathbf{v}_1,\dots,\mathbf{v}_h
\in\mathbb{E}_n
$.\\
All that points $\overline{\mathbf{a}_{i}}$ can be obtained by mirroring vectors $\mathbf{a}_{i}-\mathbf{a}_{i+1}$ through the $(n-i)$-dimensional subspace $span \{\mathbf{a}_{i+2}-\mathbf{a}_{i+1}, \dots, \mathbf{a}_{n+1}-\mathbf{a}_{i+1}\}$, respectively, by using the corresponding $(n-i)$-blade 
\begin{center}
$
L_i
=
(\mathbf{a}_{i+2}-\mathbf{a}_{i+1})
\wedge\cdots\wedge
(\mathbf{a}_{n+1}-\mathbf{a}_{i+1})
$.
\end{center}
as computed in the two following paragraphs.

\subsubsection{Vector mirrored by a $k$-dimensional linear subspace}
We recall that to $k$ linearly independent vectors $\mathbf{u}_1$,\dots $\mathbf{u}_k \in\mathbb{E}_n$ we can associate the $k$-dimensional linear subspace $span\{\mathbf{u}_1,\dots ,\mathbf{u}_k\}\subseteq \mathbb{E}_n$.
Moreover, $\mathbf{u}_1\wedge\cdots \wedge\mathbf{u}_k$ is always a $k$-blade. As a matter of fact, there always exists an orthogonal basis $\{\mathbf{g}_1,\dots,\mathbf{g}_k\}$ of  $span\{\mathbf{u}_1,\dots,\mathbf{u}_k\}$, and you can verify that $\mathbf{u}_1\wedge\cdots\wedge\mathbf{u}_k$ is a non zero multiple of the geometric product $\mathbf{g}_1\cdots\mathbf{g}_k$.
We recall that 
\begin{itemize}
	\item the square of every $k$-blade is a non zero scalar, 
	\item every $k$-blade~$B$ is invertible in $\mathbb{G}_n$, and $\displaystyle B^{-1}=\frac{1}{B^2}B$.
\end{itemize}
Let us recall that, if $\mathbf{v}\in\mathbb{E}_n$, and $B=\mathbf{u}_1\wedge\cdots\wedge\mathbf{u}_k$ is a $k$-blade, then
\begin{equation}
\label{eq:k-blade vector decomposition}
\mathbf{v}
=
\mathbf{v}B B^{-1}
=
(\mathbf{v}B)B^{-1}
=
(\mathbf{v}\circ B + \mathbf{v}\Wedge B)B^{-1}
=
(\mathbf{v}\circ B)B^{-1}+(\mathbf{v}\Wedge B)B^{-1}
\end{equation}
where
\begin{align*}
\mathbf{v} \circ B
& =
\frac{1}{2}\big(\mathbf{v}B-(-1)^k B\mathbf{v}\big)
=
(-1)^{k+1} B \circ\mathbf{v}
\in\mathbb{G}_{n\choose k-1}\ ,\\
\mathbf{v} \Wedge B
& =
\frac{1}{2}\big(\mathbf{v}B+(-1)^k B\mathbf{v}\big)
=
(-1)^{k} B \Wedge\mathbf{v}
\in\mathbb{G}_{n\choose k+1}\ .
\end{align*}
We also recall that
\arraycolsep=1.pt\def\arraystretch{0.5}
\begin{eqnarray}
\label{eq:circ}
\mathbf{v}\circ (\mathbf{u}_1\wedge \cdots\wedge\mathbf{u}_k)
& =
\displaystyle
\sum_{i=1}^k (-1)^{i+1}(\mathbf{v}\cdot \mathbf{u}_i) 
\bigwedge_{
\begin{array}{c}
\scriptstyle j=1 \\
\scriptstyle j\ne i
\end{array}}^k
\mathbf{u}_j\\
\mathbf{v}\Wedge (\mathbf{u}_1\wedge \cdots\wedge\mathbf{u}_k)
& = 
\mathbf{v}\wedge \mathbf{u}_1\wedge \cdots\wedge\mathbf{u}_k\ .
\end{eqnarray}
So, we have that
\[
B\mathbf{v}
=
2(B\circ\mathbf{v})+(-1)^k\mathbf{v}B
=
(B\circ \mathbf{v})
+
(B\Wedge \mathbf{v})
\ .
\]
\begin{Rem}\label{rem:rejection is orthogonal k}
If $B=\mathbf{u}_1\wedge\cdots\wedge\mathbf{u}_k$ is a $k$-blade, then $(\mathbf{v}\Wedge B)B^{-1}$ is a vector which is orthogonal to all $\mathbf{u}_1$,\dots, $\mathbf{u}_k$ that is, to all $k$-dimensional $span\{ \mathbf{u}_1,\dots, \mathbf{u}_k\}$. As a matter of fact, for $i=1,\dots,k$ we have that 
\[
2\big[(\mathbf{v}\Wedge B)B^{-1}\big]\cdot \mathbf{u}_i
=
(\mathbf{v}\Wedge B)B^{-1} \mathbf{u}_i
+
\mathbf{u}_i(\mathbf{v}\Wedge B)B^{-1}
\]
You can verify that $B\mathbf{u}_i = B\circ\mathbf{u}_i = (-1)^{k+1}\mathbf{u}_i \circ B= (-1)^{k+1}B \mathbf{u}_i$ (for $i=1,\dots,k$). So we have that 
\[
2\big[(\mathbf{v}\Wedge B)B^{-1}\big]\cdot \mathbf{u}_i
=
\big[
(-1)^{k+1}(\mathbf{v}\Wedge B) \mathbf{u}_i 
+
\mathbf{u}_i(\mathbf{v}\Wedge B)
\big]B^{-1}
\]
As $\mathbf{v}\Wedge B= \mathbf{v}\wedge \mathbf{u}_1 \wedge\cdots\wedge  \mathbf{u}_k$ is either zero or a $(k+1)$-blade, then we can write
\[
4\big[(\mathbf{v}\Wedge B)B^{-1}\big]\cdot \mathbf{u}_i
=
\big[\mathbf{u}_i\Wedge (\mathbf{v}\Wedge B)
\big]B^{-1}
=
\big(\mathbf{u}_i\wedge \mathbf{v}\wedge \mathbf{u}_1\wedge\cdots\wedge \mathbf{u}_k
\big)B^{-1}
=
0\ .
\]
\end{Rem}
In a similar way you can verify that for every vector $\mathbf{v}\in\mathbb{E}_n$ we have that 
\begin{center}
$\mathbf{v}_\parallel= (\mathbf{v}\circ B)B^{-1}\in span\{\mathbf{u}_1,\dots,\mathbf{u}_k\}$.
\end{center}
So, relation~(\ref{eq:k-blade vector decomposition}) corresponds to the orthogonal  decomposition of vector $\mathbf{v}$ with respect to the $k$-dimensional linear subspace $span\{\mathbf{u}_1,\dots,\mathbf{u}_k\}=\mathcal{L}_{(\mathbf{0},\mathbf{u}_1,\dots,\mathbf{u}_k)}$.\\
Then, the vector $\mathbf{\hat{v}}$, mirrored of vector $\mathbf{v}$ by the $k$-dimensional linear subspace $span\{\mathbf{u}_1,\dots,\mathbf{u}_k\}$, can be written as
\begin{align*}
\mathbf{\hat{v}}
& = 
\mathbf{v}_{\parallel}-\mathbf{v}_{\bot}
=
(\mathbf{v}\circ B)B^{-1}
-
(\mathbf{v}\Wedge B)B^{-1}
=
\big[(\mathbf{v}\circ B)-(\mathbf{v}\Wedge B)\big]B^{-1}\\
& =
\big[(-1)^{k+1}(B\circ\mathbf{v})-(-1)^{k}(B\Wedge \mathbf{v})\big]B^{-1}
=
(-1)^{k+1}\big[(B\circ\mathbf{v})+(B\Wedge \mathbf{v})\big]B^{-1}\\
& =
(-1)^{k+1}B\mathbf{v}B^{-1}
=
(-1)^{k+1}\Big[2(B\circ\mathbf{v})+(-1)^{k}\mathbf{v}B\big]B^{-1}\\
& = 
(-1)^{k+1}2(B\circ\mathbf{v})B^{-1}-\mathbf{v}\ .
\end{align*}
Moreover, $|\mathbf{\hat{v}}|= |\mathbf{v}|$; as a matter of fact,
\begin{center}
$
|\mathbf{\hat{v}}|^2
=
\mathbf{\hat{v}}\mathbf{\hat{v}}
=
B\mathbf{v}B^{-1}
B\mathbf{v}B^{-1}
=
B\mathbf{v}\mathbf{v}B^{-1}
=
|\mathbf{v}|^2BB^{-1}
=
|\mathbf{v}|^2
$.
\end{center}
\begin{Rem*}
Notice that, $\mathbf{\hat{v}}$ is in $span\{\mathbf{u}_1,\dots,\mathbf{u}_k,\mathbf{v}\}$.
\end{Rem*}

\subsubsection{Point mirrored by a multi-dimensional hyperplane in $\mathbb{E}_n$}
Let $n>3$, we consider $n+1$ points $\mathbf{a}_1$, \dots, $\mathbf{a}_{n+1}$ in $\mathbb{E}_n$ such that $\Delta_{(\mathbf{a}_1,\dots,\mathbf{a}_{n+1})}\ne 0$.\\
While $1\le i\le n-1$, we consider the a hyperplane~$\mathcal{L}_{(\mathbf{a}_{i+1},\dots, \mathbf{a}_{n+1})}$ passing through $n+1-i$ points~$\mathbf{a}_{i+1}$,\dots  $\mathbf{a}_{n+1}$. Then, we mirror the point~$\mathbf{a}_i$ through that hyperplane~$\mathcal{L}_{(\mathbf{a}_{i+1},\dots, \mathbf{a}_{n+1})}$. We denote $\overline{\mathbf{a}_i}$ that mirrored point.\\
We can express  $\overline{\mathbf{a}_i}$ by computing in $\mathbb{G}_n$ the vector $\mathbf{\hat{v}}$, obtained by mirroring vector $\mathbf{v}=\mathbf{a}_i-\mathbf{a}_{i+1}$ by the $(n-i)$-blade 
\[
L_i
=
L_{(\mathbf{a}_{i+1},\dots,\mathbf{a}_{n+1})}
=
(\mathbf{a}_{i+2}-\mathbf{a}_{i+1})\wedge\cdots\wedge(\mathbf{a}_{n+1}-\mathbf{a}_{i+1})
=
\bigwedge_{j=2}^{n+1-i}
(\mathbf{a}_{i+j}-\mathbf{a}_{i+1})\ ,
\]
for $i=1,\dots,n-1$.
More precisely, the  mirrored point $\overline{\mathbf{a}_1}$
can be computed by using the geometric Clifford product in $\mathbb{G}_n$ as follows
\begin{align*}
\overline{\mathbf{a}_i}
= &
\mathbf{a}_{i+1}
+
(-1)^{n-i+1}
L_i(\mathbf{a}_i-\mathbf{a}_{i+1})(L_i)^{-1}\\
= &
\mathbf{a}_{i+1}
+
(-1)^{n-i+1}
2\big[L_i\circ (\mathbf{a}_i-\mathbf{a}_{i+1})\big](L_i)^{-1}
-(\mathbf{a}_i-\mathbf{a}_{i+1})\\
= &
2\mathbf{a}_{i+1}
+
(-1)^{n-i+1}
2\big[L_i\circ (\mathbf{a}_i-\mathbf{a}_{i+1})\big](L_i)^{-1}
-\mathbf{a}_i\in span \{ \mathbf{a}_i,\mathbf{a}_{i+2},\dots,\mathbf{a}_{n+1}\}\ .
\end{align*}
Thus, 
\begin{align*}
\overline{\mathbf{a}_i}-\mathbf{a}_i
& =
2(\mathbf{a}_{i+1}-\mathbf{a}_i)
+
(-1)^{n-i+1}
2\big[L_i\circ (\mathbf{a}_i-\mathbf{a}_{i+1})\big](L_i)^{-1}\\
& =
2(\mathbf{a}_{i+1}-\mathbf{a}_i)
+
2\big[(\mathbf{a}_i-\mathbf{a}_{i+1})\circ L_i\big](L_i)^{-1}\\
& =
2(\mathbf{a}_{i+1}-\mathbf{a}_i)L_i(L_i)^{-1}
-
2\big[(\mathbf{a}_{i+1}-\mathbf{a}_i)\circ L_i\big](L_i)^{-1}\\
& =
2\Big\{
(\mathbf{a}_{i+1}-\mathbf{a}_i)L_i
-
\big[(\mathbf{a}_{i+1}-\mathbf{a}_i)\circ L_i\big]
\Big\}(L_i)^{-1}\\
& =
2
\big[(\mathbf{a}_{i+1}-\mathbf{a}_i)\Wedge L_i\big](L_i)^{-1}\ ,
\end{align*}
and you can verify that $|\overline{\mathbf{a}_i}-\mathbf{a}_{i+1}|=|\mathbf{a}_i-\mathbf{a}_{i+1}|$.
Let us proof a key property in the construction of the mirrored points.
\begin{Lem}\label{lem:telescopic property}
By keeping the notations and hypothesis of this section, we have that
\[
(\mathbf{a}_{i+1}-\mathbf{a}_i)\Wedge L_i
=
L_{i-1}
\ ,
\]
for each $i=2,\dots,n-1$.
\end{Lem}
{\bf Proof of Lemma~\ref{lem:telescopic property}}
\begin{align*}
& (\mathbf{a}_{i+1}-\mathbf{a}_i)\Wedge L_i =\\
= &
(\mathbf{a}_{i+1}-\mathbf{a}_i)
\wedge
(\mathbf{a}_{i+2}-\mathbf{a}_{i+1})\wedge\cdots\wedge(\mathbf{a}_{n+1}-\mathbf{a}_{i+1})\\
= &
(\mathbf{a}_{i+1}-\mathbf{a}_i)
\wedge
(\mathbf{a}_{i+2}-\mathbf{a}_i+\mathbf{a}_i-\mathbf{a}_{i+1})\wedge\cdots\wedge(\mathbf{a}_{n+1}-\mathbf{a}_i+\mathbf{a}_i-\mathbf{a}_{i+1})\\
= &
(\mathbf{a}_{i+1}-\mathbf{a}_i)
\wedge
(\mathbf{a}_{i+2}-\mathbf{a}_i)
\wedge\cdots\wedge
(\mathbf{a}_{n+1}-\mathbf{a}_i)=L_{i-1}\ . \ \square
\end{align*}
\begin{Rem}\label{rem:ortho-diagonals k}
By using Remark~\ref{rem:rejection is orthogonal k}, you can also verify that each vector $\overline{\mathbf{a}_i}-\mathbf{a}_i$ is orthogonal to all $\mathbf{a}_{i+2}-\mathbf{a}_{i+1}$,\dots,$\mathbf{a}_{n+1}-\mathbf{a}_{i+1}$. 
\end{Rem}
At the end of that process of mirroring the $n-1$ points $\mathbf{a}_1$,\dots, $\mathbf{a}_{n-1}$ we have that vectors
\begin{center}
$\overline{\mathbf{a}_1}-\mathbf{a}_1$,\dots , $\overline{\mathbf{a}_{n-1}}-\mathbf{a}_{n-1}$, and $\mathbf{a}_{n+1}-\mathbf{a}_n$
\end{center}
are mutually orthogonal. Moreover the following key property holds
\begin{Prop}\label{prop:hyperoctahedron}
By keeping the notations and hypothesis of this section, we have that
\begin{equation*}
\Big[
(\overline{\mathbf{a}_1}-\mathbf{a}_1)
\wedge
\cdots
\wedge
(\overline{\mathbf{a}_{n-1}}-\mathbf{a}_{n-1})
\Big]
\Wedge
(\mathbf{a}_{n+1}-\mathbf{a}_{n})
=
2^{n-1}
\Delta_{(\mathbf{a}_1,\dots,\mathbf{a}_{n+1})}
\end{equation*}
\end{Prop}
{\bf Proof of Proposition~\ref{prop:hyperoctahedron}}\\
By Lemma~\ref{lem:telescopic property}, we can write 
$
\overline{\mathbf{a}_i}-\mathbf{a}_i
=
2
L_{i-1}(L_i)^{-1}
$. So, by the foregoing orthogonality property we can arrive at a telescopic product
\begin{align*}
& 
\big[(\overline{\mathbf{a}_1}-\mathbf{a}_1)
\wedge
(\overline{\mathbf{a}_2}-\mathbf{a}_2)
\wedge
(\overline{\mathbf{a}_3}-\mathbf{a}_3)
\wedge\cdots\wedge
(\overline{\mathbf{a}_{n-1}}-\mathbf{a}_{n-1})
\big]
\Wedge
(\mathbf{a}_{n+1}-\mathbf{a}_{n})
=\\
= & 
(\overline{\mathbf{a}_1}-\mathbf{a}_1)
\wedge
(\overline{\mathbf{a}_2}-\mathbf{a}_2)
\wedge
(\overline{\mathbf{a}_3}-\mathbf{a}_3)
\wedge\cdots\wedge
(\overline{\mathbf{a}_{n-1}}-\mathbf{a}_{n-1})
\wedge
(\mathbf{a}_{n+1}-\mathbf{a}_{n})
\\
= & 
(\overline{\mathbf{a}_2}-\mathbf{a}_2)
(\overline{\mathbf{a}_3}-\mathbf{a}_3)
\cdots
(\overline{\mathbf{a}_{n-1}}-\mathbf{a}_{n-1})
(\mathbf{a}_{n+1}-\mathbf{a}_{n})=\\
= &
2^{n-1}
\big[(\mathbf{a}_{2}-\mathbf{a}_1)\Wedge L_1\big](L_1)^{-1}
L_{1}(L_2)^{-1}
L_{2}(L_3)^{-1}
\cdots
L_{n-2}(L_{n-1})^{-1}
(\mathbf{a}_{n+1}-\mathbf{a}_{n})
\\
= &
2^{n-1}
\big[(\mathbf{a}_{2}-\mathbf{a}_1)\Wedge L_1\big]
(L_{n-1})^{-1}
(\mathbf{a}_{n+1}-\mathbf{a}_{n})\\
= &
2^{n-1}
\big[(\mathbf{a}_{2}-\mathbf{a}_1)\Wedge L_1\big]
(\mathbf{a}_{n+1}-\mathbf{a}_n)^{-1}
(\mathbf{a}_{n+1}-\mathbf{a}_{n})\\
= &
2^{n-1}
(\mathbf{a}_{2}-\mathbf{a}_1)\Wedge L_1\\
= &
2^{n-1}
(\mathbf{a}_{2}-\mathbf{a}_1)
\wedge
(\mathbf{a}_{3}-\mathbf{a}_{2})\wedge\cdots\wedge(\mathbf{a}_{n+1}-\mathbf{a}_{2})\\
= &
2^{n-1}
(\mathbf{a}_{2}-\mathbf{a}_1)
\wedge
(\mathbf{a}_{3}-\mathbf{a}_{1})\wedge\cdots\wedge(\mathbf{a}_{n+1}-\mathbf{a}_{1})
=
2^{n-1}
\Delta_{(\mathbf{a}_1,\dots ,\mathbf{a}_{n+1})}\ .\ \square
\end{align*}

\subsection{Convergence of the mean secant hyperplane to the tangent hyperplane II}
Let us consider the $n+1$ points $\mathbf{a}_1$,\dots,$\mathbf{a}_{n+1}$ in $\mathbb{E}_n$ as in the foregoing paragraphs.
\begin{Def}\label{def:multi diff pseudo-vector n-1}
Let us define the following blades
\arraycolsep=1.pt\def\arraystretch{0.5}
\begin{align*}
&
\Delta
=
\Delta_{(\mathbf{a}_1,\dots,\mathbf{a}_{n+1})}
=
(\mathbf{a}_2-\mathbf{a}_1)
\wedge
\cdots
\wedge
(\mathbf{a}_{n+1}-\mathbf{a}_{1})
=
\bigwedge_{i=2}^{n+1}
(\mathbf{a}_i-\mathbf{a}_1)\in \mathbb{G}_{n\choose n}\simeq \mathbb{R}I_n\\
&
\overline{\Delta}
=
\overline{\Delta}_{(\mathbf{a}_1,\dots,\mathbf{a}_{n+1})}
=
(\overline{\mathbf{a}_1}-\mathbf{a}_1)
\wedge
\cdots
\wedge
(\overline{\mathbf{a}_{n-1}}-\mathbf{a}_{n-1})
=
\bigwedge_{i=1}^{n-1}
(\overline{\mathbf{a}_i}-\mathbf{a}_i)\in \mathbb{G}_{n\choose n-1}\\
&
\overline{\Delta^i}
=
\overline{\Delta^i}_{(\mathbf{a}_1,\dots,\mathbf{a}_{n+1})}
=
\left(
\bigwedge_{\scriptstyle
\begin{array}{c}
\scriptstyle j=1	\\
\scriptstyle j\ne i
\end{array}
}^{n-1}
(\overline{\mathbf{a}_j}-\mathbf{a}_j)
\right)
\Wedge
(\mathbf{a}_{n+1}-\mathbf{a}_{n})\in \mathbb{G}_{n\choose n-1}\ ,
\end{align*}
where $i=1,\dots,n-1$.\\
Then, we can define the ``mean multi-difference $\mathbb{G}_{n\choose n-1}$-pseudo-vector'',
\begin{align*}
& 
\overline{\Delta}f_{(\mathbf{a}_1,\dots,\mathbf{a}_{n+1})}
=\\
= &
\frac{1}{2^{n-1}}
\left\{
\kern-3pt
\left(
\sum_{i=1}^{n-1}
(-1)^{i+1}
\big[
f(\overline{\mathbf{a}_i})
-
f(\mathbf{a}_i)
\big]
\overline{\Delta^i}_{(\mathbf{a}_1,\dots,\mathbf{a}_{n+1})}
\kern-3pt
\right)
\kern-3pt +
(-1)^{n+1}
\big[
f(\mathbf{a}_{n+1})
-
f(\mathbf{a}_n)
\big]
\overline{\Delta}_{(\mathbf{a}_1,\dots,\mathbf{a}_{n+1})}
\kern-3pt
\right\}
\end{align*}
\end{Def}
\begin{Rem*}
Contrary to the pseudo-scalar $\Delta$, the foregoing blades $\overline{\Delta}$ and $\overline{\Delta^i}$ can always be written as geometric products
\arraycolsep=1.pt\def\arraystretch{0.5}
\begin{align*}
&
\overline{\Delta}
=
\overline{\Delta}_{(\mathbf{a}_1,\dots,\mathbf{a}_{n+1})}
=
(\overline{\mathbf{a}_1}-\mathbf{a}_1)
\cdots
(\overline{\mathbf{a}_{n-1}}-\mathbf{a}_{n-1})
=
\prod_{i=1}^{n-1}
(\overline{\mathbf{a}_i}-\mathbf{a}_i)\\
&
\overline{\Delta^i}
=
\overline{\Delta^i}_{(\mathbf{a}_1,\dots,\mathbf{a}_{n+1})}
=
\left(
\prod_{
\begin{array}{c}
\scriptstyle j=1	\\
\scriptstyle j\ne i
\end{array}
}^{n-1}
(\overline{\mathbf{a}_j}-\mathbf{a}_j)
\right)
(\mathbf{a}_{n+1}-\mathbf{a}_{n})\ ,
\end{align*}
(where $i=1,\dots,n-1$), because vectors $\overline{\mathbf{a}_1}-\mathbf{a}_1$, \dots, $\overline{\mathbf{a}_{n-1}}-\mathbf{a}_{n-1}$ and $\mathbf{a}_{n+1}-\mathbf{a}_n$ are mutually orthogonal. 
These are the key properties that allow us to prove the next general Theorem. 
\end{Rem*}
\begin{Thm}\label{thm:n}
If the function $f:\Omega\subseteq\mathbb{E}_n\to\mathbb{R}$ is strongly differentiable at $\mathbf{x}_0$ (a point internal to $\Omega$), then
\[
\lim_{\begin{array}{c}
	\scriptstyle (\mathbf{a}_1,\dots,\mathbf{a}_{n+1})
	\to (\mathbf{x_0},\dots,\mathbf{x_0})\in\mathbb{E}_n^{n+1}\\
	\scriptstyle \Delta_{(\mathbf{a}_1,\dots ,\mathbf{a}_{n+1})}\ne 0
\end{array}
}
\Big(\overline{\Delta}f_{(\mathbf{a}_1,\dots,\mathbf{a}_{n+1})}\Big)
\big(\Delta_{(\mathbf{a}_1,\dots,\mathbf{a}_{n+1})}\big)^{-1}
=\nabla f(\mathbf{x}_0)\ .
\]
\end{Thm}
\begin{Rem*}
The foregoing result state that the ``mean secant hyperplane''
\[
z= f(\mathbf{a}_1)\ +\ \overline{\mathbf{r}}_{f_{(\mathbf{a}_1,\dots,\mathbf{a}_{n+1})}} \cdot (\mathbf{x}-\mathbf{a}_1)\ ,
\]
where 
$\overline{\mathbf{r}}_{f_{(\mathbf{a}_1,\dots,\mathbf{a}_{n+1})}}
=
\Big(\overline{\Delta}f_{(\mathbf{a}_1,\dots,\mathbf{a}_{n+1})}\Big)
\big(\Delta_{(\mathbf{a}_1,\dots,\mathbf{a}_{n+1})}\big)^{-1}
$, always converges to the tangent hyperplane  
\[
z= f(\mathbf{x}_0)\ +\ \nabla f(\mathbf{x}_0) \cdot (\mathbf{x}-\mathbf{x}_0)\ ,
\]
as the non degenerate simplex having vertices $\mathbf{a}_1$, \dots, $\mathbf{a}_{n+1}$ contracts to the point~$\mathbf{x}_0$.
\end{Rem*}
\begin{Prop}\label{prop:mean multi-difference quotient n}
\[
\overline{\mathbf{r}}_{f_{(\mathbf{a}_1,\dots,\mathbf{a}_{n+1})}}
= 
\left(
\sum_{i=1}^{n-1}
\frac{f(\overline{\mathbf{a}_i})
-
f(\mathbf{a}_i)
}{\overline{\mathbf{a}_i}-\mathbf{a}_i}
\right)
+
\frac{f(\mathbf{a}_{n+1})
-
f(\mathbf{a}_n)}{\mathbf{a}_{n+1}-\mathbf{a}_n}\ .
\]
\end{Prop}
{\bf Proof of Proposition~\ref{prop:mean multi-difference quotient n}.}\\
Let us recall, from Remark~\ref{rem:ortho-diagonals k}, that vectors
$\overline{\mathbf{a}_1}-\mathbf{a}_1$,\dots, $\overline{\mathbf{a}_{n-1}}-\mathbf{a}_{n-1}$, $\mathbf{a}_{n+1}-\mathbf{a}_n$ are mutually orthogonal. So, by Proposition~\ref{prop:hyperoctahedron}, we can write
\begin{align*}
&
\overline{\Delta}_{(\mathbf{a}_1,\dots,\mathbf{a}_{n+1})}
\Wedge
(\mathbf{a}_{n+1}-\mathbf{a}_n)
= 
2^{n-1}
\Delta_{(\mathbf{a}_1,\dots,\mathbf{a}_{n+1})}\\
&
\Delta_{(\mathbf{a}_1,\dots,\mathbf{a}_{n+1})}
= 
\frac{1}{2^{n-1}}
\overline{\Delta}_{(\mathbf{a}_1,\dots,\mathbf{a}_{n+1})}
\Wedge
(\mathbf{a}_{n+1}-\mathbf{a}_n) \\
& \phantom{\hspace{55pt}} 
=
\frac{1}{2^{n-1}}
(\overline{\mathbf{a}_1}-\mathbf{a}_1)
\cdots
(\overline{\mathbf{a}_{n-1}}-\mathbf{a}_{n-1})
(\mathbf{a}_{n+1}-\mathbf{a}_n)\\
&
\big(\Delta_{(\mathbf{a}_1,\dots,\mathbf{a}_{n+1})}\big)^{-1}
= 
2^{n-1}
(\mathbf{a}_{n+1}-\mathbf{a}_n)^{-1}
(\overline{\mathbf{a}_{n-1}}-\mathbf{a}_{n-1})^{-1}
\cdots
(\overline{\mathbf{a}_1}-\mathbf{a}_1)^{-1}\ .
\end{align*}
Then, 
\begin{align*}
&
\overline{\mathbf{r}}_{f_{(\mathbf{a}_1,\dots,\mathbf{a}_{n+1})}}
= 
\Big(\overline{\Delta}f_{(\mathbf{a}_1,\dots,\mathbf{a}_{n+1})}\Big)
\big(\Delta_{(\mathbf{a}_1,\dots,\mathbf{a}_{n+1})}\big)^{-1}
=\\
= &
\left(
\sum_{i=1}^{n-1}
(-1)^{i+1}
\big[
f(\overline{\mathbf{a}_i})
-
f(\mathbf{a}_i)
\big]
\overline{\Delta^i}_{(\mathbf{a}_1,\dots,\mathbf{a}_{n+1})}
\kern-3pt
\right)\big(\Delta_{(\mathbf{a}_1,\dots,\mathbf{a}_{n+1})}\big)^{-1}+\\
& +
(-1)^{n+1}
\big[
f(\mathbf{a}_{n+1})
-
f(\mathbf{a}_n)
\big]
\overline{\Delta}_{(\mathbf{a}_1,\dots,\mathbf{a}_{n+1})}
\big(\Delta_{(\mathbf{a}_1,\dots,\mathbf{a}_{n+1})}\big)^{-1}\\
= &
\left(
\sum_{i=1}^{n-1}
\big[
f(\overline{\mathbf{a}_i})
-
f(\mathbf{a}_i)
\big]
(\overline{\mathbf{a}_i}-\mathbf{a}_i)^{-1}
\right)
+
\big[
f(\mathbf{a}_{n+1})
-
f(\mathbf{a}_n)
\big]
(\mathbf{a}_{n+1}-\mathbf{a}_n)^{-1}\ ,
\end{align*}
as you can verify that 
$
\overline{\Delta}\Delta^{-1}
=
(-1)^{n+1}
2^{n-1}
(\mathbf{a}_{n+1}-\mathbf{a}_{n})^{-1}
$ 
and, 
\begin{center}
$
\overline{\Delta^i}\Delta^{-1}
=
(-1)^{i+1}
2^{n-1}
(\overline{\mathbf{a}_{i}}-\mathbf{a}_{i})^{-1}\ .
$, for $i=1\dots,n-1$.~$\square$
\end{center}
{\bf Proof of Theorem~\ref{thm:n}.}\\

Let us recall that $\Delta_{(\mathbf{a}_1,\dots,\mathbf{a}_{n+1})}$ is a pseudo-scalar; so, by equation~(\ref{eq:circ}), we have that
\begin{align*}
& \nabla f(\mathbf{x}_0)\Delta_{(\mathbf{a}_1,\dots,\mathbf{a}_{n+1})}
=
\nabla f(\mathbf{x}_0)\circ\Delta_{(\mathbf{a}_1,\dots,\mathbf{a}_{n+1})}\\
= &
\frac{1}{2^{n-1}}
\left[
\sum_{i=1}^{n-1} (-1)^{i+1}
\big(\nabla f(\mathbf{x}_0)\cdot 
(\overline{\mathbf{a}_i}-\mathbf{a}_i)\big) 
\overline{\Delta^i}
\ +\
(-1)^{n+1}
\big(\nabla f(\mathbf{x}_0)\cdot 
(\mathbf{a}_{n+1}-\mathbf{a}_n)\big) 
\overline{\Delta}
\right]\ ,
\end{align*}
and we can write
\begin{align*}
& 
\nabla f(\mathbf{x}_0)
=
\nabla f(\mathbf{x}_0)
\Delta_{(\mathbf{a}_1,\dots,\mathbf{a}_{n+1})}
\big(\Delta_{(\mathbf{a}_1,\dots,\mathbf{a}_{n+1})}\big)^{-1}
=\\
= & 
\big(
\nabla f(\mathbf{x}_0)\circ\Delta_{(\mathbf{a}_1,\dots,\mathbf{a}_{n+1})}\big)
\big(\Delta_{(\mathbf{a}_1,\dots,\mathbf{a}_{n+1})}\big)^{-1}\\
= &
\frac{1}{2^{n-1}}
\left[
\sum_{i=1}^{n-1} (-1)^{i+1}
\big(\nabla f(\mathbf{x}_0)\cdot 
(\overline{\mathbf{a}_i}-\mathbf{a}_i)\big) 
\overline{\Delta^i}
\ +\
(-1)^{n+1}
\big(\nabla f(\mathbf{x}_0)\cdot 
(\mathbf{a}_{n+1}-\mathbf{a}_n)\big) 
\overline{\Delta}
\right]\Delta^{-1}\\
= &
\left(
\sum_{i=1}^{n-1} 
\big(\nabla f(\mathbf{x}_0)\cdot 
(\overline{\mathbf{a}_i}-\mathbf{a}_i)\big) 
(\overline{\mathbf{a}_i}-\mathbf{a}_i)^{-1}
\right)
\ +\
\big[\nabla f(\mathbf{x}_0)\cdot 
(\mathbf{a}_{n+1}-\mathbf{a}_n)\big] 
(\mathbf{a}_{n+1}-\mathbf{a}_n)^{-1}\\
= &
\left(
\sum_{i=1}^{n-1} 
\frac{\nabla f(\mathbf{x}_0)\cdot 
(\overline{\mathbf{a}_i}-\mathbf{a}_i)}
{\overline{\mathbf{a}_i}-\mathbf{a}_i}
\right)
\ +\
\frac{\nabla f(\mathbf{x}_0)\cdot (\mathbf{a}_{n+1}-\mathbf{a}_n)}
{\mathbf{a}_{n+1}-\mathbf{a}_n}\ .
\end{align*}
So,
\begin{align*}
& 
\overline{\mathbf{r}}_{f_{(\mathbf{a}_1,\dots,\mathbf{a}_{n+1})}}
-
\nabla f(\mathbf{x}_0)
=\\
= &
\left(
\sum_{i=1}^{n-1} 
\frac{f(\overline{\mathbf{a}_i})
-
f(\mathbf{a}_i)
-
\nabla f(\mathbf{x}_0)\cdot 
(\overline{\mathbf{a}_i}-\mathbf{a}_i)
}
{\overline{\mathbf{a}_{i}}-\mathbf{a}_{i}} 
\right)
+
\frac{f(\mathbf{a}_{n+1})
-
f(\mathbf{a}_n)
-
\nabla f(\mathbf{x}_0)\cdot 
(\mathbf{a}_{n+1}-\mathbf{a}_n)}
{\mathbf{a}_{n+1}-\mathbf{a}_{n}}\ ,
\end{align*} 
and
\begin{align*}
\Big|
\overline{\mathbf{r}}_{f_{(\mathbf{a}_1,\dots,\mathbf{a}_{n+1})}}
-
\nabla f(\mathbf{x}_0)
\Big|
\le &
\sum_{i=1}^{n-1} 
\frac{
\big|
f(\overline{\mathbf{a}_i})
-
f(\mathbf{a}_i)
-
\nabla f(\mathbf{x}_0)\cdot 
(\overline{\mathbf{a}_i}-\mathbf{a}_i)
\big|}{\big|\overline{\mathbf{a}_{i}}-\mathbf{a}_{i}\big|}+\\
& +
\frac{
\big|
f(\mathbf{a}_{n+1})
-
f(\mathbf{a}_n)
-
\nabla f(\mathbf{x}_0)\cdot 
(\mathbf{a}_{n+1}-\mathbf{a}_n)
\big|
}{\big|\mathbf{a}_{n+1}-\mathbf{a}_{n}\big|}\ .
\end{align*}
As $f$ is strongly differentiable at $\mathbf{x}_0$, we know that, given $\epsilon>0$ there exists $\delta_\epsilon>0$ such that if  $|\mathbf{u}-\mathbf{x}_0|<\delta_\epsilon$ and $|\mathbf{v}-\mathbf{x}_0|<\delta_\epsilon$, then
\begin{center}
$
|f(\mathbf{u})-f(\mathbf{v})-\nabla f(\mathbf{x}_0)\cdot(\mathbf{u}-\mathbf{v})|<\epsilon |\mathbf{u}-\mathbf{v}|\ .
$
\end{center}
So, if we choose $\mathbf{a}_1$, \dots, $\mathbf{a}_{n+1}$ non collinear, and such that
\begin{center}
 $|\mathbf{a}_i-\mathbf{x}_0|
	< \frac{1}{3}\left(\delta_{\frac{\epsilon}{n}}\right)$, for $i=1,\dots,n+1$,
\end{center}
then 
\[
\big|
\overline{\mathbf{r}}_{f_{(\mathbf{a}_1,\dots,\mathbf{a}_{n+1})}}
-
\nabla f(\mathbf{x}_0)
\big|<\epsilon\ ,
\]
as, for $i=1,\dots,n-1$, we have
\begin{itemize}
	\item $|\mathbf{a}_i-\mathbf{a}_{i+1}|\le |\mathbf{a}_i-\mathbf{x}_0|+ |\mathbf{a}_{i+1}-\mathbf{x}_0|< \frac{2}{3}\left(\delta_{\frac{\epsilon}{n}}\right)$,
	\item $|\overline{\mathbf{a}_i}-\mathbf{x}_0|\le|\overline{\mathbf{a}_i}-\mathbf{a}_{i+1}|+|\mathbf{a}_{i+1}-\mathbf{x}_0|=|\mathbf{a}_i-\mathbf{a}_{i+1}|+|\mathbf{a}_{i+1}-\mathbf{x}_0| < \delta_{\frac{\epsilon}{3}}$.~$\square$
\end{itemize}

\section{Conclusions and perspectives}\label{sec:conclusions}

Here, we sketch some possible applications of the foregoing results in pure and applied mathematics.

\subsection{Geometric Calculus}
At page~45 of~\cite{HestenesSob} it is written:

\begin{center}
\textit{Leibniz's notation $dF/d\tau$ or $\partial F/\partial\tau$ emphasizes the definition of derivative as the limit of a difference quotient. It will be seen that differentiation by a general multivector cannot be defined by a difference quotient, so Leibniz's notation is appropriate only for scalar variables.} 
\end{center}

\noindent
In this work we have shown that, using as product the Clifford geometric product, as numerator our mean multi-difference pseudo-vector~$\overline{\Delta}f$, and as denominator our pseudo-scalar~$\Delta$, Leibniz's quotient notation  
\[
\overline{\Delta}f/\Delta
=
\Big(\overline{\Delta}f_{(\mathbf{a}_1,\dots,\mathbf{a}_{n+1})}\Big)
\big(\Delta_{(\mathbf{a}_1,\dots,\mathbf{a}_{n+1})}\big)^{-1}\ ,
\]
remains appropriate to approximate  (or define) the gradient of multi-variable scalar functions (provided they are strongly differentiable). 

\subsection{Fundamental Theorem of Calculus for a generalized Riemann integral}
The foregoing observation can suggest the possibility of extending to higher dimensions the elegant definition\footnote{See \cite{Mawhin}.} of non absolute integration given by Kurzweil and Henstock by using our results. In particular, we guess the possibility of defining a multi-vector valued generalized Riemann integral on orientable hyper-volumes and hyper-surfaces such that a Fundamental theorem of Calculus of the following form
\begin{equation}
\label{eq:KH integrals}
\underset{
\begin{array}{c}
	\scriptstyle\text{oriented}\\
	\scriptstyle\text{hyper-volume}\\
	\Omega
\end{array}
}{\iint}
\underbrace{
\underbrace{\nabla f(\mathbf{x})}_{\textrm{vector}}\underbrace{d\mathbf{x}}_{\textrm{pseudo-scalar}}}_{\textrm{pseudo-vector}}
=
\underset{
\begin{array}{c}
	\scriptstyle\text{oriented}\\
	\scriptstyle\text{rectifiable}\\
	\scriptstyle\text{hyper-surface}\\
	\partial\Omega
\end{array}
}{\int}
\underbrace{
\underbrace{f(\mathbf{x})}_{\textrm{scalar}}\underbrace{d\mathbf{x}}_{\textrm{pseudo-vector}}}_{\textrm{pseudo-vector}}
\end{equation}
may hold, once $f:\Omega\subset\mathbb{E}_n\to\mathbb{R}$ is everywhere strongly differentiable. If this would be the case, we could aim to obtain a notion of Kurzweil-Henstock $n$-dimensional integrability which would be invariant by rotation, while, to my knowledge, the known extensions of Kurzweil-Henstock integral to $n$-dimensions are not rotation invariant (that is, it can happen that rotations of integrable functions are no more integrable). Moreover, once defined the generalized integrals satisfying~(\ref{eq:KH integrals}), one could obtain a subsequent non absolute change of variable formula for multiple integrals.

\subsection{Finite Element Method}
Delaunay triangulations maximize the minimum angle of all angles of the triangles in a triangulation of a plane point cloud. This extremal property allows in Finite Element Method analysis to prevent divergence phenomena due to the Schwarz paradox. The algorithms we have used do not suffers such divergence phenomena; so, they could be used in Finite Element Method on arbitrary triangulations (non necessarily of Delaunay type).

One can profit further of the above possibility to neglect Delaunay triangulations. As a matter of fact,  for $n$-dimensional point clouds (with $n\ge 3$) it is more and more difficult to generate  a Delaunay ``tretrahedralization''. On the contrary, we claim that our convergent algorithms (independent from the maximal Delaunay property) can be applied to efficiently discretize differential operators on any $n$-dimensional smooth manifold.

\subsection{Teaching}
The classical global Schwarz paradox deals with the area of surfaces, and its presentation is often omitted from Advanced Calculus curricula. In this work we have shown that the local form of the Schwarz paradox  involves the very definition of gradient, and could be presented at the beginning of a multi-variable differential calculus course. In this sense, it provides a further motivation to introduce the Clifford Geometric Algebra as a unifying language for Geometry and Analysis.

\end{document}